\documentclass{amsart}
\usepackage{amsmath}
\usepackage{graphicx}
\usepackage{latexsym}
\usepackage{color}
\usepackage{amscd}
\usepackage{enumerate}
\usepackage{hyperref}
\usepackage{subfigure}
\usepackage{soul}
\usepackage{comment}
\usepackage{epsfig}
\usepackage{epstopdf}
\usepackage{afterpage}

\newtheorem{thm}{Theorem}

\newtheorem{theorem}[thm]{Theorem}

\newtheorem{question}[thm]{Question}

\theoremstyle{definition}

\newtheorem*{definition*}{Definition}

\newtheorem{rmk}[thm]{Remark}

\newtheorem{remark}[thm]{Remark}

\newcommand{\rmI}{I}
\newcommand{\rmII}{I\hspace{-.1em}I}
\newcommand{\Mod}{\mathrm{Mod}}

\newcommand{\Z}{\mathbb{Z}}

\allowdisplaybreaks[4]
\begin{document}

\title[Sections of the Matsumoto-Cadavid-Korkmaz Lefschetz fibration] 
{Sections of the Matsumoto-Cadavid-Korkmaz Lefschetz fibration}

\author[N. Hamada]{Noriyuki Hamada}
\address{Graduate School of Mathematical Sciences, The University of Tokyo, Komaba, Meguro-ku, Tokyo 153-8914, Japan}
\email{nhamada@ms.u-tokyo.ac.jp}

\begin{abstract}
We give a maximal set of disjoint $(-1)$-sections of the well-known Lefschetz fibration constructed by Matsumoto, Cadavid and Korkmaz.
In fact, we obtain several such sets for a fixed genus, which implies that the Matsumoto-Cadavid-Korkmaz Lefschetz fibration has more than one supporting minimal Lefschetz pencils. 
We also determine the diffeomorphism types of the obtained supporting minimal Lefschetz pencils.
\end{abstract}

\maketitle

\section{Introduction} 
In smooth $4$-dimensional topology, Lefschetz fibrations and pencils have been of great interest due to its close relationship to symplectic $4$-manifolds and its combinatorial description via mapping class groups.
As the notion of Lefschetz fibration came from that of Lefschetz pencil, blowing up at the base points of a given Lefschetz pencil naturally yields a Lefschetz fibration in which the exceptional spheres become disjoint sections of self-intersection $-1$ (called \textit{$(-1)$-sections}).
In such a situation, we say that the Lefschetz pencil \textit{supports} the resulting Lefschetz fibration.
Conversely, if a Lefschetz fibration has disjoint $(-1)$-sections it can be blown down to obtain a Lefschetz pencil that supports the original fibration.
In this way we can interchangeably think of the base points of a Lefschetz pencil as a set of disjoint $(-1)$-sections of a Lefschetz fibration.
Therefore, to investigate disjoint $(-1)$-sections of a given Lefschetz fibration, or equivalently to investigate its supporting Lefschetz pencils, has been a fundamental concern.
In particular, it has been extensively studied in the case of the Lefschetz fibration whose monodromy is the hyperelliptic relation~\cite{KorkmazOzbagci2008, Onaran2010, Tanaka2012}.

Regarding this motivation, there is a remarkably important example of Lefschetz fibration, which we would like to explore.
In~\cite{Matsumoto1996}, Matsumoto originally constructed a genus-$2$ Lefschetz fibration on $T^2 \times S^2 \# 4\overline{\mathbb{CP}}{}^{2}$ with $8$ critical points.
Then Cadavid~\cite{Cadavid1998} and Korkmaz~\cite{Korkmaz2001} independently generalized it to higher genera:
a genus-$g$ Lefschetz fibration on $\Sigma_{g/2} \times S^2 \# 4\overline{\mathbb{CP}}{}^{2}$ with  $2g+4$ critical points for even $g \geq 2$, 
or on $\Sigma_{(g-1)/2} \times S^2 \# 8\overline{\mathbb{CP}}{}^{2}$ with $2g+10$ critical points for odd $g \geq 1$, respectively.
We call this Lefschetz fibration \textit{the Matsumoto-Cadavid-Korkmaz Lefschetz fibration} (MCK for short) in this paper.
The MCK Lefschetz fibration has become one of the most basic examples in the theory of Lefschetz fibrations and played great roles, especially as a powerful source to construct new Lefschetz fibrations, surface bundles, Stein fillings, symplectic $4$-manifolds, and so on, with various interesting features~\cite{OzbagciStipsicz2000, Korkmaz2001, Korkmaz2009, Stipsicz2002, OzbagciStipsicz2004, Gurtas_preprint2005, Yun2006, ParkYun2009, ParkYun2011, Baykur2012, AkhmedovOzbagci2014, AkhmedovSaglam2015, AkhmedovSakalli2016, Baykur2016, BaykurHayano2016, HamadaKobayashiMonden_preprint, Kobayashi2016, KobayashiMonden2016, OkudaTakamura_preprint}. 
The MCK Lefschetz fibration itself has several remarkable features such as having quite small number of critical points (the smallest among the known examples for $g\geq4$), large $b_1$ (the largest among the known for even $g$ \footnote{For odd $g$ it had been also the largest until Baykur~\cite{Baykur_preprint} recently found a Lefschetz fibration with $b_1$ one larger than that of the MCK Lefschetz fibration.}), high symmetricity of the vanishing cycles, in particular, it is hyperelliptic when $g$ is even.
Besides, the MCK Lefschetz fibration can be also viewed as a generalization of the well-known elliptic Lefschetz fibration $E(1)= \mathbb{CP}{}^{2} \# 9\overline{\mathbb{CP}}{}^{2} \rightarrow S^2$; the genus-$1$ MCK Lefschetz fibration is isomorphic to this elliptic fibration.

Although a set of \textit{two} disjoint $(-1)$-sections of the MCK Lefschetz fibration has been already known (see below), the problem asking the maximal number of disjoint $(-1)$-sections it can admit has been still unsolved.
Since the underlying $4$-manifold has four or eight exceptional spheres (depending on the parity of genus), we might expect that the fibration can also admit four or eight  disjoint $(-1)$-sections, respectively.
The main aim of this paper is to show that this expectation is quite right,
giving such a maximal set of disjoint $(-1)$-sections for arbitrary genus by explicitly constructing a monodromy factorization that locates the required sections.

In fact, we do not give \textit{only one} such maximal set, but also give \textit{several}.
This is because such a maximal set of $(-1)$-sections provides a \textit{minimal} Lefschetz pencil and we are also interested in what kind of supporting minimal Lefschetz pencils for the MCK Lefschetz fibration there exist.
We will study the pencil structures of the obtained supporting minimal pencils and see that some of them are not mutually isomorphic.
Namely, the MCK Lefschetz fibration has more than one supporting minimal Lefschetz pencils (Theorem~\ref{thm:multiplepencils}).
To the best of the author's knowledge, this is the first example of a Lefschetz fibration having multiple supporting minimal Lefschetz pencils.
We also observe that the minimal Lefschetz pencils supporting the MCK Lefschetz fibrations exhaust all the diffeomorphism types of minimal symplectic $4$-manifolds with symplectic Kodaira dimension $-\infty$ (Remark~\ref{rmk:Kodaira-infty}).

For the convenience of the reader we briefly recall and summarize the correspondence between Lefschetz fibrations and relations among Dehn twists.
We refer to~\cite{GompfStipsicz1999} for details and the undefined terminology.
Throughout the paper, we assume a Lefschetz fibration to be smooth, relatively minimal, nontrivial, i.e. it has at least one critical point, and over the base space $S^2$.
Let $\Sigma_g^k$ denote a compact oriented surface of genus $g$ with $k$ boundary components and $\mathrm{Mod}(\Sigma_g^k)$ the mapping class group of $\Sigma_g^k$ whose elements (and isotopies used in the definition) are restricted to be identity on the boundary.
If $k=0$ it is dropped from the notation.
We adopt the functional notation for the product of mapping class groups: for two elements $\varphi_1, \varphi_2$ in $\mathrm{Mod}(\Sigma_g^k)$ the product $\varphi_2 \varphi_1$ means applying $\varphi_1$ first and then $\varphi_2$.
Via its monodromy representation, a genus-$g$ Lefschetz fibration with $n$ critical points gives rise to a \textit{positive factorization} (or \textit{monodromy factorization}) of the identity in $\mathrm{Mod}(\Sigma_g)$
\begin{equation*}
	t_{a_n} \cdots t_{a_2} t_{a_1}= 1
\end{equation*}
by nontrivial right-handed Dehn twists $\{ t_{a_i} \}$ (along its \textit{vanishing cycles}). 
The factorization is uniquely determined up to so-called \textit{Hurwitz equivalence}, which is the equivalence relation generated by \textit{Hurwitz moves} and \textit{simultaneous conjugations}.
Conversely, given such a positive factorization we can construct a Lefschetz fibration whose monodromy factorization is the given one.

If a genus-$g$ Lefschetz fibration with $t_{a_n} \cdots t_{a_2} t_{a_1} = 1$ has $k$ disjoint $(-1)$-sections $\{ s_j \}$, the monodromy factorization can be lifted to $\mathrm{Mod}(\Sigma_g^k)$ in the form
\begin{equation*}
	t_{\tilde{a}_n}  \cdots t_{\tilde{a}_2} t_{\tilde{a}_1} = t_{\delta_1} t_{\delta_2} \cdots t_{\delta_k},
\end{equation*}
where $\tilde{a}_i$ is a lift of $a_i$ and $\delta_j$ is the curve parallel to the $j$-the boundary components: that is a positive factorization of the \textit{boundary multi-twist} $t_{\delta_1} t_{\delta_2} \cdots t_{\delta_k}$.
To be precise, here we identify the reference fiber $F$ of the Lefschetz fibration with $\Sigma_g$ and $F \setminus \cup (\text{small open disk neighborhood of } F \cap s_j)$ with $\Sigma_g^k$.
Again, conversely, given a positive factorization $t_{\tilde{a}_n}  \cdots t_{\tilde{a}_2} t_{\tilde{a}_1} = t_{\delta_1} t_{\delta_2} \cdots t_{\delta_k}$ of the boundary multi-twist in $\mathrm{Mod}(\Sigma_g^k)$ such that $\tilde{a}_i$ descends to a homotopically non-trivial curve $a_i$ in $\Sigma_g$, we can construct a Lefschetz fibration with the desired $(-1)$-sections $\{ s_j \}$ (cf.~\cite{EKKOS2002}).
In terms of this description of $(-1)$-sections of Lefschetz fibrations---hence Lefschetz pencils---, the isomorphism class of a Lefschetz pencil is uniquely determined by the positive factorization of the boundary multi-twist up to the \textit{generalized Hurwitz equivalence} in the sense of Baykur-Hayano~\cite{BaykurHayano_preprint} (see also~\cite{Matsumoto1996}).

\begin{figure}[htbp]
	\centering
	\subfigure[For odd genus $g$. \label{F:MCKLF_odd}]
	{\includegraphics[height=90pt]{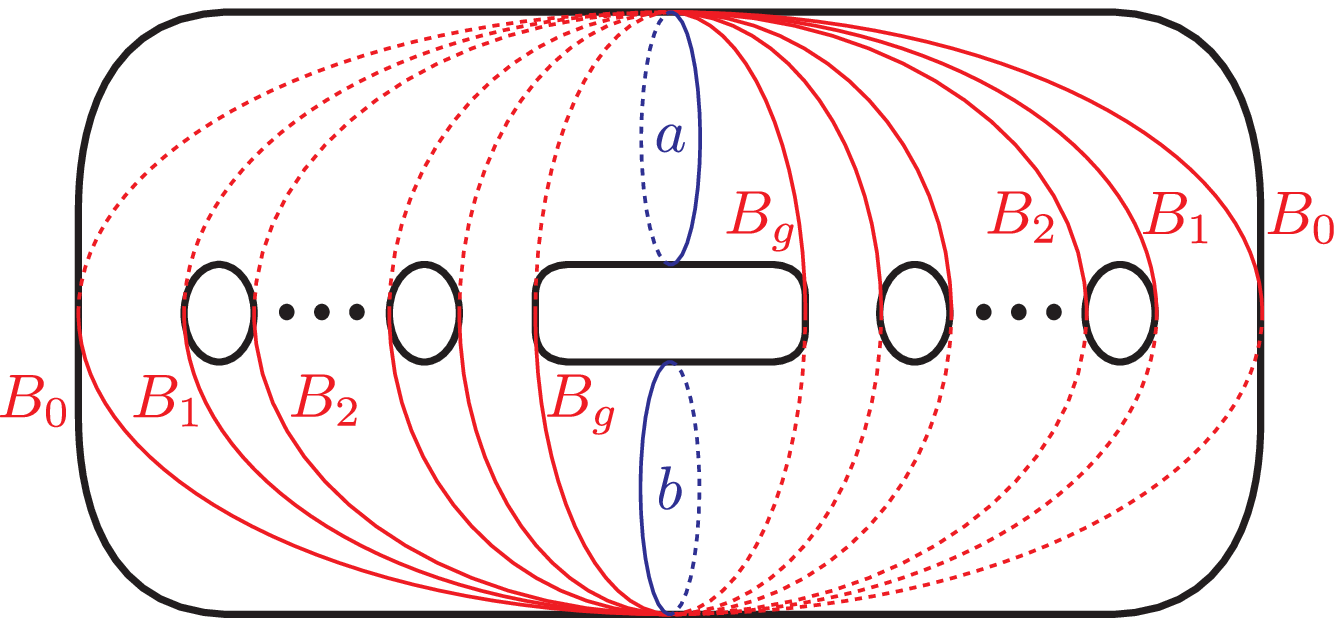}} 
	\hspace{0pt}	
	\subfigure[For even genus $g$. \label{F:MCKL_even}]
	{\includegraphics[height=90pt]{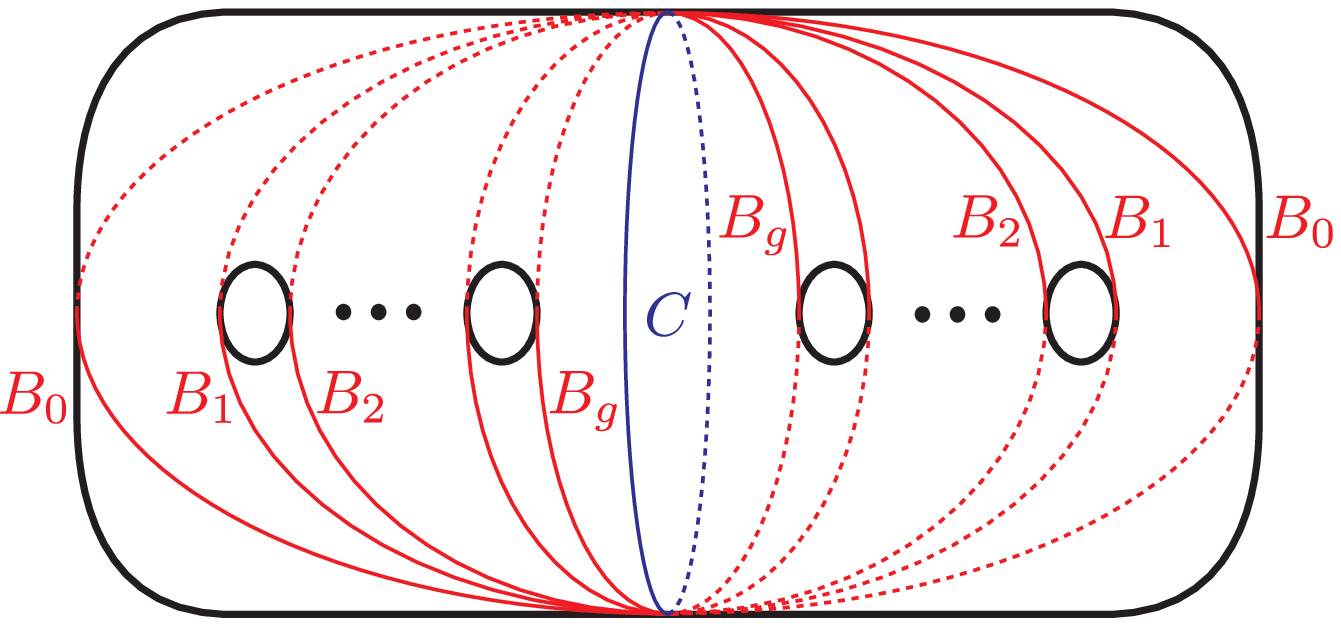}} 
	\caption{The vanishing cycles for the MCK Lefschetz fibration.} 	
	\label{F:MCKLF}
\end{figure}

Now we can introduce the MCK Lefschetz fibration.
With the curves in $\Sigma_g$ as depicted in Figure~\ref{F:MCKLF}, we set
\begin{equation}
	W = 
	\begin{cases}
		(t_{B_0} t_{B_1} \dots t_{B_g} t_{C})^2 & \text{for even $g$} \\
		(t_{B_0} t_{B_1} \dots t_{B_g} t_{a}^2 t_{b}^2)^2 & \text{for odd $g$}.
	\end{cases}
\end{equation}
Then the relation $W=1$ holds in $\Mod(\Sigma_g)$.
Thus we can obtain a Lefschetz fibration $f_W : X^4_W \rightarrow S^2$ with the monodromy factorization $W=1$; we call $f_W$ \textit{the Matsumoto-Cadavid-Korkmaz Lefschetz fibration}.
We also refer to the relation $W=1$ as \textit{the MCK relation}.
We remark that when $g=1$, one can directly check that the MCK relation is Hurwitz equivalent to the relation $(t_{\alpha} t_{\beta})^6=1$ where $\alpha$ is the meridian and $\beta$ is the longitude of the torus, that is, a monodromy factorization of the elliptic Lefschetz fibration $E(1) \rightarrow S^2$.
In this sense, we can regard the MCK Lefschetz fibration as a generalization of the elliptic Lefschetz fibration.
As we already mentioned, the diffeomorphism type of the total space is as follows:
\begin{equation*}
	X_W = 
	\begin{cases}
		\Sigma_{g/2} \times S^2 \# 4\overline{\mathbb{CP}}{}^{2} & \text{for even $g$} \\
		\Sigma_{(g-1)/2} \times S^2 \# 8\overline{\mathbb{CP}}{}^{2} & \text{for odd $g$}.
	\end{cases}
\end{equation*}
\begin{figure}[htbp]
	\centering
	\subfigure[For odd genus $g$. \label{F:MCKLF_odd_twosections}]
	{\includegraphics[height=90pt]{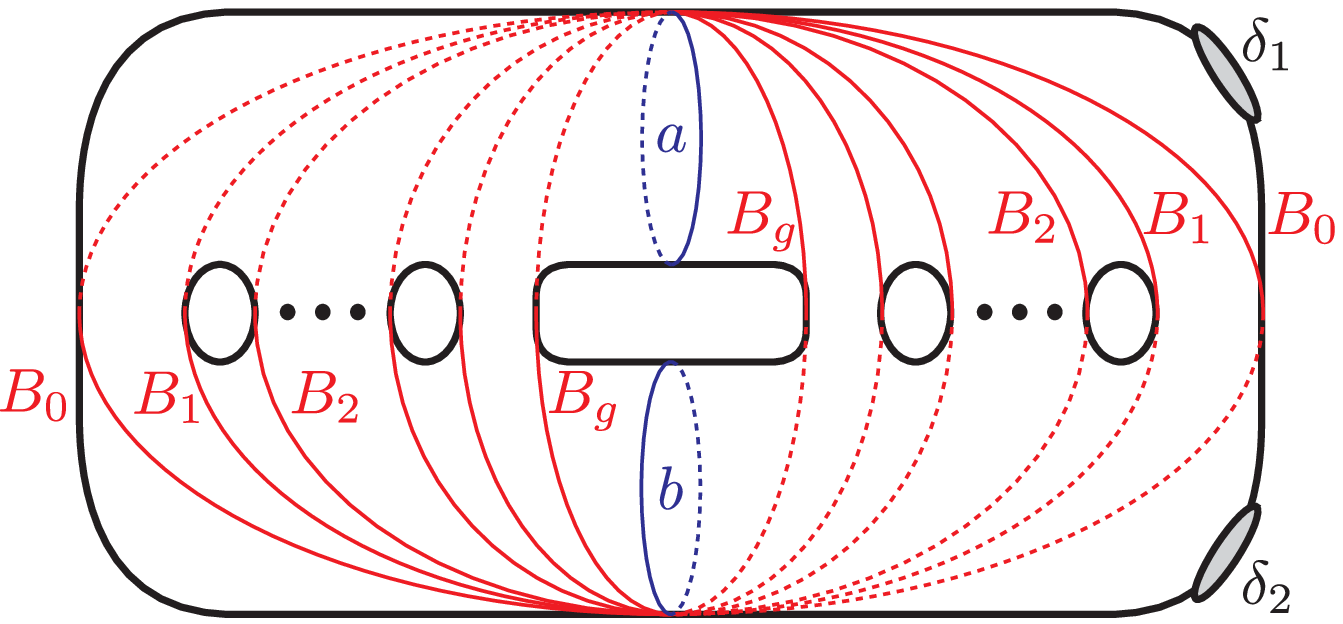}} 
	\hspace{0pt}	
	\subfigure[For even genus $g$. \label{F:MCKLF_even_twosections}]
	{\includegraphics[height=90pt]{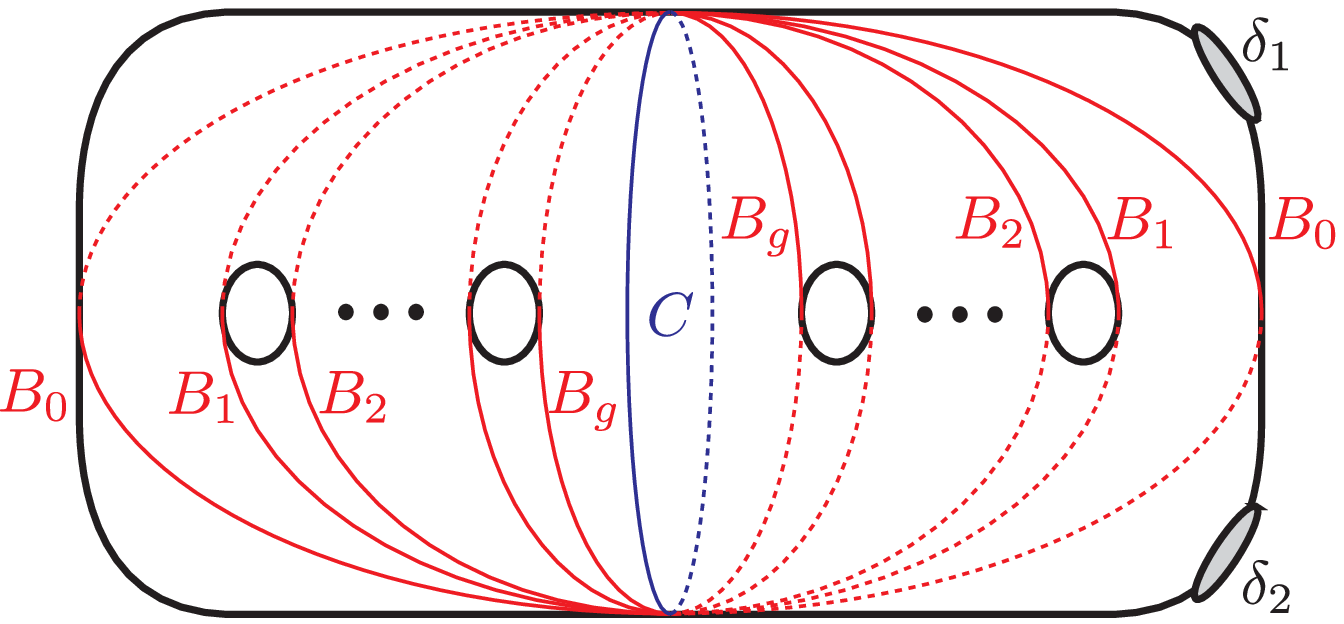}} 
	\caption{The curves for the known lift of the MCK relation.}
	\label{F:MCKLF_twosections} 	
\end{figure}

Moreover, from Korkmaz's construction via the Birman-Hilden double-covering, the relation $W=1$ can be naturally lifted to $\Mod(\Sigma_g^2)$:
\begin{equation} \label{eq:knownlifts}
	W= t_{\delta_1} t_{\delta_2},
\end{equation}
where the lifted curves are as shown in Figure~\ref{F:MCKLF_twosections}~\cite{Korkmaz2001, Korkmaz2009}
(Here we used the same symbols as in $\Sigma_g$, by abuse of notation).
It follows that the MCK Lefschetz fibration has at least two disjoint $(-1)$-sections.

In summary, the main arguments in the paper, which appear in Section~\ref{S:LiftsMCK}, consist of establishing positive factorizations of the boundary multi-twist that are lifts of the monodromy factorization $W=1$ of the MCK Lefschetz fibration
to $\mathrm{Mod}(\Sigma_g^k)$ with the largest possible $k$ (specifically, $k=4$ for even $g$ and $k=8$ for odd $g$).
We will see that some of the new lifts are also further lifts of the known one~\eqref{eq:knownlifts}.
By blowing down the resulting sections we obtain \textit{minimal} Lefschetz pencils.
In Section~\ref{S:sections-LPs}, we will determine the diffeomorphism types of the total spaces of those pencils.

\section{Basic relations}

Our method to construct positive factorizations of the boundary multi-twist is purely combinatorial: combining known relations to get a new relation.
We also utilize relations in the braid groups via the Birman-Hilden double-covering, as Korkmaz originally did, to make the argument clearer.
For convenience's sake, in this Section, we gather several known relations among half twists and Dehn twists that will be used later.
For a comprehensive reference, consult~\cite{FarbMargalit2012}. 

Let $B_n$ be the braid group on $n$ strands with the standard generators $\sigma_1, \sigma_2, \cdots, \sigma_{n-1}$ where $\sigma_i$ is the braid that has only one crossing at which the $i$-th strand passes in front of the (i+1)-st strand.
When we think of $B_n$ as the mapping class group of a closed disk with horizontally arranged $n$ marked points,
we can also describe the braid $\sigma_i$ as the right-handed half-twist $\tau_{\gamma_i}$ about the standard arc $\gamma_i$ connecting the $i$-th puncture and the $(i+1)$-st puncture in the disk.
We again use the functional notation for the product of braid groups as we will mainly see them as mapping class groups.

\subsection{Braid relation}
We will denote by $\tau_{\alpha}$ the right-handed half-twist about a simple proper arc $\alpha$ connecting two punctures in a punctured disk.
Let $\beta$ be another such arc, then we have the relation
\begin{equation*}
	\tau_{\alpha} \tau_{\beta} \tau_{\alpha}^{-1} = \tau_{\tau_{\alpha}(\beta)},
	\qquad \text{or equivalently,} \qquad
	\tau_{\alpha} \tau_{\beta} = \tau_{\tau_{\alpha}(\beta)} \tau_{\alpha},
\end{equation*}
which we call the \textit{braid relation}.
This relation implies the \textit{usual} braid relation: $\tau_{\alpha} \tau_{\beta} \tau_{\alpha} = \tau_{\beta} \tau_{\alpha} \tau_{\beta} $ when $\alpha$ and $\beta$ intersect exactly once at an end point.
It also implies that $\tau_{\alpha} \tau_{\beta} = \tau_{\beta} \tau_{\alpha} $ when $\alpha$ and $\beta$ are disjoint.

As for Dehn twists, for any simple closed curves $a$ and $b$ (on any surface), we have the relation
\begin{equation*}
	t_{a} t_{b} = t_{t_{a} (b)} t_{a},
\end{equation*}
which is a variation of the braid relation.
Similarly as the half-twists, we can induce that 
$t_{a} t_{b} t_{a} = t_{b} t_{a} t_{b} $ when $a$ and $b$ intersect transversely at one point and that $t_{a} t_{b} = t_{b} t_{a} $ when $a$ and $b$ are disjoint. 
In the context of monodromy factorizations, exchanging a subword $t_{a} t_{b}$ to $t_{t_{a} (b)} t_{a}$ or vice versa is called a \textit{Hurwitz move} or an \textit{elementary transformation}.

We will freely use braid relations (often without mentioning it) in the calculations.

\subsection{The chain relation}
In the braid group $B_n$, we have the \textit{chain relation}
\begin{equation*}
	(\sigma_1 \sigma_2 \cdots \sigma_{n-1})^{n} = t_{\delta},
\end{equation*}
where $\delta$ is the curve parallel to the boundary of the punctured disk.
The left-hand side can be seen as the full-twist about the $n$ strands, which in turn can be regarded as the Dehn twist along the boundary-parallel curve.

\subsection{Lantern relation}
\begin{figure}[htbp]
	\centering
	\includegraphics[height=80pt]{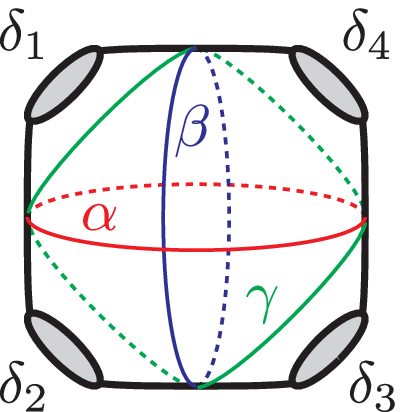}
	\caption{Four-holed sphere with boundary $\{ \delta_1, \delta_2, \delta_3, \delta_4 \}$.} \label{F:LanternRelation}	
\end{figure}

Consider the four-holed sphere in Figure~\ref{F:LanternRelation}.
Then the relation
\begin{equation*}
	t_{\alpha} t_{\beta} t_{\gamma}= t_{\delta_1} t_{\delta_2} t_{\delta_3} t_{\delta_4}
\end{equation*}
holds in $\Mod(\Sigma_0^4)$.
This is called the \textit{lantern relation}.

\subsubsection*{Lantern breeding}
\begin{figure}[htbp]
	\centering
	\includegraphics[height=180pt]{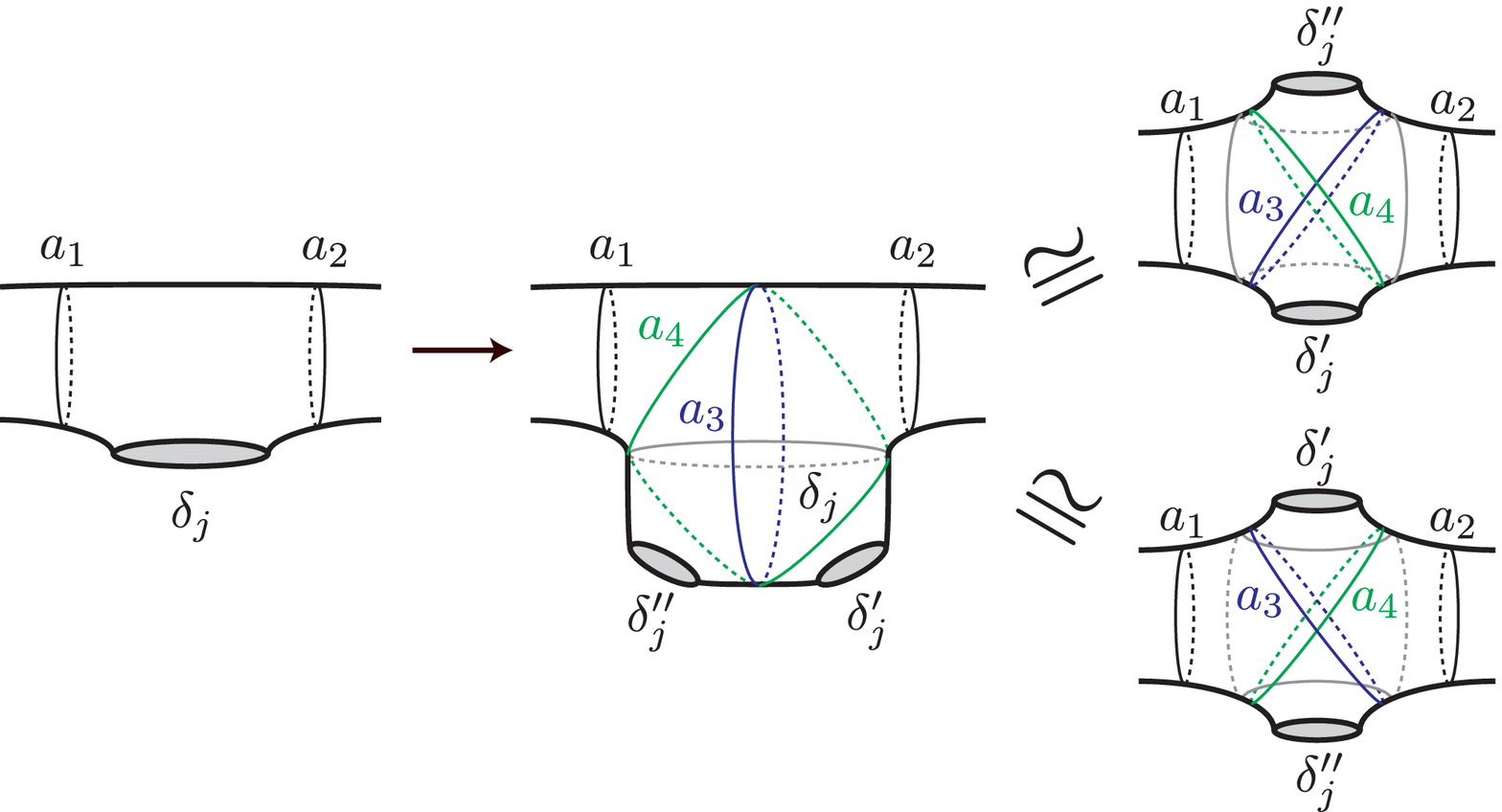}
	\caption{Lantern breeding.} \label{F:LanternBreeding}	
\end{figure}

We will often encounter the following situation, where we can effectively find a new section from an old one:
suppose that we have a positive factorization of the boundary multi-twist in $\mathrm{Mod}(\Sigma_g^k)$ in the form
\begin{equation} \label{eq:LanternBreedingBefore}
W_1 t_{a_1} t_{a_2} W_2 = t_{\delta_1} \cdots t_{\delta_j} \cdots t_{\delta_k},
\end{equation}
where the curves $a_1$, $a_2$ and $\delta_j$ bound a three-holed sphere as in the leftmost of Figure~\ref{F:LanternBreeding}.
By the commutativity, we can reform it to
\begin{equation*}
W_1 t_{\delta_j}^{-1} t_{a_1} t_{a_2} W_2 = t_{\delta_1} \cdots t_{\delta_{j-1}} t_{\delta_{j+1}} \cdots t_{\delta_k}.
\end{equation*}
Then we embed $\Sigma_g^k$ to $\Sigma_g^{k+1}$ as indicated in the middle of Figure~\ref{F:LanternBreeding} so that the curves $a_1$, $a_2$, $\delta_j^{\prime}$ and  $\delta_j^{\prime\prime}$ bound a four-holed sphere.
Thus we have a lantern relation 
$t_{a_1} t_{a_2} t_{\delta_j^{\prime}} t_{\delta_j^{\prime\prime}} = t_{\delta_j} t_{a_3} t_{a_4}$, 
which is equivalent to 
$ t_{\delta_j}^{-1} t_{a_1} t_{a_2} = t_{a_3} t_{a_4} t_{\delta_j^{\prime}}^{-1} t_{\delta_j^{\prime\prime}}^{-1}$.
By substituting this to the above equation, we get
\begin{equation*}
W_1 t_{a_3} t_{a_4} t_{\delta_j^{\prime}}^{-1} t_{\delta_j^{\prime\prime}}^{-1} W_2 = t_{\delta_1} \cdots t_{\delta_{j-1}} t_{\delta_{j+1}} \cdots t_{\delta_k},
\end{equation*}
which can be rewritten as
\begin{equation*}
W_1 t_{a_3} t_{a_4} W_2 = t_{\delta_1} \cdots t_{\delta_{j-1}} t_{\delta_j^{\prime}} t_{\delta_j^{\prime\prime}} t_{\delta_{j+1}} \cdots t_{\delta_k}. 
\end{equation*}
After renaming the boundary components, we obtain
\begin{equation} \label{eq:LanternBreedingAfter}
W_1 t_{a_3} t_{a_4} W_2 = t_{\delta_1} \cdots t_{\delta_k} t_{\delta_{k+1}}
\end{equation}
in $\mathrm{Mod}(\Sigma_g^{k+1})$, which is again a positive factorization of the boundary multi-twist with the number of boundary components increased.
Moreover, this is a further lift of the original factorization~\eqref{eq:LanternBreedingBefore}.
We call this technique to obtain the new factorization~\eqref{eq:LanternBreedingAfter} from an old factorization~\eqref{eq:LanternBreedingBefore} the \textit{lantern breeding} (with respect to $\{ a_1, a_2, \delta_j \}$).
The lantern breeding is the most basic technique to find $(-1)$-sections, which has been repeatedly used in the literature.

\section{Lifts of the MCK relation} \label{S:LiftsMCK} 

From now on, we construct several positive factorizations of the boundary multi-twist in $\Mod(\Sigma_g^k)$ where $k=4$ for even $g$ or $k=8$ for odd $g$, respectively, so that each of them is a lift of the monodromy factorization $W=1$ of the MCK Lefschetz fibration.
Let us give an outline of the construction procedure, which consists of a number of steps.
We first find a lift of $W=1$ for odd genus $g=2h-1$ to $\Mod(\Sigma_{2h-1}^2)$ by constructing a relation in $B_{4h}$ and then translating it into a relation in $\Mod(\Sigma_{2h-1}^2)$ via the Birman-Hilden double-covering.
We perform two lantern breedings to the resulting lift to obtain a further lift to $\Mod(\Sigma_{2h-1}^4)$: $W_0=\delta_1 \cdots \delta_4$.
By performing two other lantern breedings in two different ways we will have two lifts to $\Mod(\Sigma_{2h-1}^6)$: $W_{\mathrm{\rmI}}=\delta_1 \cdots \delta_6$ and $W_{\mathrm{\rmII}}=\delta_1 \cdots \delta_6$.
We will find further room to perform two more lantern breedings with several combinations, which gives six lifts to $\Mod(\Sigma_{2h-1}^8)$: 
$W_{\mathrm{\rmI}A}=\delta_1 \cdots \delta_8$, 
$W_{\mathrm{\rmI}B\sharp}=\delta_1 \cdots \delta_8$, 
$W_{\mathrm{\rmI}B\flat}=\delta_1 \cdots \delta_8$, 
$W_{\mathrm{\rmII}A}=\delta_1 \cdots \delta_8$, 
$W_{\mathrm{\rmII}B\sharp}=\delta_1 \cdots \delta_8$ and 
$W_{\mathrm{\rmII}B\flat}=\delta_1 \cdots \delta_8$.
For even genus $g=2h$, we embed $\Sigma_{2h-1}^6$ into $\Sigma_{2h}^4$ and make the best of the already-established relations $W_{\mathrm{\rmI}}=\delta_1 \cdots \delta_6$ and $W_{\mathrm{\rmII}}=\delta_1 \cdots \delta_6$ to construct four lifts of $W=1$ for $g=2h$ to $\Mod(\Sigma_{2h}^4)$: 
$W_{\mathrm{\rmI}A}=\delta_1 \cdots \delta_4$, 
$W_{\mathrm{\rmI}B}=\delta_1 \cdots \delta_4$, 
$W_{\mathrm{\rmII}A}=\delta_1 \cdots \delta_4$ and 
$W_{\mathrm{\rmII}B}=\delta_1 \cdots \delta_4$.

\subsection{A lift to $\mathrm{Mod}(\Sigma_{2h-1}^4)$ via the braid group} \
We first follow Korkmaz's original construction of the MCK relation via the braid group, but in a slightly different manner so that we can find configurations of curves for lantern breedings.

\begin{figure}[t]
	\centering
	\includegraphics[height=100pt]{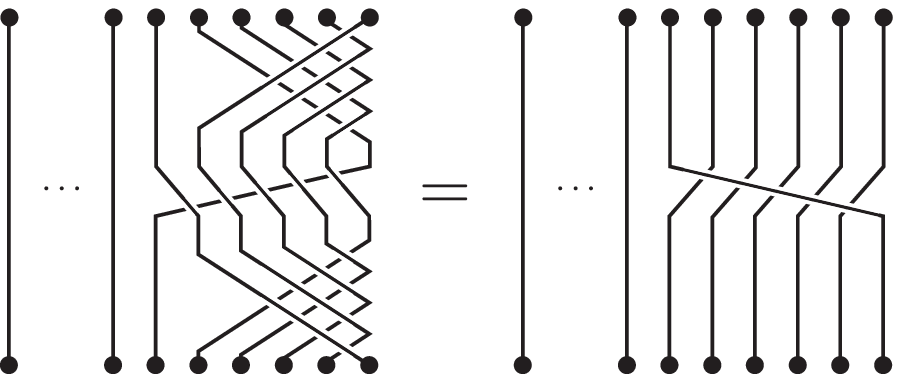}
	\caption{A conjugation of $\sigma_{2h+1} \cdots \sigma_{4h-1}$ coincides with $\sigma_{4h-1} \cdots \sigma_{2h+1}$. We read braid pictures up to down as we read braid words right to left.} \label{F:Rightleft2Leftright}	
\end{figure}

\begin{figure}[htbp]
	\centering
	\subfigure[The arcs $\Gamma(\sigma_{2h}), \Gamma^2(\sigma_{2h}), \cdots, \Gamma^{2h}(\sigma_{2h})$ and the curve $d_2$. \label{F:prebetai}]
	{\includegraphics[height=100pt]{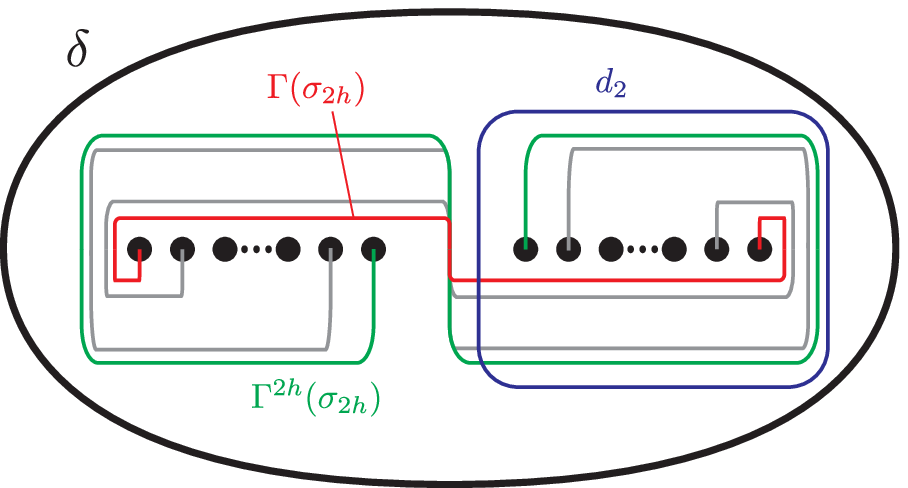}} 
	\hspace{.8em}	
	\subfigure[The arcs $\beta_0, \beta_1, \cdots, \beta_{2h-1}$ and the curves $d_1, d_2$.
	\label{F:betai}]
	{\includegraphics[height=100pt]{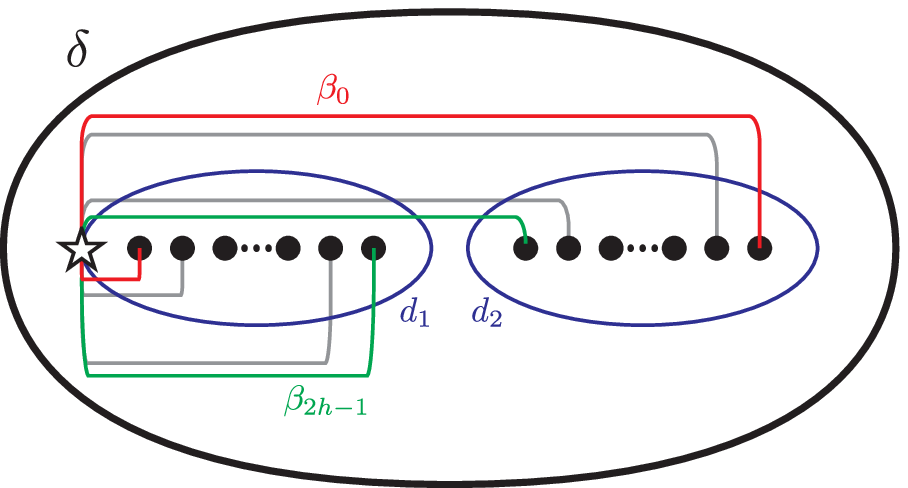}} 
	\caption{The curves for the relations in $B_{4h}$.} 	
\end{figure}

\subsubsection \
Consider the chain relation in $B_{4h}$:
\begin{equation*}
	t_{\delta} = (\sigma_1 \sigma_2 \cdots \sigma_{4h-1})^{4h}.
\end{equation*}
By a cyclic permutation and commutativity relations, we alter it as
\begin{align*}
	t_{\delta} &= (\sigma_{2h+1} \cdots \sigma_{4h-1} \cdot \sigma_1 \cdots \sigma_{2h-1} \cdot \sigma_{2h} )^{4h} \\
	&= (\sigma_1 \cdots \sigma_{2h-1} \cdot \sigma_{2h+1} \cdots \sigma_{4h-1} \cdot \sigma_{2h} )^{4h} .
\end{align*}
Then we take the simultaneous conjugation of the last expression by 
\begin{align*}
	(\sigma_{4h-1} \cdots \sigma_{2h+2}) (\sigma_{4h-1} \cdots \sigma_{2h+3}) \cdots (\sigma_{4h-1} \sigma_{4h-2}) (\sigma_{4h-1}).
\end{align*}
Since this does not affect $\sigma_1 \cdots \sigma_{2h-1}$ and $\sigma_{2h}$, we only need to see the effect on $\sigma_{2h+1} \cdots \sigma_{4h-1}$:
\begin{align*}
	& (\sigma_{4h-1} \cdots \sigma_{2h+2}) (\sigma_{4h-1} \cdots \sigma_{2h+3}) \cdots (\sigma_{4h-1} \sigma_{4h-2}) (\sigma_{4h-1}) \\
	&\qquad \cdot (\sigma_{2h+1} \cdots \underline{\sigma_{4h-1}) \cdot (\sigma_{4h-1}^{-1}}) (\sigma_{4h-2}^{-1} \sigma_{4h-1}^{-1}) \cdots (\sigma_{2h+3}^{-1} \cdots \sigma_{4h-1}^{-1}) (\sigma_{2h+2}^{-1} \cdots \sigma_{4h-1}^{-1}) \\
	=\;& (\sigma_{4h-1} \cdots \sigma_{2h+2}) (\sigma_{4h-1} \cdots \sigma_{2h+3}) \cdots (\sigma_{4h-1} \sigma_{4h-2}) (\underline{\sigma_{4h-1}}) \\
	&\qquad \cdot (\sigma_{2h+1} \cdots \underline{\sigma_{4h-2}) \cdot (\sigma_{4h-2}^{-1}} \; \underline{\sigma_{4h-1}^{-1}}) \cdots (\sigma_{2h+3}^{-1} \cdots \sigma_{4h-1}^{-1}) (\sigma_{2h+2}^{-1} \cdots \sigma_{4h-1}^{-1}) \\
	=\;& (\sigma_{4h-1} \cdots \sigma_{2h+2}) (\sigma_{4h-1} \cdots \sigma_{2h+3}) \cdots (\underline{\sigma_{4h-1} \sigma_{4h-2}}) \\
	&\qquad \cdot (\sigma_{2h+1} \cdots \underline{\sigma_{4h-3}) \cdot (\sigma_{4h-3}^{-1}} \;  \underline{\sigma_{4h-2}^{-1} \sigma_{4h-1}^{-1}}) \cdots (\sigma_{2h+3}^{-1} \cdots \sigma_{4h-1}^{-1}) (\sigma_{2h+2}^{-1} \cdots \sigma_{4h-1}^{-1}) \\
	=\;& \cdots \\
	=\;& (\sigma_{4h-1} \cdots \sigma_{2h+2}) (\underline{\sigma_{4h-1} \cdots \sigma_{2h+3}}) \cdot (\sigma_{2h+1} \underline{\sigma_{2h+2}) \cdot (\sigma_{2h+2}^{-1}} \; \underline{\sigma_{2h+3}^{-1} \cdots \sigma_{4h-1}^{-1}}) \\
	=\;& \sigma_{4h-1} \cdots \sigma_{2h+2} \sigma_{2h+1}
\end{align*}
(cf. Figure~\ref{F:Rightleft2Leftright}). 
Hence, the above simultaneous conjugation results in
\begin{align*}
	t_{\delta} 
	&= (\sigma_1 \cdots \sigma_{2h-1} \cdot \sigma_{4h-1} \cdots \sigma_{2h+1} \cdot \sigma_{2h} )^{4h}.
\end{align*}
Setting $\Gamma = \sigma_1 \cdots \sigma_{2h-1} \cdot \sigma_{4h-1} \cdots \sigma_{2h+1} = \tau_{\gamma_1} \cdots \tau_{\gamma_{2h-1}} \cdot \tau_{\gamma_{4h-1}} \cdots \tau_{\gamma_{2h+1}}$, we modify the equation as follows
\begin{align*}
	t_{\delta} 
	&= (\Gamma \sigma_{2h} )^{4h} 
	= ((\Gamma \tau_{\gamma_{2h}} )^{2h})^2 \\
	&= (\tau_{\Gamma(\gamma_{2h})} \cdot \Gamma^2 \tau_{\gamma_{2h}} \cdot  (\Gamma \tau_{\gamma_{2h}})^{2h-2})^2 \\
	&= \cdots \\
	&= (\tau_{\Gamma(\gamma_{2h})} \tau_{\Gamma^2(\gamma_{2h})} \cdots \tau_{\Gamma^{2h}(\gamma_{2h})} \cdot \Gamma^{2h})^2
\end{align*}
(cf. Figure~\ref{F:prebetai}).
In addition, by the chain relations we see 
\begin{align*}
	\Gamma^{2h}  = (\tau_{\gamma_1} \cdots \tau_{\gamma_{2h-1}})^{2h}  (\tau_{\gamma_{4h-1}} \cdots \tau_{\gamma_{2h+1}})^{2h} = t_{d_1} t_{d_2},
\end{align*}
where $d_1$ (or $d_2$) is the boundary of a small regular neighborhood of $\gamma_1 \cup \cdots \cup \tau_{\gamma_{2h-1}}$ (or $\gamma_{4h-1} \cup \cdots \cup \tau_{\gamma_{2h+1}}$), respectively. 
We thus have
\begin{align*}
	t_{\delta} = (\tau_{\Gamma(\gamma_{2h})} \tau_{\Gamma^2(\gamma_{2h})} \cdots \tau_{\Gamma^{2h}(\gamma_{2h})} t_{d_1} t_{d_2} )^2.
\end{align*}
Finally, taking the simultaneous conjugation by $t_{d_2}^{-1}$ and putting $\beta_i = t_{d_2}^{-1} \Gamma^{i+1}(\sigma_{2h})$ ($i= 0, 1, \cdots, 2h-1$), we obtain the following relation in $B_{4h}$:
\begin{align} \label{eq:MCKinBraidgroup}
	t_{\delta} = (\tau_{\beta_0} \tau_{\beta_1} \cdots \tau_{\beta_{2h-1}} t_{d_1} t_{d_2} )^2,
\end{align}
where the arcs and curves are as depicted in Figure~\ref{F:betai}.
\qed

\begin{figure}[htbp]
	\centering
	\includegraphics[height=120pt]{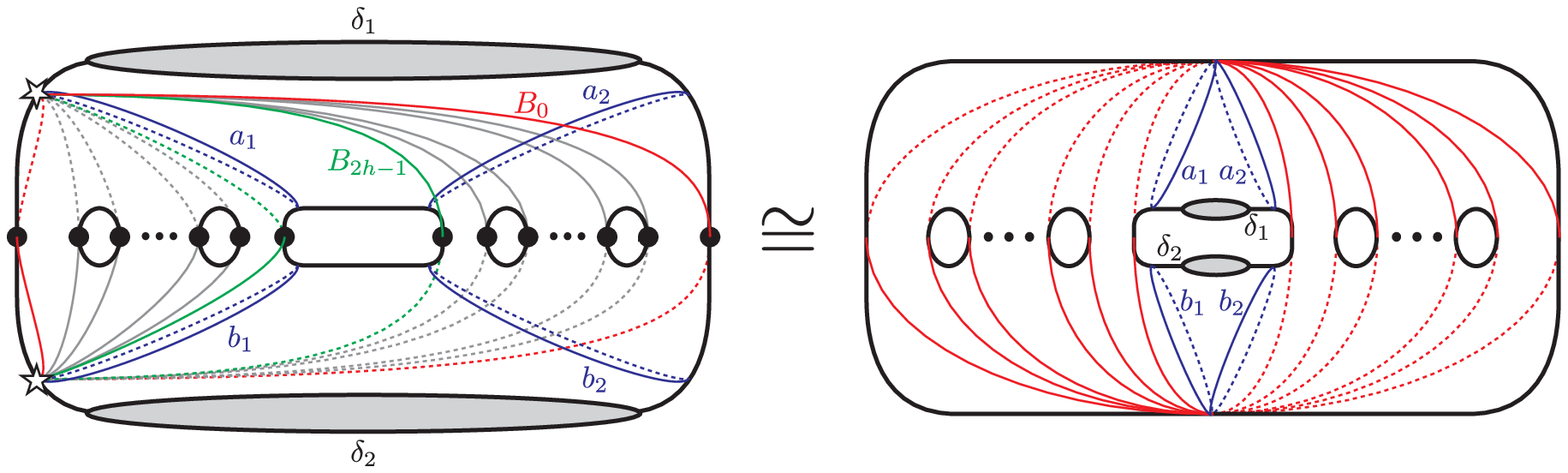}
	\caption{The curves $B_0, \cdots, B_{2h-1}, a_1, b_1, a_2, b_2$ on $\Sigma_{2h-1}^2$: the lifts of $\beta_0, \cdots, \beta_{2h-1}, d_1, d_2$.} 
	\label{F:SectionsMCK_viaBraidgroup}
\end{figure}

\subsubsection \
Nextly, we project the relation~\eqref{eq:MCKinBraidgroup} to $\Mod(\Sigma_{2h-1}^2)$ via the Birman-Hilden double covering.
We briefly recall this notion, but for the detail see~\cite{FarbMargalit2012}.
Consider the surface $\Sigma_{2h-1}^2$ in the left of Figure~\ref{F:SectionsMCK_viaBraidgroup} and the involution $\iota$ that is the rotation about the horizontal axis by $\pi$ (which is a lift of the hyperelliptic involution).
Taking the quotient yields a double branched covering $p: \Sigma_{2h-1}^2 \rightarrow D^2 = \Sigma_{2h-1}^2/\iota$ with $4h$ branched points.
We mark (the images of) these $4h$ points on $D^2$ and identify the mapping class group of the marked disk with $B_{4h}$. 
For a diffeomorphism $f : D^2 \rightarrow D^2$ preserving the marked points setwise and the boundary component pointwise, we can find a unique diffeomorphism $f^{\prime} : \Sigma_{2h-1}^2 \rightarrow \Sigma_{2h-1}^2$ such that $f \circ p  = p \circ f^{\prime}$. 
Note that, then $f^{\prime}$ commutes with $\iota$ and preserves the boundary components pointwise.
Moreover, it turns out that $[f^{\prime}] \in \Mod(\Sigma_{4h-1}^2)$ depends only on $[f] \in B_{4h}$.
In this way, we can define a homomorphism
\begin{align*}
	\Phi : B_{4h} \rightarrow \Mod(\Sigma_{4h-1}^2),
\end{align*}
which is known to be injective (and the image is the \textit{symmetric mapping class group}).
Note that the half-twist $\tau_{\gamma}$ about a proper arc $\gamma$ connecting two distinct punctures projects to the Dehn twist $t_c$ along the simple closed curve $c$ that is the lift of $\gamma$ by $p$.

We project the relation~\eqref{eq:MCKinBraidgroup} to $\Mod(\Sigma_{2h-1}^2)$ by $\Phi$.
The resulting relation is
\begin{align} \label{eq:MCKviaBraidgroup}
	t_{\delta_1} t_{\delta_2} = (t_{B_0} t_{B_1} \cdots t_{B_{2h-1}} t_{a_1}t_{b_1} t_{a_2}t_{b_2} )^2,
\end{align}
where the curve $B_i$ is the lift of $\beta_i$ and $a_j \cup b_j$ is the lift of $d_j$ as illustrated in the left of Figure~\ref{F:SectionsMCK_viaBraidgroup} (the picture can be seen in a more symmetrical manner as the right of the same Figure).
\qed
\begin{figure}[htbp]
	\centering
	\includegraphics[height=150pt]{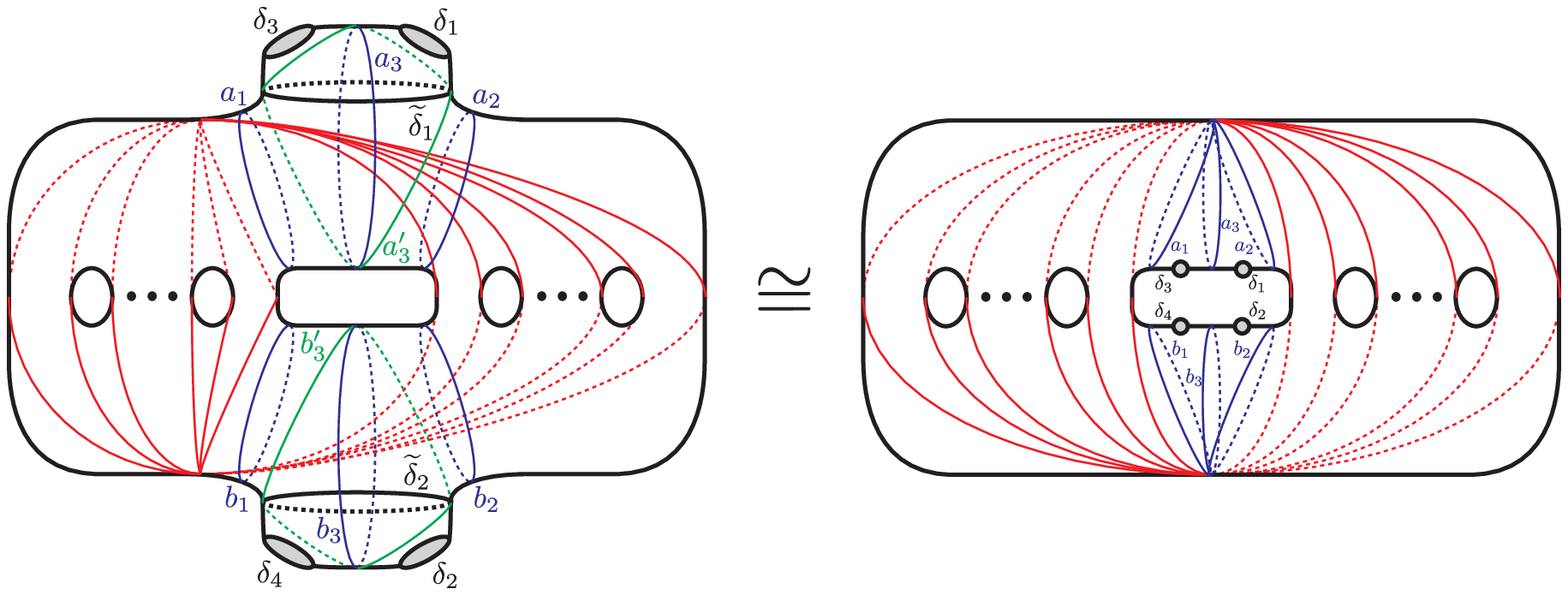}
	\caption{The curves $B_0, \cdots, B_{2h-1}, a_i, b_i$ for $W_{0}$.} \label{F:SectionsMCK_W0}
\end{figure}

\subsubsection \
In the relation~\eqref{eq:MCKviaBraidgroup}, we find two chances for lantern breeding, i.e., $\{a_1, a_2, \delta_1 \}$ and $\{b_1, b_2, \delta_2\}$.
We perform these lantern breedings as follows:
\begin{align}
	t_{\delta_1} t_{\delta_2} 
	&= t_{B_0} \cdots t_{B_{2h-1}} \underline{t_{a_1} t_{a_2}} \; \underline{t_{b_1}t_{b_2}} \cdot t_{B_0} \cdots t_{B_{2h-1}} t_{a_1}t_{b_1} t_{a_2}t_{b_2} \notag \\
	&\downarrow \notag \\
	t_{\delta_1} \cdots t_{\delta_4}
	&= t_{B_0} \cdots t_{B_{2h-1}} t_{a_3} t_{a_3^{\prime}} t_{b_3} t_{b_3^{\prime}} \cdot t_{B_0}  \cdots t_{B_{2h-1}} t_{a_1}t_{b_1} t_{a_2}t_{b_2} \notag \\
	&= t_{B_0} \cdots t_{B_{2h-1}} t_{a_3} t_{b_3} t_{a_3^{\prime}} t_{b_3^{\prime}} \cdot t_{B_0} \cdots t_{B_{2h-1}} t_{a_1}t_{b_1} t_{a_2}t_{b_2}. \label{eq:preW0}
\end{align}
Here the new relation holds in $\Mod(\Sigma_{2h-1}^4)$, where the involved curves are as in Figure~\ref{F:SectionsMCK_W0}.
\qed

\begin{figure}[h!]
	\centering
	\subfigure[Actions of $\Phi^{-1}(\Delta)$ on $d_1, d_2, d_3$. \label{F:ActionsPhiinvDelta}]
	{\includegraphics[height=80pt]{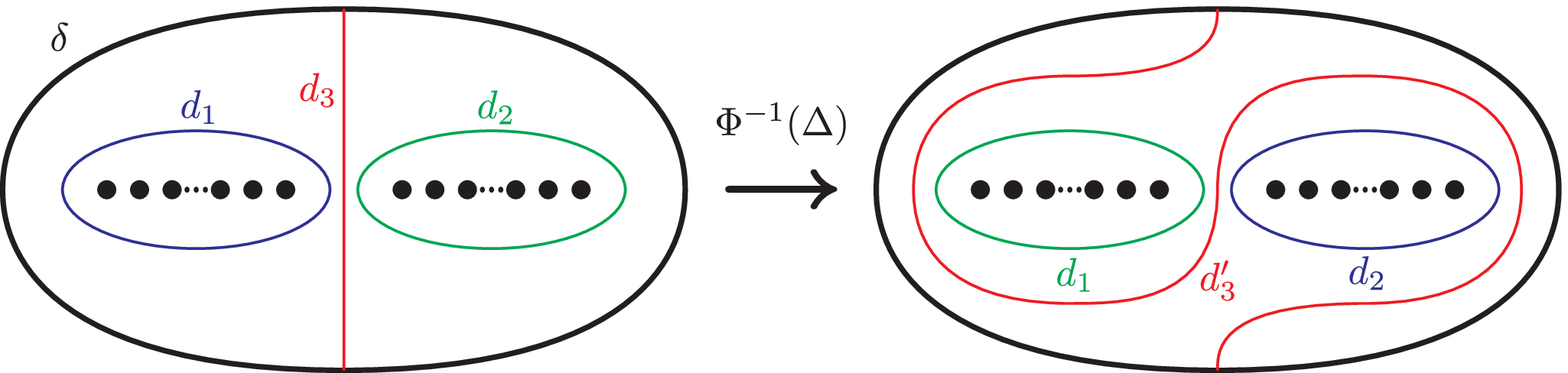}} 
	\hspace{.8em}	
	\subfigure[Actions of $\Delta$ on $a_1, a_2, a_3$. \label{F:ActionsDelta}]
	{\includegraphics[height=75pt]{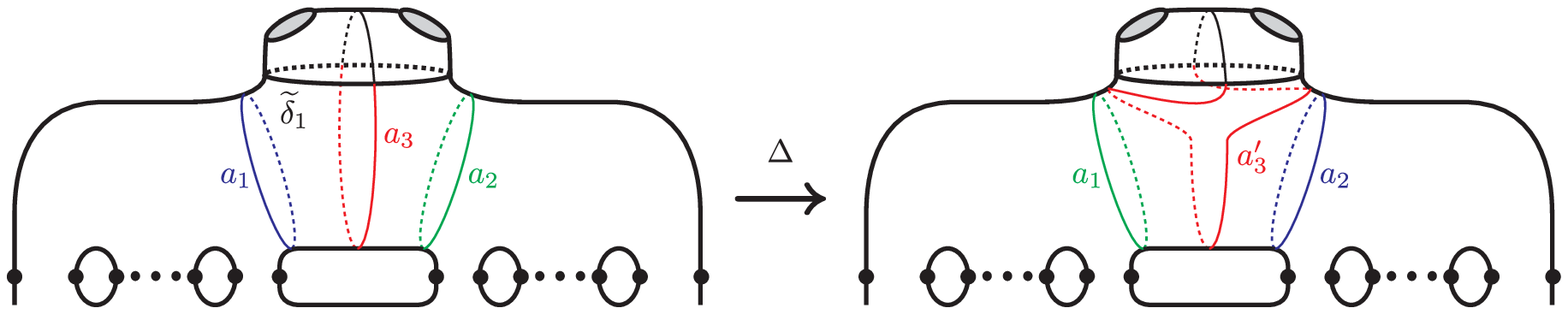}} 
	\caption{Analysis of actions of $\Delta = t_{B_0} \cdots t_{B_{2h-1}}$.}
\end{figure}

\subsubsection \
We can indeed simplify the relation~\eqref{eq:preW0}, yet, to do so we need to examine the effects of some Hurwitz moves.
Let us write $\Delta = t_{B_0} \cdots t_{B_{2h-1}}$ and examine the actions of $\Delta$ on the curves $a_i, b_i$ and $a_3^{\prime}, b_3^{\prime}$.
To study these actions, we can reduce the argument in the mapping class group to that in the braid group; we observe the actions of 
$\Phi^{-1}({\Delta}) = \tau_{\beta_0} \tau_{\beta_1} \cdots \tau_{\beta_{2h-1}}$ on the curves $d_1$ and $d_2$ and the arc $d_3$ in Figure~\ref{F:ActionsPhiinvDelta}.
Here we regard $\Phi$ as the natural map $B_{4h} \rightarrow \Mod(\Sigma_{2h-1}^2) \rightarrow \Mod(\Sigma_{2h-1}^4)$ where $\Sigma_{2h-1}^2$ is the subsurface of $\Sigma_{2h-1}^4$ bounded by $\widetilde{\delta}_1$ and $\widetilde{\delta}_2$ (the original boundary before the lantern breedings) as indicated in the left of Figure~\ref{F:SectionsMCK_W0}.
As long as considering actions on $d_1$, $d_2$ and $d_3$, the action of $\tau_{\beta_0} \tau_{\beta_1} \cdots \tau_{\beta_{2h-1}}$ has the same effect as that of $\omega = \tau_{\beta_0} \tau_{\beta_1} \cdots \tau_{\beta_{2h-1}} t_{d_1} t_{d_2}$ does since $t_{d_1} t_{d_2}$ acts trivially on $d_1$, $d_2$ and $d_3$.
From the construction of the relation~\eqref{eq:MCKinBraidgroup}, we see
\begin{align*}
	\omega = \tau_{\beta_0} \tau_{\beta_1} \cdots \tau_{\beta_{2h-1}} t_{d_1} t_{d_2} = t_{d_2}^{-1} (\Gamma \tau_{\gamma_{2h}} )^{2h} t_{d_2} = t_{d_2}^{-1} (\sigma_1 \cdots \sigma_{2h-1} \cdot \sigma_{4h-1} \cdots \sigma_{2h+1} \cdot \sigma_{2h} )^{2h} t_{d_2}.
\end{align*}
\begin{figure}[htbp]
	\centering
	\includegraphics[height=100pt]{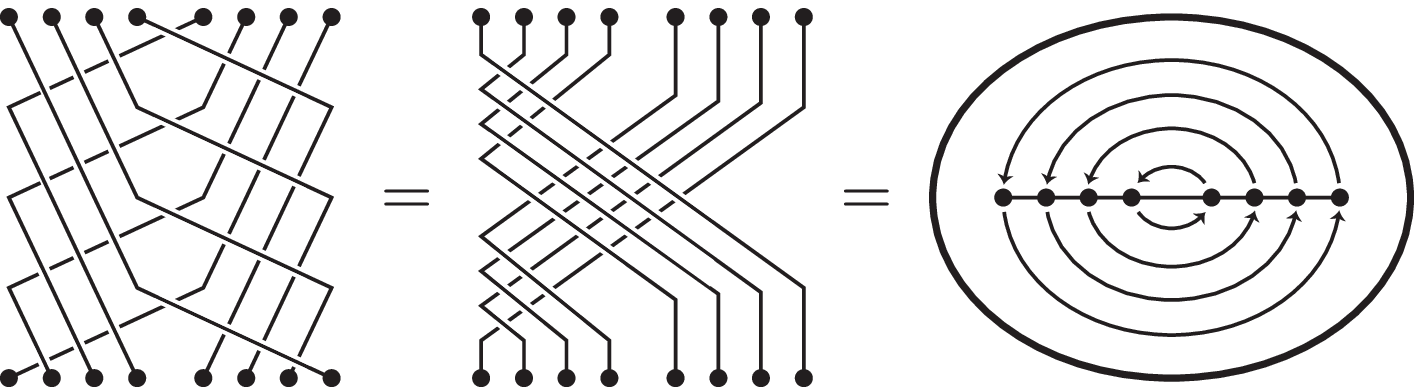}
	\caption{$(\sigma_1 \sigma_2 \sigma_3 \cdot \sigma_7 \sigma_6 \sigma_5 \cdot \sigma_4)^4 = (\sigma_1) (\sigma_2 \sigma_1) \cdots (\sigma_7 \cdots \sigma_2 \sigma_1) = \text{the half-twist}$.} 
	\label{F:GlobalHalftwist}
\end{figure}

\noindent
We now consider the action of the part 
$ (\sigma_1 \cdots \sigma_{2h-1} \cdot \sigma_{4h-1} \cdots \sigma_{2h+1} \cdot \sigma_{2h} )^{2h}$.
First, we claim that 
\begin{align*}
	(\sigma_1 \cdots \sigma_{2h-1} \cdot \sigma_{4h-1} \cdots \sigma_{2h+1} \cdot \sigma_{2h} )^{2h} 
	= (\sigma_1) (\sigma_2 \sigma_1) (\sigma_3 \sigma_2 \sigma_1) \cdots (\sigma_{4h-1} \cdots \sigma_1)
\end{align*}
(cf. Figure~\ref{F:GlobalHalftwist}).
The right hand side is ``the half-twist about the $4h$ strands'' and as a diffeomorphism it can be described as in the rightmost of Figure~\ref{F:GlobalHalftwist}.
We can verify the above equation as follows: 
\begin{align*} 
	& \text{(LHS)} = \textstyle \left(\prod\limits_{i=1}^{2h-1}\sigma_i \prod\limits_{j=1}^{2h-1}\sigma_{4h-j} \cdot \sigma_{2h} \right)^{2h} 
	= \left(\prod\limits_{i=1}^{2h-1}\sigma_i \prod\limits_{j=1}^{2h}\sigma_{4h-j} \right)^{2h} \\
	&= \textstyle \prod\limits_{i=1}^{2h-1}\sigma_i \underline{\prod\limits_{j=1}^{2h}\sigma_{4h-j} \cdot \prod\limits_{i=1}^{2h-2}\sigma_i} \cdot \sigma_{2h-1} \cdot \prod\limits_{j=1}^{2h}\sigma_{4h-j} \left(\prod\limits_{i=1}^{2h-1}\sigma_i \prod\limits_{j=1}^{2h}\sigma_{4h-j} \right)^{2h-2} \\
	&= \textstyle \prod\limits_{i=1}^{2h-1}\sigma_i \prod\limits_{i=1}^{2h-2}\sigma_i \left( \prod\limits_{j=1}^{2h+1}\sigma_{4h-j} \prod\limits_{j=1}^{2h}\sigma_{4h-j} \right) \left(\prod\limits_{i=1}^{2h-1}\sigma_i \prod\limits_{j=1}^{2h}\sigma_{4h-j} \right)^{2h-2} \\
	&= \textstyle \prod\limits_{i=1}^{2h-1}\sigma_i \prod\limits_{i=1}^{2h-2}\sigma_i \underline{\left( \prod\limits_{j=1}^{2h+1}\sigma_{4h-j} \prod\limits_{j=1}^{2h}\sigma_{4h-j} \right) \cdot \prod\limits_{i=1}^{2h-3}\sigma_i} \cdot \underline{\sigma_{2h-2}} \cdot \sigma_{2h-1} \cdot \prod\limits_{j=1}^{2h}\sigma_{4h-j} \\
	&\qquad \textstyle \cdot \left(\prod\limits_{i=1}^{2h-1}\sigma_i \prod\limits_{j=1}^{2h}\sigma_{4h-j} \right)^{2h-3} \\
	&= \textstyle \prod\limits_{i=1}^{2h-1}\sigma_i \prod\limits_{i=1}^{2h-2}\sigma_i \prod\limits_{i=1}^{2h-3}\sigma_i \left( \prod\limits_{j=1}^{2h+2}\sigma_{4h-j} \prod\limits_{j=1}^{2h+1}\sigma_{4h-j} \prod\limits_{j=1}^{2h}\sigma_{4h-j} \right) \left(\prod\limits_{i=1}^{2h-1}\sigma_i \prod\limits_{j=1}^{2h}\sigma_{4h-j} \right)^{2h-3} \\
	&= \cdots \\
	&= \textstyle \prod\limits_{i=1}^{2h-1}\sigma_i \prod\limits_{i=1}^{2h-2}\sigma_i \prod\limits_{i=1}^{2h-3}\sigma_i \cdots \prod\limits_{i=1}^{1}\sigma_i \left( \underline{\prod\limits_{j=1}^{4h-1}\sigma_{4h-j} \cdot \prod\limits_{j=1}^{4h-2}\sigma_{4h-j}} \cdot \prod\limits_{j=1}^{4h-3}\sigma_{4h-j} \cdots \prod\limits_{j=1}^{2h}\sigma_{4h-j} \right) \\
	&= \textstyle \prod\limits_{i=1}^{2h-1}\sigma_i \prod\limits_{i=1}^{2h-2}\sigma_i \cdots \prod\limits_{i=1}^{1}\sigma_i \left( \underline{\prod\limits_{j=1}^{4h-2}\sigma_{4h-1-j} \prod\limits_{j=1}^{4h-1}\sigma_{4h-j} \cdot \prod\limits_{j=1}^{4h-3}\sigma_{4h-j}} \cdot \prod\limits_{j=1}^{4h-4}\sigma_{4h-j} \cdots \prod\limits_{j=1}^{2h}\sigma_{4h-j} \right) \\
	&= \textstyle \prod\limits_{i=1}^{2h-1}\sigma_i \prod\limits_{i=1}^{2h-2}\sigma_i \cdots \prod\limits_{i=1}^{1}\sigma_i \left( \prod\limits_{j=1}^{4h-3}\sigma_{4h-2-j} \prod\limits_{j=1}^{4h-2}\sigma_{4h-1-j} \prod\limits_{j=1}^{4h-1}\sigma_{4h-j} \cdot \prod\limits_{j=1}^{4h-4}\sigma_{4h-j} \cdots \prod\limits_{j=1}^{2h}\sigma_{4h-j} \right) \\
	&= \cdots \\
	&= \textstyle \prod\limits_{i=1}^{2h-1}\sigma_i \prod\limits_{i=1}^{2h-2}\sigma_i \cdots \prod\limits_{i=1}^{1}\sigma_i \left( \prod\limits_{j=1}^{2h}\sigma_{2h+1-j} \cdots \prod\limits_{j=1}^{4h-3}\sigma_{4h-2-j} \prod\limits_{j=1}^{4h-2}\sigma_{4h-1-j} \prod\limits_{j=1}^{4h-1}\sigma_{4h-j} \right) \\
	&= \textstyle \prod\limits_{i=1}^{2h-2}\sigma_i \cdot \underline{\sigma_{2h-1}} \cdot \prod\limits_{i=1}^{2h-3}\sigma_i \cdot \underline{\sigma_{2h-2}} \cdots \prod\limits_{i=1}^{1}\sigma_i \cdot \underline{\sigma_2} \cdot \underline{\sigma_1} \cdot \left( \prod\limits_{j=1}^{2h}\sigma_{2h+1-j} \cdots \prod\limits_{j=1}^{4h-2}\sigma_{4h-1-j} \prod\limits_{j=1}^{4h-1}\sigma_{4h-j} \right) \\
	&= \textstyle \prod\limits_{i=1}^{2h-2}\sigma_i \prod\limits_{i=1}^{2h-3}\sigma_i \cdots \prod\limits_{i=1}^{1}\sigma_i \left( \prod\limits_{j=1}^{2h-1}\sigma_{2h-j} \cdot \prod\limits_{j=1}^{2h}\sigma_{2h+1-j} \cdots \prod\limits_{j=1}^{4h-2}\sigma_{4h-1-j} \prod\limits_{j=1}^{4h-1}\sigma_{4h-j} \right) \\
	&= \cdots \\
	&= \textstyle \prod\limits_{j=1}^{1}\sigma_{2-j} \cdots \prod\limits_{j=1}^{4h-2}\sigma_{4h-1-j} \prod\limits_{j=1}^{4h-1}\sigma_{4h-j}
	= \text{(RHS)}. \qed
\end{align*}

From its geometric description in the rightmost of Figure~\ref{F:GlobalHalftwist}, it is obvious that the map 
$(\sigma_1) (\sigma_2 \sigma_1) \cdots (\sigma_{4h-1} \cdots \sigma_1)$ sends $d_1$, $d_2$, $d_3$ to $d_2$, $d_1$, $d_3^{\prime}$, respectively.
So does $\omega$ and then $\Phi^{-1}({\Delta})$:
\begin{align*}
	\Phi^{-1}({\Delta})(d_1) &= d_2 & \Phi^{-1}({\Delta})(d_2) &= d_1 & \Phi^{-1}({\Delta})(d_3) &= d_3^{\prime}
\end{align*}
(cf. Figure~\ref{F:ActionsPhiinvDelta}).
Considering the geometric description again, it is easy to see that $\Delta$ sends the upper half of the surface $\Sigma_{2h-1}^4$ to the upper half and the lower half to the lower half.
Since $a_i \cup b_i$ is the lift of $d_i$ ($i=1,2$), $\Delta(a_1)$ must be $a_2$ or $b_2$.
Meanwhile $a_1$ and $a_2$ are in the upper half while $b_2$ are not. We can thus conclude that $\Delta(a_1) = a_2$.
The curve $a_3$ contains one component of the lift of $d_3$ and the rest arc of $a_3$ lies out of the support of $\Delta$ (see Figure~\ref{F:ActionsDelta}).
It follows that $\Delta(a_3)$ consists of one component of the lift of $d_3^{\prime}$ and the same arc above.
The resulting curve is indeed $a_3^{\prime}$.
By the similar arguments, we can observe the actions of $\Delta$ on the other curves.
Summarizing, we have the following properties
\begin{align*}
	\Delta(a_1) &= a_2,  &  \Delta(b_1) &= b_2, \\
	\Delta(a_2) &= a_1,  &  \Delta(b_2) &= b_1, \\
	\Delta(a_3) &= a_3^{\prime},  &  \Delta(b_3) &= b_3^{\prime}.
\end{align*}

\subsubsection \
Using the above properties, we can simplify the relation~\eqref{eq:preW0} as follows
\begin{align*}
	t_{\delta_1} \cdots t_{\delta_4} 
	&= t_{B_0} \cdots t_{B_{2h-1}} t_{a_3} t_{b_3} \underline{t_{a_3^{\prime}} t_{b_3^{\prime}} \cdot t_{B_0} \cdots t_{B_{2h-1}}} t_{a_1}t_{b_1} t_{a_2}t_{b_2} \\
	&= t_{B_0} \cdots t_{B_{2h-1}} t_{a_3} t_{b_3}\cdot t_{B_0} \cdots t_{B_{2h-1}} t_{\Delta^{-1}(a_3^{\prime})} t_{\Delta^{-1}(b_3^{\prime})}  t_{a_1}t_{b_1} t_{a_2}t_{b_2} \\
	&= t_{B_0} \cdots t_{B_{2h-1}} t_{a_3} t_{b_3}\cdot \underline{t_{B_0} \cdots t_{B_{2h-1}}} t_{a_3} t_{b_3}  t_{a_1}t_{b_1} \underline{t_{a_2}t_{b_2}} \\
	&= t_{B_0} \cdots t_{B_{2h-1}}  t_{a_3} t_{b_3} t_{\Delta(a_2)}t_{\Delta(b_2)} \cdot t_{B_0} \cdots t_{B_{2h-1}} t_{a_3} t_{b_3}  t_{a_1}t_{b_1} \\
	&= (t_{B_0} \cdots t_{B_{2h-1}} t_{a_3} t_{b_3} t_{a_1}t_{b_1} )^2.
\end{align*} 
By similar modifications, we obtain 
\begin{align}
	t_{\delta_1} \cdots t_{\delta_4} 
	&= (t_{B_0} \cdots t_{B_{2h-1}} t_{a_1} t_{a_3} t_{b_1} t_{b_3})^2, \label{eq:W0forWI} \\ 
	t_{\delta_1} \cdots t_{\delta_4} 
	&= (t_{B_0} \cdots t_{B_{2h-1}} t_{a_2} t_{a_3} t_{b_2} t_{b_3})^2, \\
	t_{\delta_1} \cdots t_{\delta_4} 
	&= (t_{B_0} \cdots t_{B_{2h-1}} t_{a_1} t_{a_3} t_{b_2} t_{b_3})^2 \label{eq:W0forWII}
\end{align}
(cf. Figure~\ref{F:SectionsMCK_W0}).
Note that the above three factorizations are Hurwitz equivalent since we only used Hurwitz moves.
We will refer to these Hurwitz equivalent words as $W_0$.
\qed

\subsection{Two lifts to $\mathrm{Mod}(\Sigma_{2h-1}^6)$}
We further perform two more lantern breedings to $W_0 = t_{\delta_1} \cdots t_{\delta_4}$ using two different configurations.

\begin{figure}[htbp]
	\centering
	\includegraphics[height=120pt]{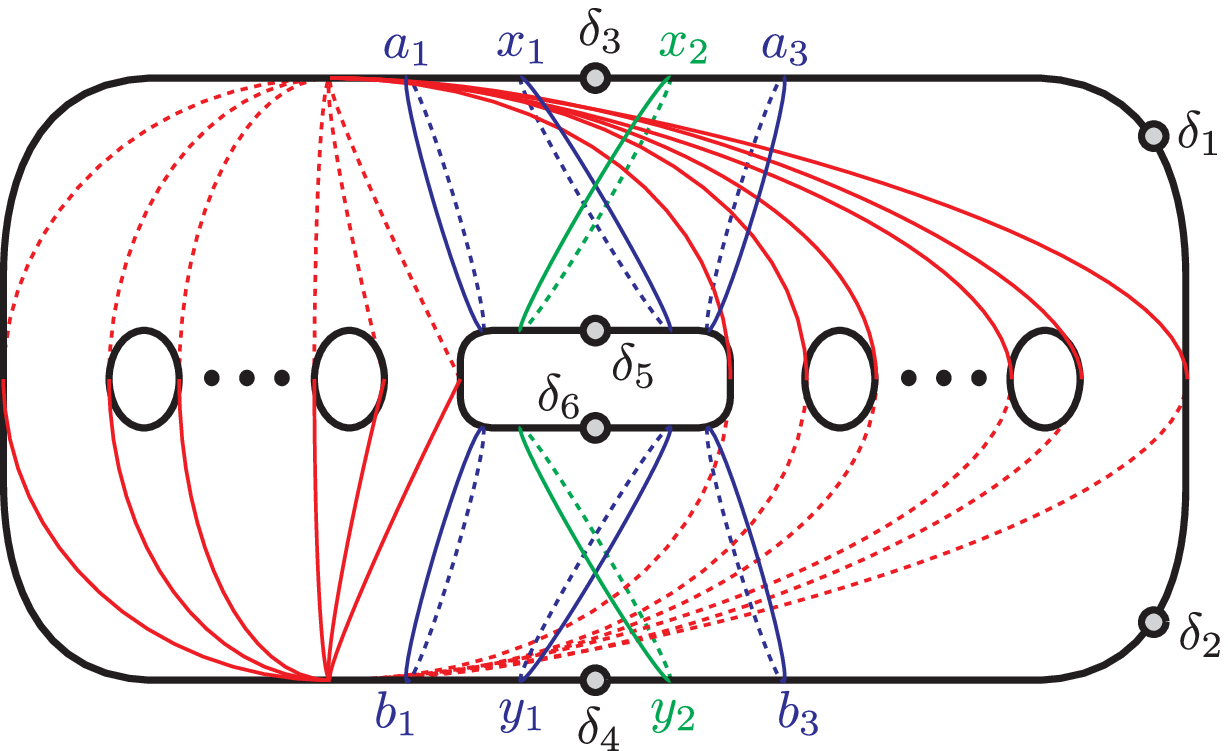}
	\caption{The curves $B_0, \cdots, B_{2h-1}$, $a_i,b_i,x_i,y_i$ for $W_{\mathrm{\rmI}}$.} \label{F:WI}
\end{figure}

\subsubsection \
We first use $\{ a_1, a_3, \delta_3 \}$ and $\{ b_1, b_3, \delta_4 \}$ with the expression~\eqref{eq:W0forWI} for $W_0$:
\begin{align}
	t_{\delta_1} \cdots t_{\delta_4} 
	&= t_{B_0} \cdots t_{B_{2h-1}} t_{a_1} t_{a_3} t_{b_1} t_{b_3} \cdot t_{B_0} \cdots t_{B_{2h-1}} \underline{t_{a_1} t_{a_3}} \; \underline{t_{b_1} t_{b_3}} \notag \\ 
	&\downarrow \notag \\ 
	t_{\delta_1} \cdots t_{\delta_6} 
	&= t_{B_0} \cdots t_{B_{2h-1}} t_{a_1} t_{a_3} t_{b_1} t_{b_3} \cdot t_{B_0} \cdots t_{B_{2h-1}} t_{x_1} t_{x_2} t_{y_1} t_{y_2}, 
\end{align}
where the curves are shown in Figure~\ref{F:WI}.
We denote by $W_{\mathrm{\rmI}}$ the last word.
\qed

\begin{figure}[htbp]
	\centering
	\includegraphics[height=120pt]{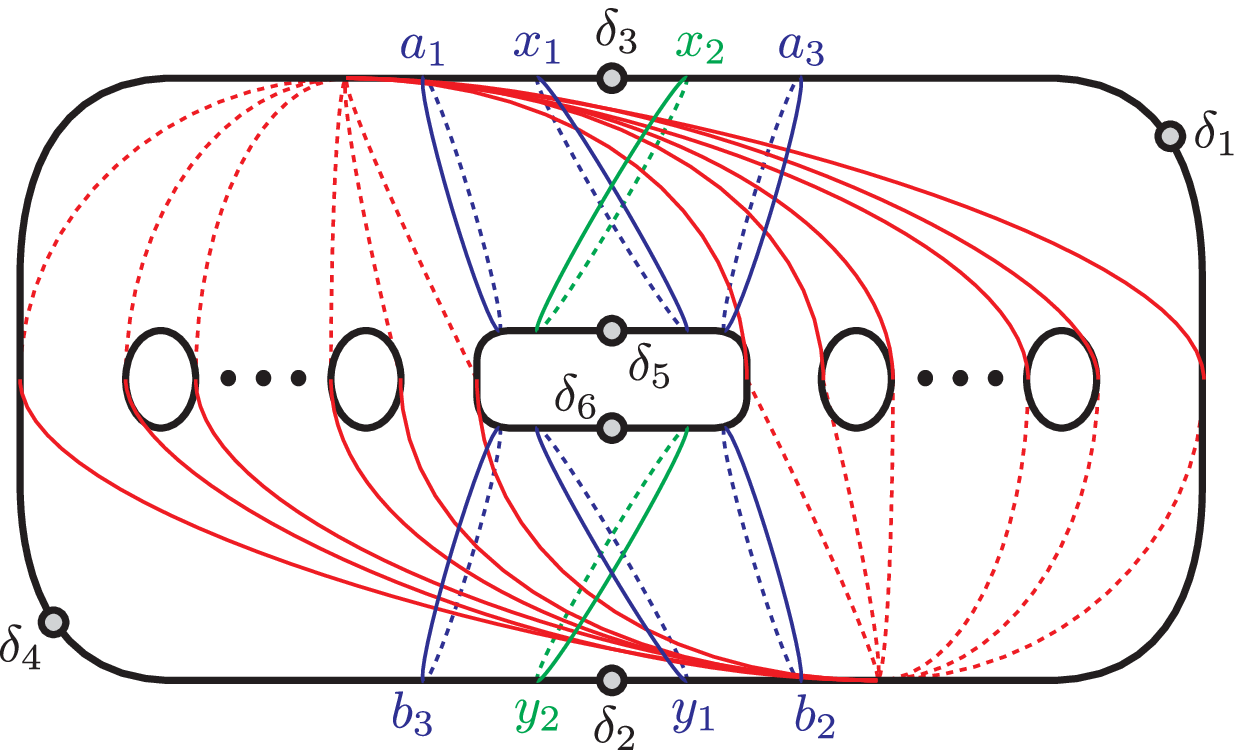}
	\caption{The curves $B_0^{\prime}, \cdots, B_{2h-1}^{\prime}$, $a_i,b_i,x_i,y_i$ for $W_{\mathrm{\rmII}}$.} \label{F:WII}
\end{figure}

\subsubsection \
Secondary, we use $\{ a_1, a_3, \delta_3 \}$ and $\{ b_2, b_3, \delta_2 \}$ with the expression~\eqref{eq:W0forWII} for $W_0$.
Before doing lantern breedings, we take the conjugation of the expression~\eqref{eq:W0forWII} by $t_{b_1}^{-1}t_{b_2}$ as
\begin{align*}
	t_{\delta_1} \cdots t_{\delta_4} 
	&= (t_{B_0^{\prime}} \cdots t_{B_{2h-1}^{\prime}} t_{a_1} t_{a_3} t_{b_2} t_{b_3})^2, 
\end{align*}
where $B_i^{\prime} = t_{b_1}^{-1}t_{b_2}(B_i)$.
Then perform the lantern breedings:
\begin{align}
	t_{\delta_1} \cdots t_{\delta_4} 
	&= t_{B_0^{\prime}} \cdots t_{B_{2h-1}^{\prime}} t_{a_1} t_{a_3} t_{b_2} t_{b_3} \cdot t_{B_0^{\prime}} \cdots t_{B_{2h-1}^{\prime}} \underline{t_{a_1} t_{a_3}} \; \underline{t_{b_2} t_{b_3}} \notag \\
	&\downarrow \notag \\
	t_{\delta_1} \cdots t_{\delta_6} 
	&= t_{B_0^{\prime}} \cdots t_{B_{2h-1}^{\prime}} t_{a_1} t_{a_3} t_{b_2} t_{b_3} \cdot t_{B_0^{\prime}} \cdots t_{B_{2h-1}^{\prime}} t_{x_1} t_{x_2} t_{y_1} t_{y_2}, 
\end{align}
where the curves are shown in Figure~\ref{F:WII}.
We denote by $W_{\mathrm{\rmII}}$ the last word.
\qed

\subsection{The case of odd genus: $g= 2h-1$} \label{subsect:odgenus}
For odd $g = 2h-1$, we finally construct lifts of the MCK relation to $\Sigma_g^8$.
There are several combinations to take further lantern breedings to $W_{\mathrm{\rmI}} = t_{\delta_1} \cdots t_{\delta_6}$ and $W_{\mathrm{\rmII}} = t_{\delta_1} \cdots t_{\delta_6}$.
Due to these combinations, we will have several lifts.
In the following calculations, ${}_{\varphi}(\psi)$ means the conjugate of $\psi$ by $\varphi$, i.e., ${}_{\varphi}(\psi)=\varphi \psi \varphi^{-1}$.

\subsubsection{Type $\mathrm{\rmI}A$}
Perform the lantern breedings to $W_{\mathrm{\rmI}} = t_{\delta_1} \cdots t_{\delta_6}$ using $\{  a_1, x_1, \delta_5 \}$ and $\{ b_3, y_2, \delta_6 \}$ as
\begin{align}
	t_{\delta_1} \cdots t_{\delta_6} 
	&= \underline{t_{B_0} \cdots t_{B_{g}}} t_{a_1} t_{a_3} t_{b_1} t_{b_3} \cdot \underline{t_{B_0} \cdots t_{B_{g}} t_{x_1}} t_{x_2} t_{y_1} \underline{t_{y_2}} \notag \\
	&= {}_{t_{y_2}}(t_{B_0} \cdots t_{B_{g}}) \; t_{y_2} t_{a_1} t_{a_3} t_{b_1} t_{b_3} t_{x_1} \cdot {}_{t_{x_1}^{-1}}(t_{B_0} \cdots t_{B_{g}}) \; t_{x_2} t_{y_1} \notag \\
	&= \underline{{}_{t_{y_2}}(t_{B_0} \cdots t_{B_{g}}) \; t_{a_3}} t_{a_1} t_{x_1} t_{y_2} t_{b_3} \underline{t_{b_1} \cdot {}_{t_{x_1}^{-1}}(t_{B_0} \cdots t_{B_{g}}) } \; t_{x_2} t_{y_1} \notag \\
	&= {}_{t_{a_3}^{-1}t_{y_2}}(t_{B_0} \cdots t_{B_{g}}) \;  t_{a_1} t_{x_1} t_{y_2} t_{b_3} \cdot {}_{t_{b_1} t_{x_1}^{-1}}(t_{B_0} \cdots t_{B_{g}}) \; t_{b_1} t_{x_2} t_{y_1} t_{a_3} \notag \\
 	&= {}_{t_{a_3}^{-1}t_{y_2}}(t_{B_0} \cdots t_{B_{g}}) \; \underline{t_{a_1} t_{x_1}} \; \underline{t_{y_2} t_{b_3}} \cdot {}_{t_{b_1} t_{x_1}^{-1}}(t_{B_0} \cdots t_{B_{g}}) \;  t_{x_2} t_{a_3} t_{y_1} t_{b_1}   \notag \\
	&\downarrow \notag \\
	t_{\delta_1} \cdots t_{\delta_8} 
	&= t_{B_{0,1}} t_{B_{1,1}} \cdots t_{B_{g,1}} t_{a_1} t_{a_2} t_{b_1} t_{b_2} \cdot 
	t_{B_{0,2}} t_{B_{1,2}} \cdots t_{B_{g,2}} t_{a_3} t_{a_4} t_{b_3} t_{b_4}, \label{eq:WIAodd}
\end{align}
where we have renamed the last curves (the curves $a_i$ and $b_j$ in the last formula are not necessarily the same as the previous ones) and they are depicted in Figure~\ref{F:WIAodd}.
(In the following Figures, the circle labeled $j$ represents the $j$-th boundary component $\delta_j$.)
We denote by $W_{\mathrm{\rmI}A}$ the right hand side of~\ref{eq:WIAodd}.
\qed

\begin{figure}[htbp]
	\centering
	\subfigure[$B_{0,1}, B_{1,1}, \cdots, B_{g,1}, a_1, a_2, b_1, b_2$.]
	{\includegraphics[height=110pt]{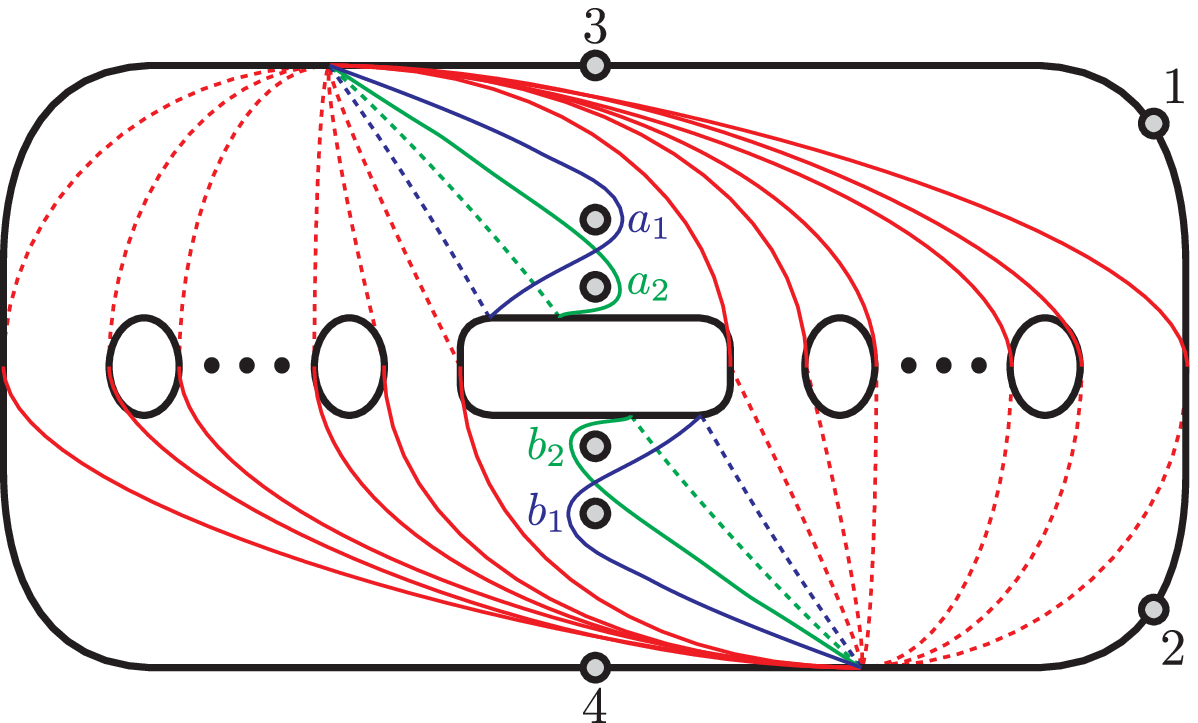}} 
	\hspace{.8em}	
	\subfigure[$B_{0,2}, B_{1,2}, \cdots, B_{g,2}, a_3, a_4, b_3, b_4$.]
	{\includegraphics[height=110pt]{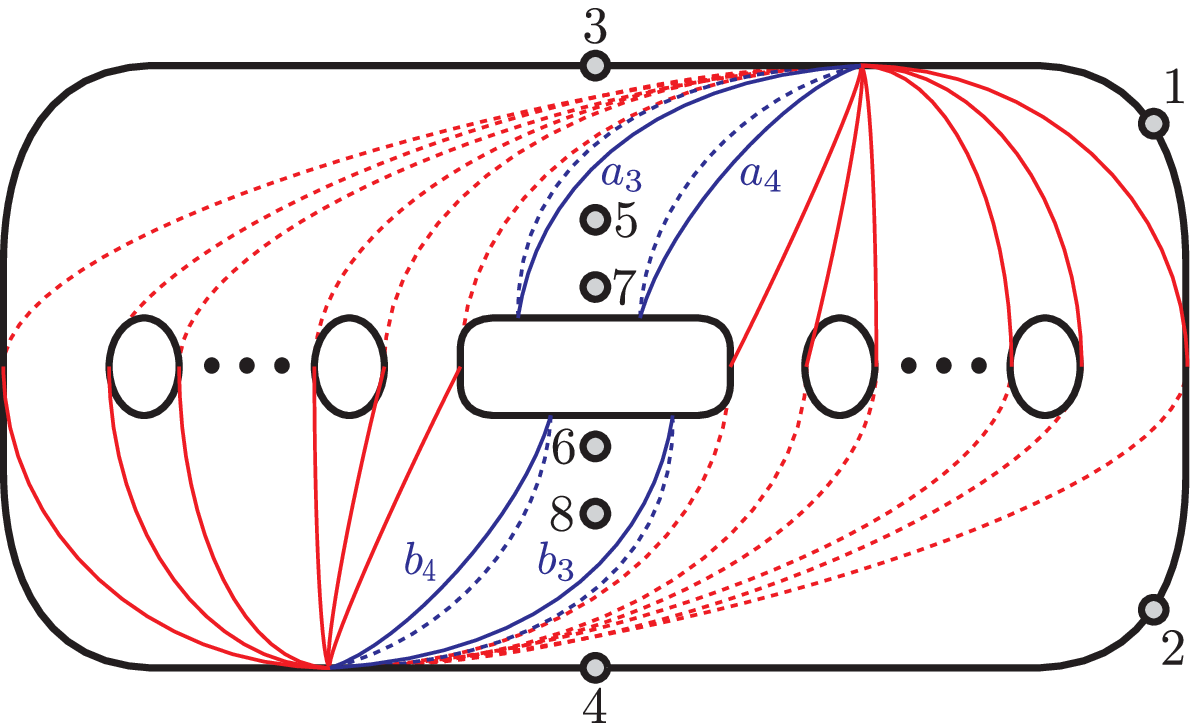}} 
	\vspace{-\baselineskip}
	\caption{The curves for $W_{\mathrm{\rmI}A}$ of odd $g$.} \label{F:WIAodd}
\end{figure}
\begin{figure}[htbp]
	\centering	 	
	\subfigure[$B_{0,1}, B_{1,1}, \cdots, B_{g,1}, a_1, a_2, b_1, b_2$.]
	{\includegraphics[height=110pt]{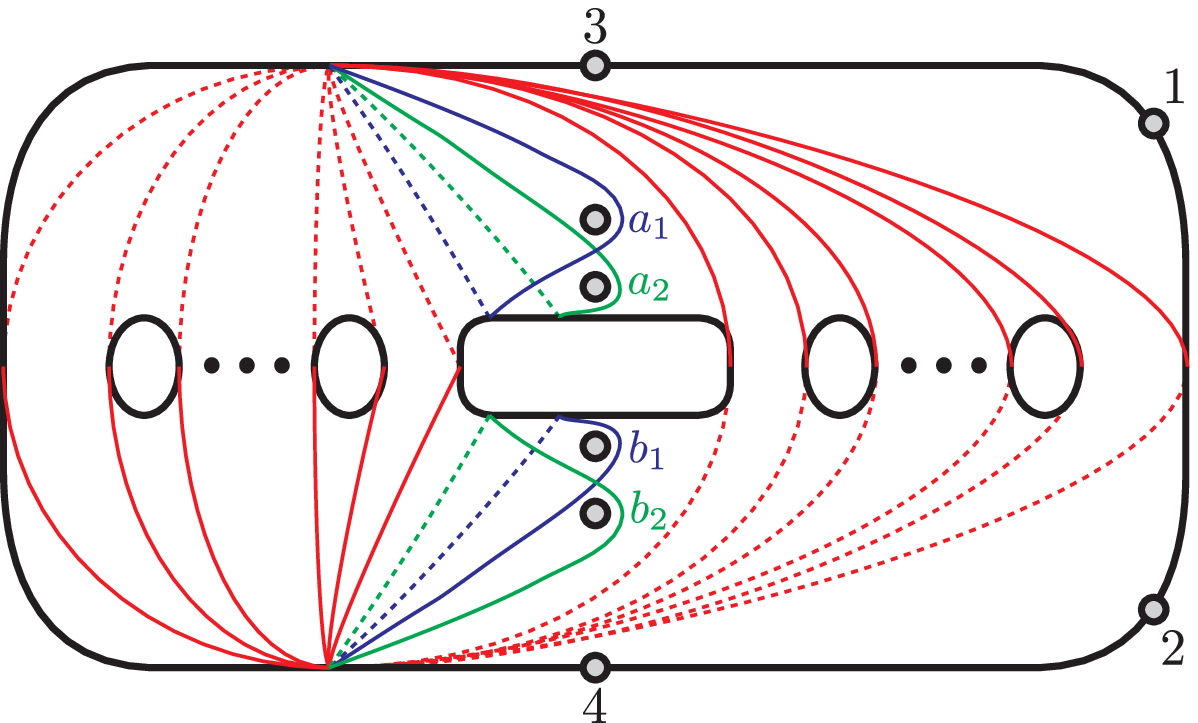}} 
	\hspace{.8em}	
	\subfigure[$B_{0,2}, B_{1,2}, \cdots, B_{g,2}, a_3, a_4, b_3, b_4$.]
	{\includegraphics[height=110pt]{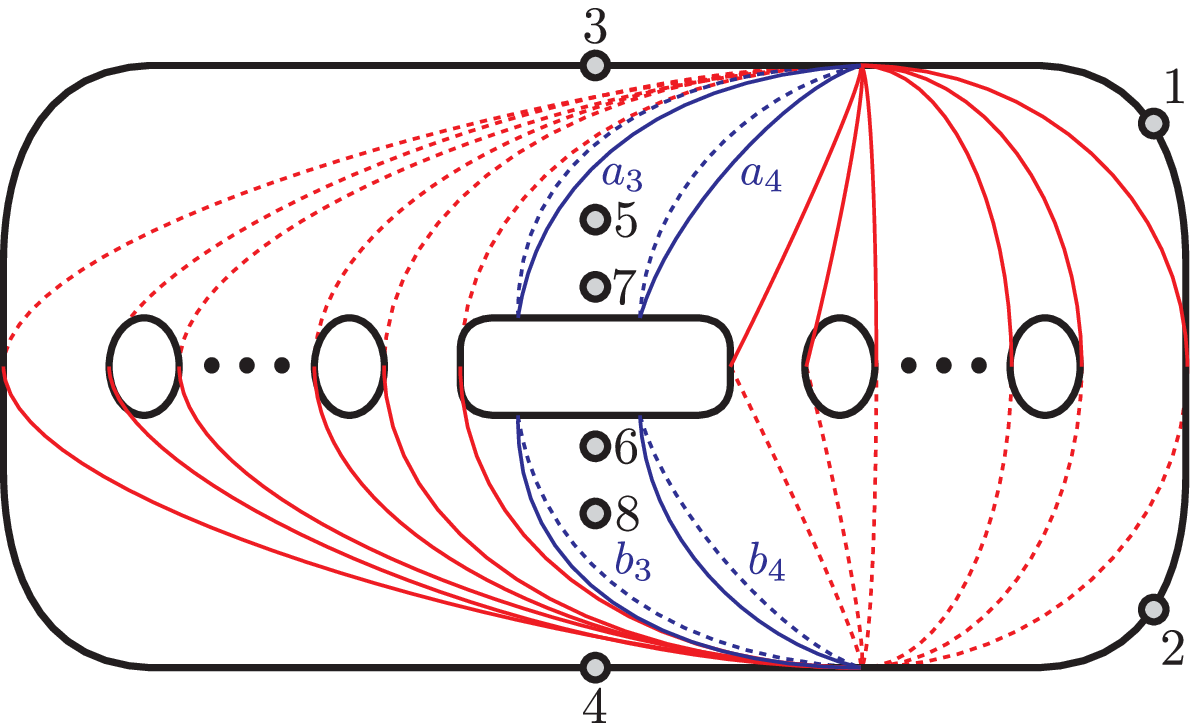}} 
	\vspace{-\baselineskip}
	\caption{The curves for $W_{\mathrm{\rmI}B\sharp}$ of odd $g$.} \label{F:WIBsharpodd}	
\end{figure}
\begin{figure}[htbp]
	\centering
	\subfigure[$B_{0,1}, B_{1,1}, \cdots, B_{g,1}, a_1, a_2, b_1, b_2$.]
	{\includegraphics[height=110pt]{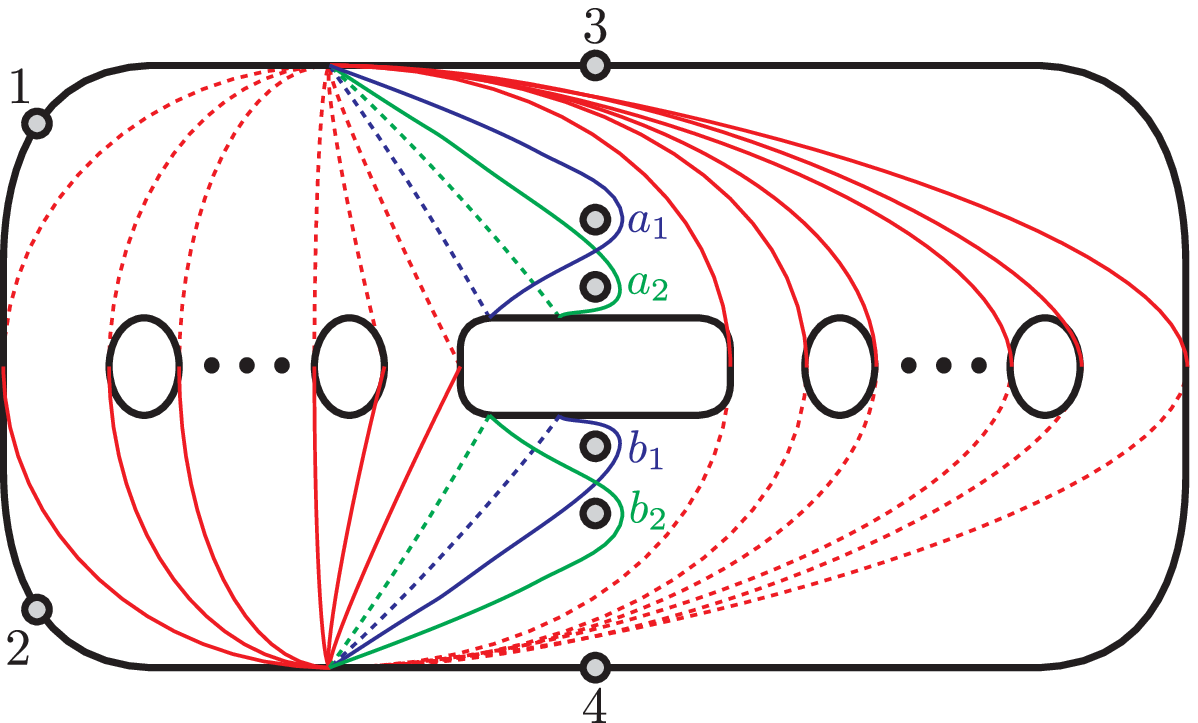}} 
	\hspace{.8em}	
	\subfigure[$B_{0,2}, B_{1,2}, \cdots, B_{g,2}, a_3, a_4, b_3, b_4$.]
	{\includegraphics[height=110pt]{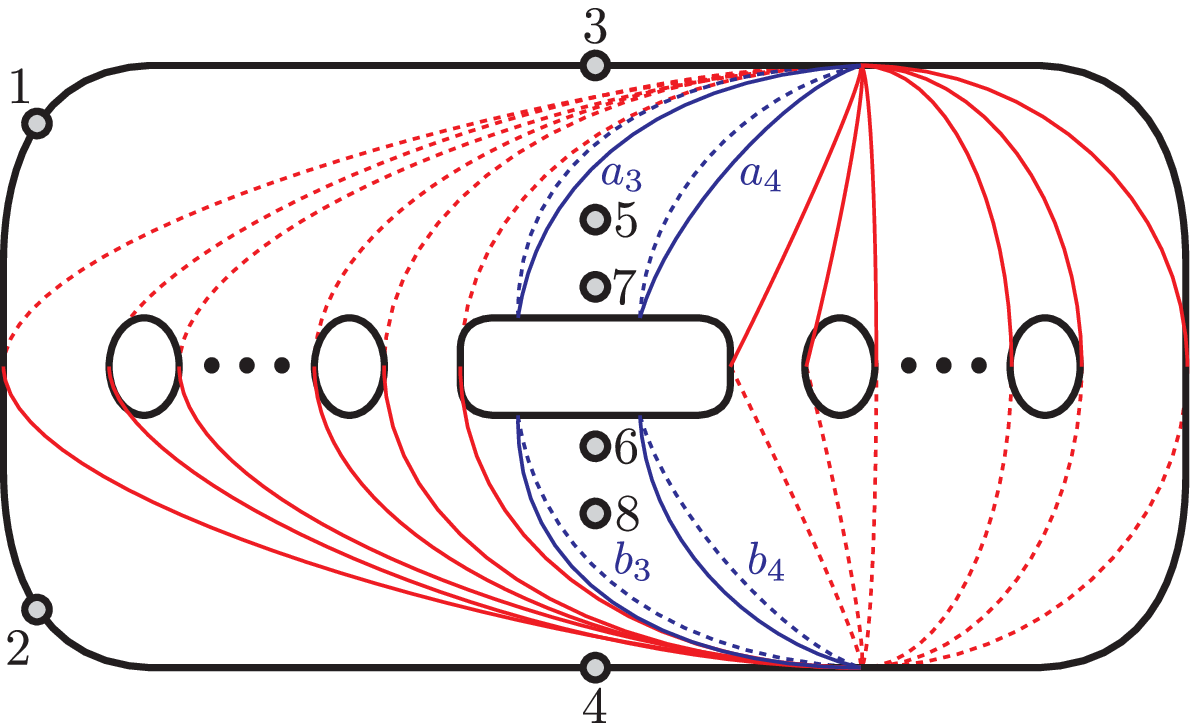}} 
	\vspace{-\baselineskip}
	\caption{The curves for $W_{\mathrm{\rmI}B\flat}$ of odd $g$.} \label{F:WIBflatodd}	
\end{figure}

The relations that come from the other lantern breedings will be constructed in much the same way as~\eqref{eq:WIAodd} so we will shorten the calculations.
We will also use the same symbols as those in~\eqref{eq:WIAodd} for the simplicity of notation.

\subsubsection{Type $\mathrm{\rmI}B\sharp$}
Perform the lantern breedings to $W_{\mathrm{\rmI}} = t_{\delta_1} \cdots t_{\delta_6}$ using $\{ a_1, x_1, \delta_5 \}$ and $\{ b_1, y_1, \delta_6 \}$:
\begin{align}
	t_{\delta_1} \cdots t_{\delta_6} 
	&= t_{B_0} \cdots t_{B_{g}} t_{a_1} t_{a_3} t_{b_1} t_{b_3} \cdot t_{B_0} \cdots t_{B_{g}} t_{x_1} t_{x_2} t_{y_1} t_{y_2} \notag \\
	&= t_{B_0} \cdots t_{B_{g}} \underline{t_{a_1} t_{x_1}} \; \underline{t_{b_1} t_{y_1}} \cdot {}_{t_{a_3}t_{b_3} t_{x_1}^{-1}t_{y_1}^{-1}}(t_{B_0} \cdots t_{B_{g}}) \;  t_{x_2} t_{a_3} t_{y_2} t_{b_3} \notag \\
	&\downarrow \notag \\
	t_{\delta_1} \cdots t_{\delta_8} 
	&= t_{B_{0,1}} t_{B_{1,1}} \cdots t_{B_{g,1}} t_{a_1} t_{a_2} t_{b_1} t_{b_2} \cdot 
	t_{B_{0,2}} t_{B_{1,2}} \cdots t_{B_{g,2}} t_{a_3} t_{a_4} t_{b_3} t_{b_4}, \label{eq:WIBsharpodd}
\end{align}
where the resulting curves are shown in Figure~\ref{F:WIBsharpodd}.
We denote by $W_{\mathrm{\rmI}B\sharp}$ the last word.
\qed

\subsubsection{Type $\mathrm{\rmI}B\flat$}
Perform the lantern breedings to $W_{\mathrm{\rmI}} = t_{\delta_1} \cdots t_{\delta_6}$ using $\{ a_3, x_2, \delta_5 \}$ and $\{ b_3, y_2, \delta_6 \}$:
\begin{align}
	t_{\delta_1} \cdots t_{\delta_6} 
	&= t_{B_0} \cdots t_{B_{2h-1}} t_{a_1} t_{a_3} t_{b_1} t_{b_3} \cdot t_{B_0} \cdots t_{B_{2h-1}} t_{x_1} t_{x_2} t_{y_1} t_{y_2} \notag \\
	&= {}_{t_{a_1}^{-1}t_{b_1}^{-1} t_{x_2}t_{y_2}}(t_{B_0} \cdots t_{B_{g}}) \; \underline{t_{x_2} t_{a_3}} \; \underline{t_{y_2} t_{b_3}} \cdot  t_{B_0} \cdots t_{B_{g}} t_{x_1} t_{a_1} t_{y_1} t_{b_1}  \notag \\
	&\downarrow \notag \\
	t_{\delta_1} \cdots t_{\delta_8} 
	&= t_{B_{0,1}} t_{B_{1,1}} \cdots t_{B_{g,1}} t_{a_1} t_{a_2} t_{b_1} t_{b_2} \cdot 
	t_{B_{0,2}} t_{B_{1,2}} \cdots t_{B_{g,2}} t_{a_3} t_{a_4} t_{b_3} t_{b_4}, \label{eq:WIBflatodd}
\end{align}
where the resulting curves are shown in Figure~\ref{F:WIBflatodd} (the figure has been rotated by $\pi$ with respect to the central vertical axis so that we may better compare it with the figures for the other lifts).
We denote by $W_{\mathrm{\rmI}B\flat}$ the last word.
\qed

\subsubsection{Type $\mathrm{\rmII}A$}
\begin{figure}[htbp]
	\centering
	\subfigure[$B_{0,1}, B_{1,1}, \cdots, B_{g,1}, a_1, a_2, b_1, b_2$.]
	{\includegraphics[height=110pt]{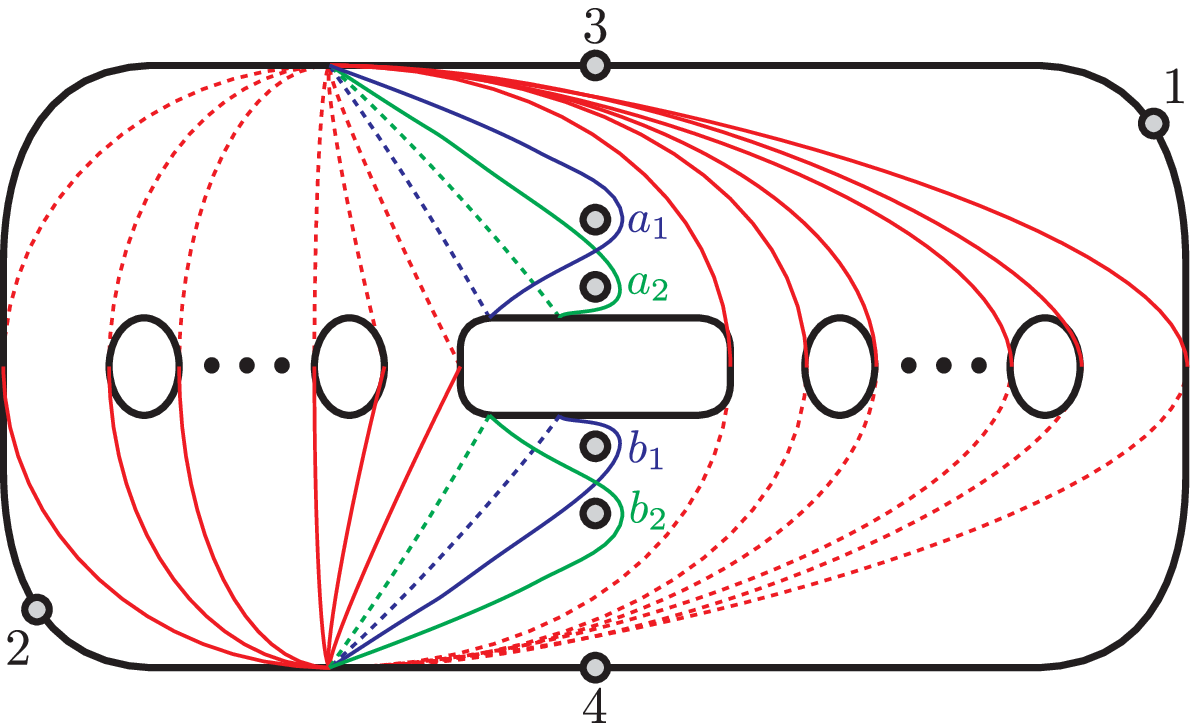}} 
	\hspace{.8em}	
	\subfigure[$B_{0,2}, B_{1,2}, \cdots, B_{g,2}, a_3, a_4, b_3, b_4$.]
	{\includegraphics[height=110pt]{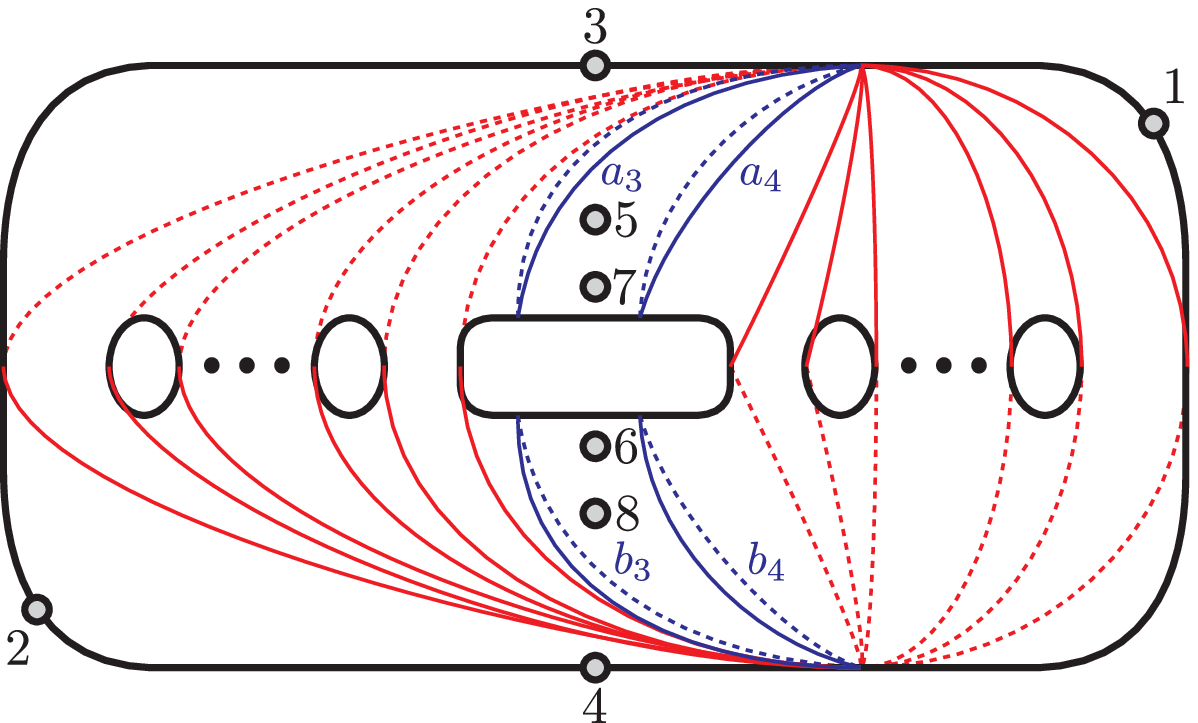}} 
	\vspace{-\baselineskip}
	\caption{The curves for $W_{\mathrm{\rmII}A}$ of odd $g$.} \label{F:WIIAodd}	
\end{figure}

Perform the lantern breedings to $W_{\mathrm{\rmII}} = t_{\delta_1} \cdots t_{\delta_6}$ using $\{ a_1, x_1, \delta_5 \}$ and $\{ b_3, y_2, \delta_6 \}$:
\begin{align}
	t_{\delta_1} \cdots t_{\delta_6} 
	&= t_{B_0^{\prime}} \cdots t_{B_{g}^{\prime}} t_{a_1} t_{a_3} t_{b_2} t_{b_3} \cdot t_{B_0^{\prime}} \cdots t_{B_{g}^{\prime}} t_{x_1} t_{x_2} t_{y_1} t_{y_2} \notag \\
	&= {}_{t_{b_2}^{-1} t_{y_2}}(t_{B_0^{\prime}} \cdots t_{B_{g}^{\prime}}) \; \underline{t_{a_1} t_{x_1}} \; \underline{t_{y_2} t_{b_3}} \cdot {}_{t_{a_3} t_{x_1}^{-1}}(t_{B_0^{\prime}} \cdots t_{B_{g}^{\prime}}) \; t_{x_2} t_{a_3} t_{y_1} t_{b_2}  \notag \\
	&\downarrow \notag \\
	t_{\delta_1} \cdots t_{\delta_8} 
	&= t_{B_{0,1}} t_{B_{1,1}} \cdots t_{B_{g,1}} t_{a_1} t_{a_2} t_{b_1} t_{b_2} \cdot 
	t_{B_{0,2}} t_{B_{1,2}} \cdots t_{B_{g,2}} t_{a_3} t_{a_4} t_{b_3} t_{b_4}, \label{eq:WIIAodd}
\end{align}
where the resulting curves are shown in Figure~\ref{F:WIIAodd}.
We denote by $W_{\mathrm{\rmII}A}$ the last word.
\qed

\subsubsection{Type $\mathrm{\rmII}B\sharp$}
\begin{figure}[htbp]
	\centering
	\subfigure[$B_{0,1}, B_{1,1}, \cdots, B_{g,1}, a_1, a_2, b_1, b_2$.]
	{\includegraphics[height=110pt]{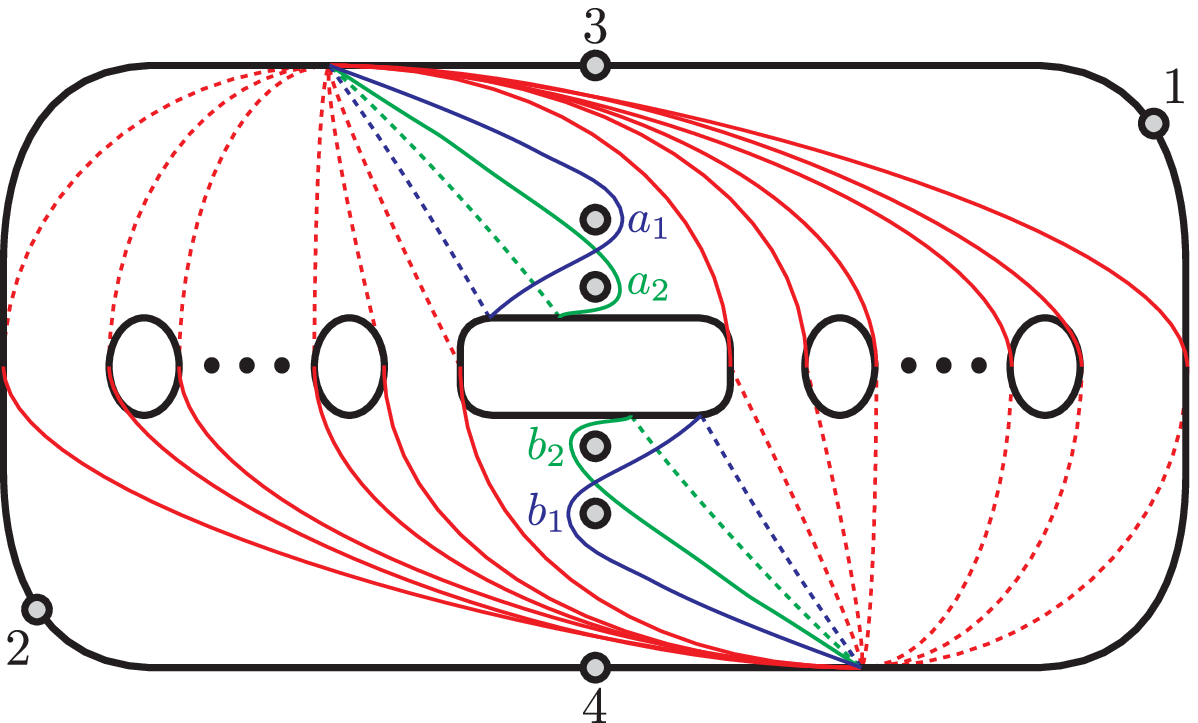}} 
	\hspace{.8em}	
	\subfigure[$B_{0,2}, B_{1,2}, \cdots, B_{g,2}, a_3, a_4, b_3, b_4$.]
	{\includegraphics[height=110pt]{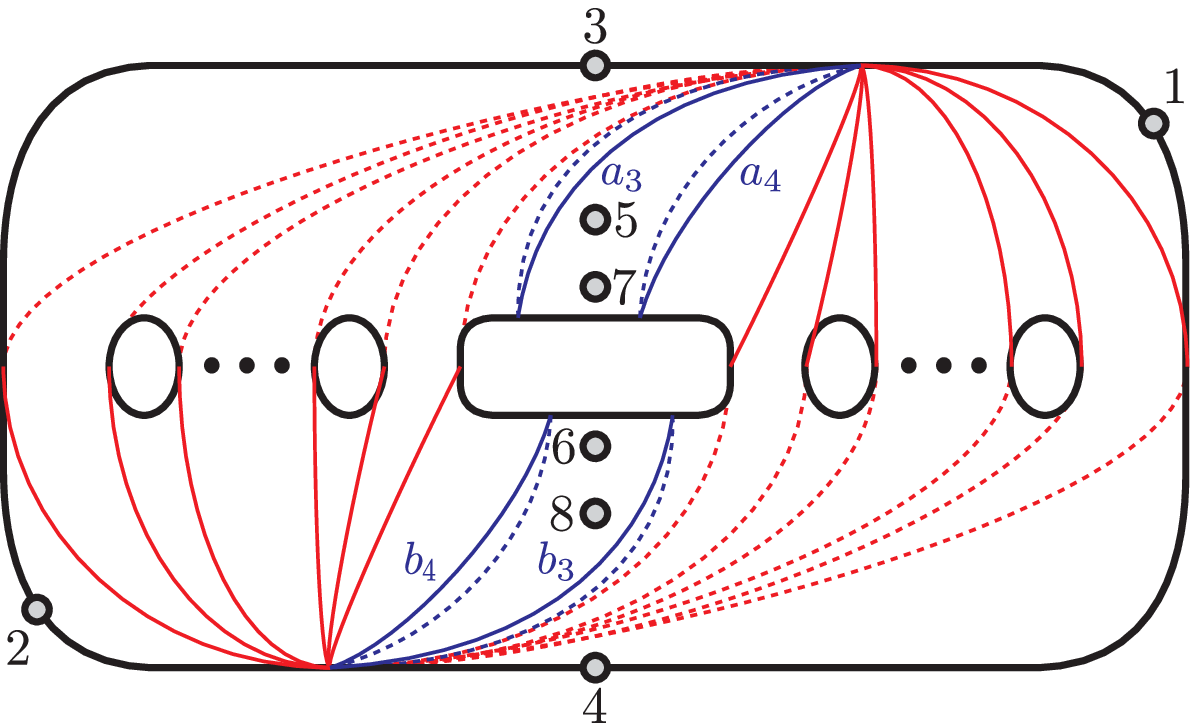}} 
	\vspace{-\baselineskip}
	\caption{The curves for $W_{\mathrm{\rmII}B\sharp}$ of odd $g$.} \label{F:WIIBsharpodd}	
\end{figure}

Perform the lantern breedings to $W_{\mathrm{\rmII}} = t_{\delta_1} \cdots t_{\delta_6}$ using $\{ a_1, x_1, \delta_5 \}$ and $\{ b_2, y_1, \delta_6 \}$:
\begin{align}
	t_{\delta_1} \cdots t_{\delta_6} 
	&= t_{B_0^{\prime}} \cdots t_{B_{g}^{\prime}} t_{a_1} t_{a_3} t_{b_2} t_{b_3} \cdot t_{B_0^{\prime}} \cdots t_{B_{g}^{\prime}} t_{x_1} t_{x_2} t_{y_1} t_{y_2} \notag \\
	&= t_{B_0^{\prime}} \cdots t_{B_{g}^{\prime}} \underline{t_{a_1} t_{x_1}} \; \underline{t_{b_2} t_{y_1}} \cdot {}_{t_{a_3} t_{b_3} t_{x_1}^{-1} t_{y_1}^{-1} }(t_{B_0^{\prime}} \cdots t_{B_{g}^{\prime}}) \; t_{x_2} t_{a_3} t_{y_2} t_{b_3}  \notag \\
	&\downarrow \notag \\
	t_{\delta_1} \cdots t_{\delta_8} 
	&= t_{B_{0,1}} t_{B_{1,1}} \cdots t_{B_{g,1}} t_{a_1} t_{a_2} t_{b_1} t_{b_2} \cdot 
	t_{B_{0,2}} t_{B_{1,2}} \cdots t_{B_{g,2}} t_{a_3} t_{a_4} t_{b_3} t_{b_4}, \label{eq:WIIBsharpodd}
\end{align}
where the resulting curves are shown in Figure~\ref{F:WIIBsharpodd}.
We denote by $W_{\mathrm{\rmII}B\sharp}$ the last word.
\qed

\subsubsection{Type $\mathrm{\rmII}B\flat$}
\begin{figure}[htbp]
	\centering
	\subfigure[$B_{0,1}, B_{1,1}, \cdots, B_{g,1}, a_1, a_2, b_1, b_2$.]
	{\includegraphics[height=110pt]{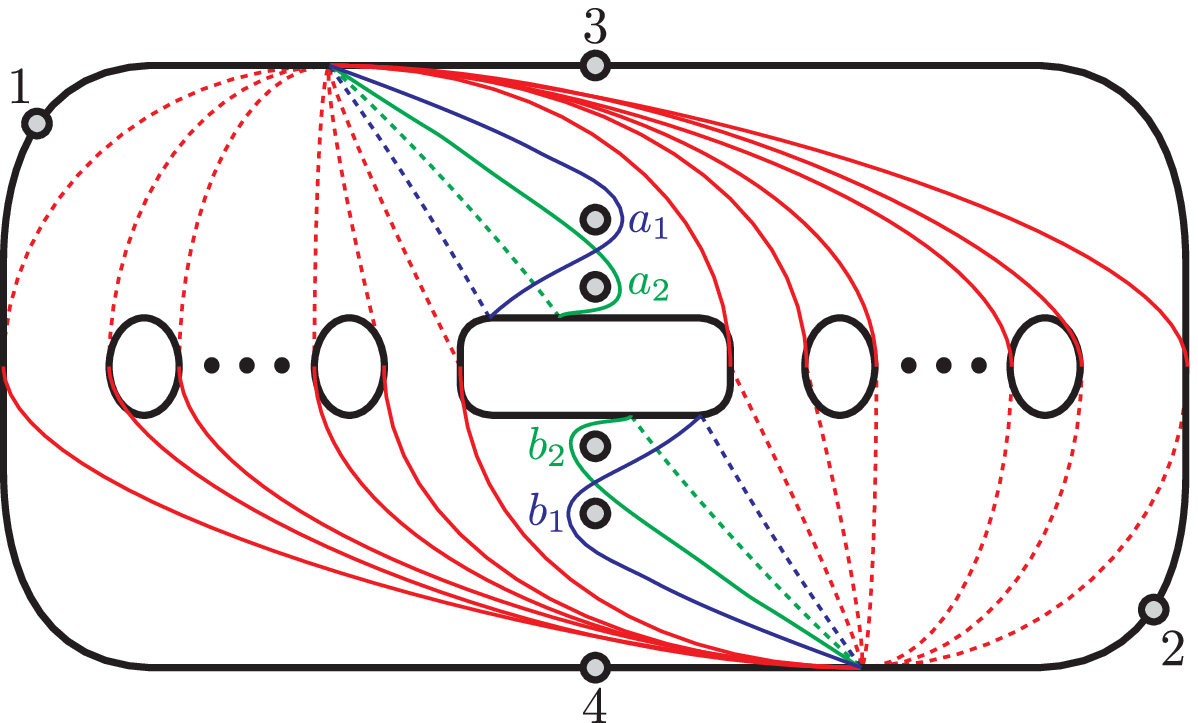}} 
	\hspace{.8em}	
	\subfigure[$B_{0,2}, B_{1,2}, \cdots, B_{g,2}, a_3, a_4, b_3, b_4$.]
	{\includegraphics[height=110pt]{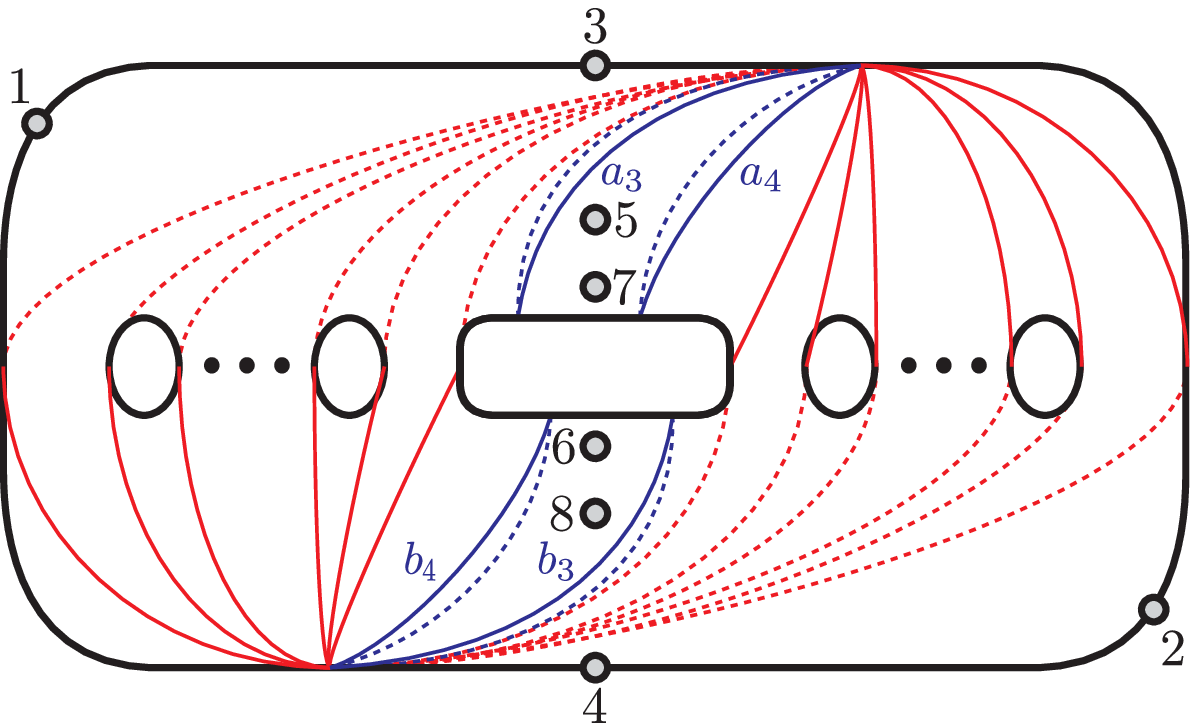}} 
	\vspace{-\baselineskip}
	\caption{The curves for $W_{\mathrm{\rmII}B\flat}$ of odd $g$.} \label{F:WIIBflatodd}	
\end{figure}

Perform the lantern breedings to $W_{\mathrm{\rmII}} = t_{\delta_1} \cdots t_{\delta_6}$ using $\{ a_3, x_2, \delta_5 \}$ and $\{ b_3, y_2, \delta_6 \}$:
\begin{align}
	t_{\delta_1} \cdots t_{\delta_6} 
	&= t_{B_0^{\prime}} \cdots t_{B_{g}^{\prime}} t_{a_1} t_{a_3} t_{b_2} t_{b_3} \cdot t_{B_0^{\prime}} \cdots t_{B_{g}^{\prime}} t_{x_1} t_{x_2} t_{y_1} t_{y_2} \notag \\
	&= {}_{t_{a_1}^{-1} t_{b_2}^{-1} t_{x_2} t_{y_2}}(t_{B_0^{\prime}} \cdots t_{B_{g}^{\prime}}) \; \underline{t_{x_2} t_{a_3}} \; \underline{t_{y_2} t_{b_3}} \cdot t_{B_0^{\prime}} \cdots t_{B_{g}^{\prime}} t_{x_1} t_{a_1} t_{y_1} t_{b_2}  \notag \\
	&\downarrow \notag \\
	t_{\delta_1} \cdots t_{\delta_8} 
	&= t_{B_{0,1}} t_{B_{1,1}} \cdots t_{B_{g,1}} t_{a_1} t_{a_2} t_{b_1} t_{b_2} \cdot 
	t_{B_{0,2}} t_{B_{1,2}} \cdots t_{B_{g,2}} t_{a_3} t_{a_4} t_{b_3} t_{b_4}, \label{eq:WIIBflatodd}
\end{align}
where the resulting curves are shown in Figure~\ref{F:WIIBflatodd}  (the figure has been rotated by $\pi$ with respect to the central vertical axis).
We denote by $W_{\mathrm{\rmII}B\flat}$ the last word.
\qed


\subsubsection*{Summary}
Each of the six relations we obtained here, namely, $W_{\mathrm{\rmI}A}=t_{\delta_1} \cdots t_{\delta_8}$, 
$W_{\mathrm{\rmI}B\sharp}=t_{\delta_1} \cdots t_{\delta_8}$, $W_{\mathrm{\rmI}B\flat}=t_{\delta_1} \cdots t_{\delta_8}$, $W_{\mathrm{\rmII}A}=t_{\delta_1} \cdots t_{\delta_8}$, $W_{\mathrm{\rmII}B\sharp}=t_{\delta_1} \cdots t_{\delta_8}$ and $W_{\mathrm{\rmII}B\flat}=t_{\delta_1} \cdots t_{\delta_8}$, is obviously a lift of the MCK relation $W=1$ for odd $g$.
Furthermore, we can observe that three of them, that is, $W_{\mathrm{\rmI}A}=t_{\delta_1} \cdots t_{\delta_8}$, 
$W_{\mathrm{\rmI}B\sharp}=t_{\delta_1} \cdots t_{\delta_8}$, $W_{\mathrm{\rmI}B\flat}=t_{\delta_1} \cdots t_{\delta_8}$ are indeed further lifts of the known lift of the MCK relation $W=t_{\delta_1} t_{\delta_2}$; just cap off the boundary components other than $\delta_1$ and $\delta_2$ (and rotate back the figure for $W_{\mathrm{\rmI}B\flat}=t_{\delta_1} \cdots t_{\delta_8}$).

\subsection{The case of even genus: $g=2h$} \label{subsect:evengenus}
\begin{figure}[h]
	\centering
	\subfigure[Embedding for $W_{\mathrm{\rmI}}$. \label{F:Embedding_WI}]
	{\includegraphics[height=100pt]{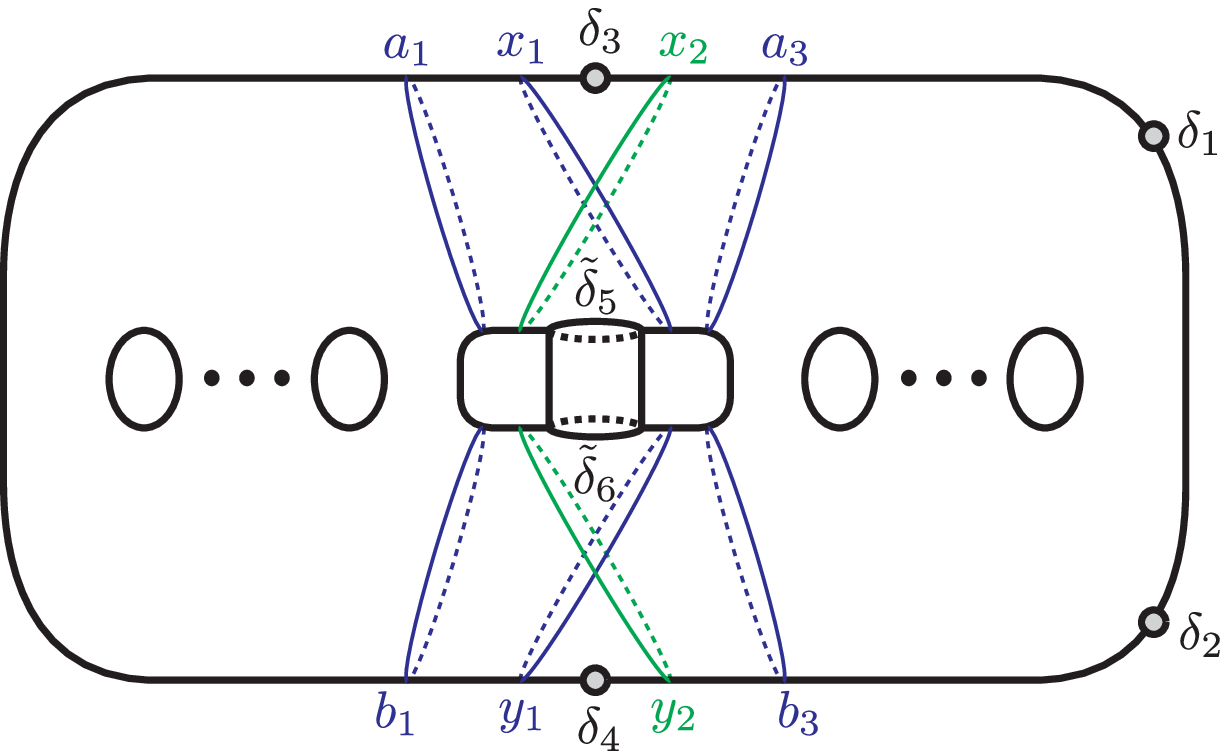}} 
	\hspace{.8em}	
	\subfigure[Embedding for $W_{\mathrm{\rmII}}$. \label{F:Embedding_WII}]
	{\includegraphics[height=100pt]{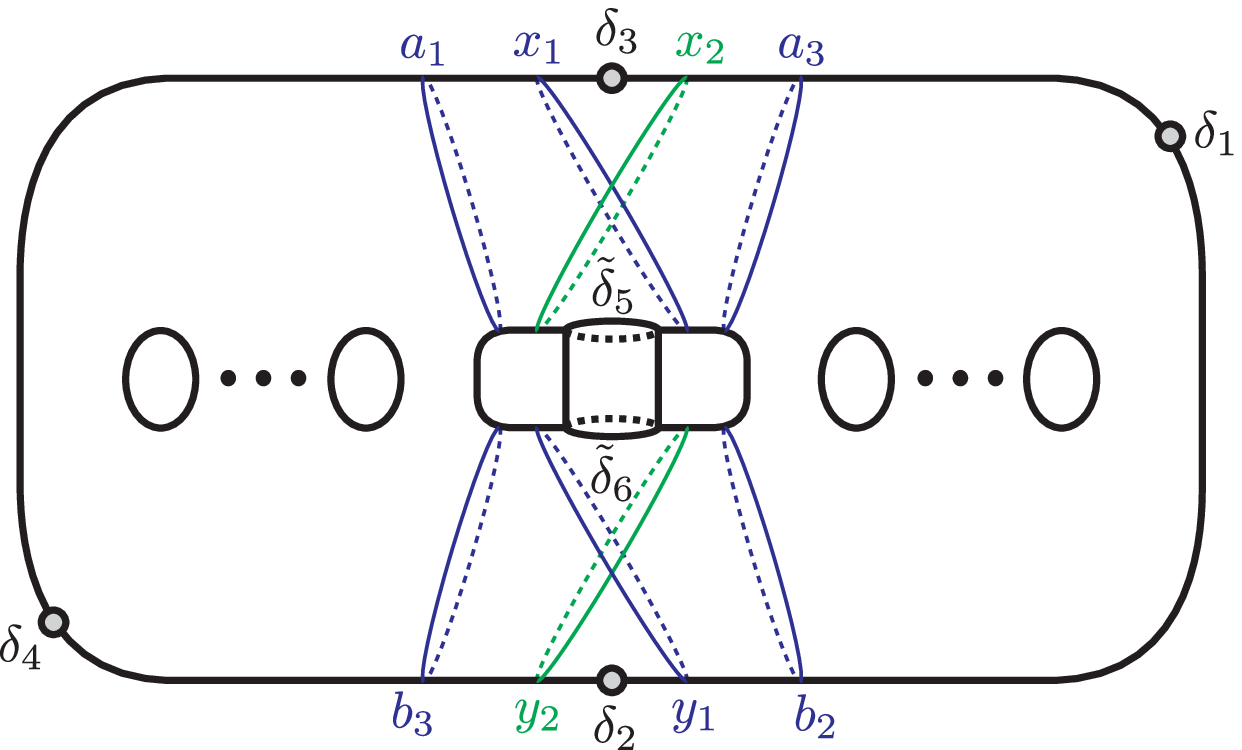}} 
	\caption{Embeddings of $\Sigma_{g-1}^6$ into $\Sigma_g^4$.} \label{F:Embedding_WIWII}	
\end{figure}

We next construct lifts of the MCK relation to $\Sigma_g^4$ for even $g=2h$.
In order to utilize the relations $W_{\mathrm{\rmI}} = t_{\delta_1} \cdots t_{\delta_6}$ and $W_{\mathrm{\rmII}} = t_{\delta_1} \cdots t_{\delta_6}$ (which hold in $\Mod(\Sigma_{2h-1}^6) = \Mod(\Sigma_{g-1}^6)$), we embed $\Sigma_{g-1}^6$ into $\Sigma_g^4$ by connecting the boundary components $\delta_5$ and $\delta_6$ of $\Sigma_{g-1}^6$ in Figure~\ref{F:WI} or Figure~\ref{F:WII} with an obvious tube as indicated in Figure~\ref{F:Embedding_WIWII}.
Then the relations become $W_{\mathrm{\rmI}} = t_{\delta_1} \cdots t_{\delta_4} t_{\tilde{\delta}_5} t_{\tilde\delta_6}$ and $W_{\mathrm{\rmII}} = t_{\delta_1} \cdots t_{\delta_4} t_{\tilde{\delta}_5} t_{\tilde\delta_6}$ (where $t_{\tilde{\delta}_5} = t_{\tilde\delta_6}$).

\subsubsection{Type $\mathrm{\rmI}A$}
\begin{figure}[h]
	\centering
	\subfigure
	{\includegraphics[height=100pt]{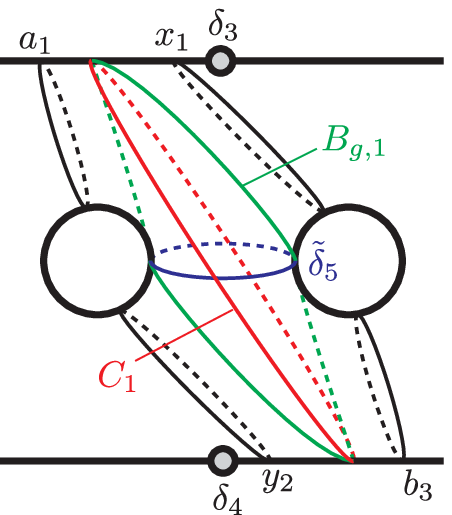}} 
	\hspace{1.8em}	
	\subfigure
	{\includegraphics[height=100pt]{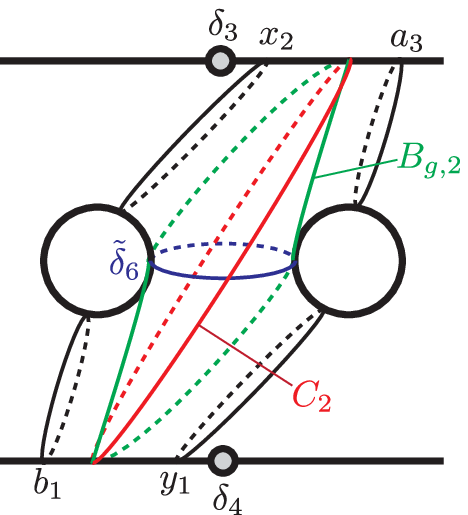}} 
	\caption{Lantern relations.} \label{F:Lantern_WIAeven}	
\end{figure}

\begin{figure}[!p] 
	\centering
	\subfigure[$B_{0,1}, B_{1,1}, \cdots, B_{g,1}, C_1$.]
	{\includegraphics[height=100pt]{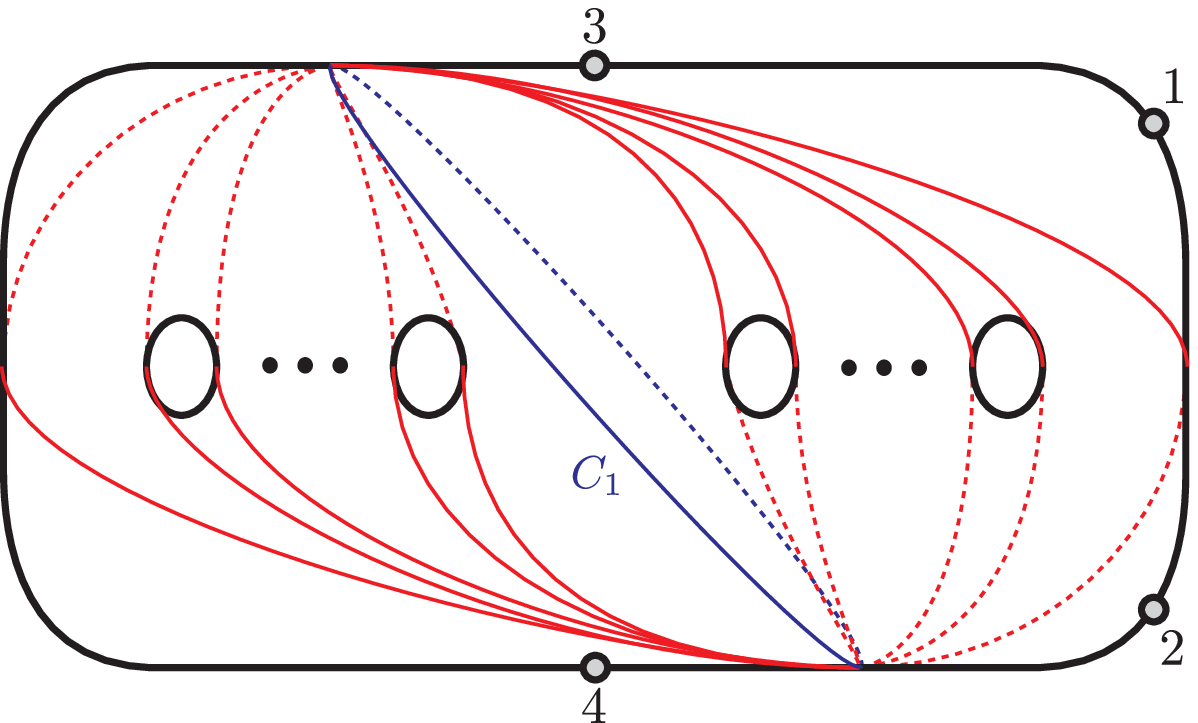}} 
	\hspace{.8em}	
	\subfigure[$B_{0,2}, B_{1,2}, \cdots, B_{g,2}, C_2$.]
	{\includegraphics[height=100pt]{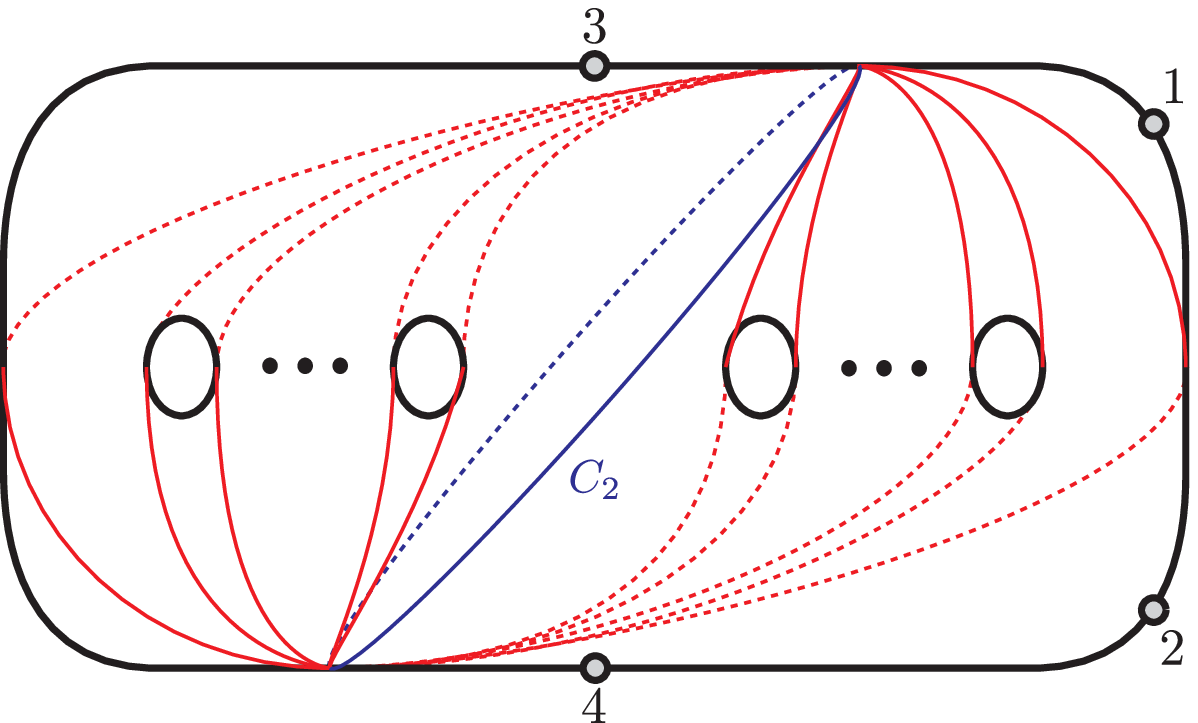}} 
	\vspace{-\baselineskip}
	\caption{The curves for $W_{\mathrm{\rmI}A}$ of even $g$.} \label{F:WIAeven}
\end{figure}

Consider the relation $W_{\mathrm{\rmI}} = t_{\delta_1} \cdots t_{\delta_4} t_{\tilde{\delta}_5} t_{\tilde\delta_6}$ in $\Mod(\Sigma_g^4)$.
We can find a four-holed sphere bounded by $\{ a_1, x_1, b_3, y_2 \}$,
and hence we have a lantern relation $t_{a_1} t_{x_1} t_{b_3} t_{y_2} = t_{B_{g,1}} t_{C_1} t_{\tilde\delta_5}$ with the curves in the left of Figure~\ref{F:Lantern_WIAeven}.
We also find another lantern relation with $\{ a_3, x_2, b_1, y_1 \}$ as in the right of Figure~\ref{F:Lantern_WIAeven}.
\\
Applying those lantern relations we alter $W_{\mathrm{\rmI}} = t_{\delta_1} \cdots t_{\delta_4} t_{\tilde{\delta}_5} t_{\tilde\delta_6}$ as
\begin{align}
	t_{\delta_1} \cdots t_{\delta_4} 
	&= \underline{t_{B_0} \cdots t_{B_{g-1}}} t_{a_1} t_{a_3} t_{b_1} t_{b_3} t_{\tilde{\delta}_5}^{-1} \cdot \underline{t_{B_0} \cdots t_{B_{g-1}} t_{x_1}} t_{x_2} t_{y_1} \underline{t_{y_2} t_{\tilde{\delta}_6}^{-1}}  \notag \\
	&= {}_{t_{y_2}}(t_{B_0} \cdots t_{B_{g-1}}) \; t_{y_2} t_{a_1} t_{a_3} t_{b_1} t_{b_3} t_{\tilde{\delta}_5}^{-1} t_{x_1} \cdot {}_{t_{x_1}^{-1}}(t_{B_0} \cdots t_{B_{g-1}}) \; t_{x_2} t_{y_1} t_{\tilde{\delta}_6}^{-1} \notag \\
	&= \underline{{}_{t_{y_2}}(t_{B_0} \cdots t_{B_{g-1}}) \; t_{b_1}} t_{a_1} t_{x_1} t_{b_3} t_{y_2} t_{\tilde{\delta}_5}^{-1} \underline{t_{a_3} \cdot {}_{t_{x_1}^{-1}}(t_{B_0} \cdots t_{B_{g-1}})} \; t_{x_2} t_{y_1} t_{\tilde{\delta}_6}^{-1} \notag \\
	&= {}_{t_{b_1}^{-1} t_{y_2}}(t_{B_0} \cdots t_{B_{g-1}}) \;  t_{a_1} t_{x_1} t_{b_3} t_{y_2}  t_{\tilde{\delta}_5}^{-1} \cdot {}_{t_{a_3} t_{x_1}^{-1}}(t_{B_0} \cdots t_{B_{g-1}}) \; t_{a_3} t_{x_2} t_{y_1} t_{\tilde{\delta}_6}^{-1} t_{b_1} \notag \\
	&= {}_{t_{b_1}^{-1} t_{y_2}}(t_{B_0} \cdots t_{B_{g-1}}) \; \underline{t_{a_1} t_{x_1} t_{b_3} t_{y_2}}  t_{\tilde{\delta}_5}^{-1} \cdot {}_{t_{a_3} t_{x_1}^{-1}}(t_{B_0} \cdots t_{B_{g-1}}) \; \underline{t_{a_3} t_{x_2} t_{b_1} t_{y_1}} t_{\tilde{\delta}_6}^{-1} \notag \\
	&= {}_{t_{b_1}^{-1} t_{y_2}}(t_{B_0} \cdots t_{B_{g-1}}) \; t_{B_{g,1}} t_{C_1} t_{\tilde\delta_5} t_{\tilde{\delta}_5}^{-1} \cdot {}_{t_{a_3} t_{x_1}^{-1}}(t_{B_0} \cdots t_{B_{g-1}}) \;  t_{B_{g,2}} t_{C_2} t_{\tilde\delta_6} t_{\tilde{\delta}_6}^{-1} \notag \\
	&= {}_{t_{b_1}^{-1} t_{y_2}}(t_{B_0} \cdots t_{B_{g-1}}) \; t_{B_{g,1}} t_{C_1} \cdot {}_{t_{a_3} t_{x_1}^{-1}}(t_{B_0} \cdots t_{B_{g-1}}) \;  t_{B_{g,2}} t_{C_2}. \notag 
\end{align}	
By renaming the resulting curves, which are depicted in Figure~\ref{F:WIAeven}, we obtain 
\begin{align}
	t_{\delta_1} \cdots t_{\delta_4} 
	&= t_{B_{0,1}} t_{B_{1,1}} \cdots t_{B_{g,1}} t_{C_1} \cdot 
	t_{B_{0,2}} t_{B_{1,2}} \cdots t_{B_{g,2}} t_{C_2}. \label{eq:WIAeven}
\end{align}
Let $W_{\mathrm{\rmI}A}$ denote the last word.
\qed

\afterpage{\clearpage} 
\begin{figure}[htbp]
	\centering
	\subfigure[$B_{0,1}, B_{1,1}, \cdots, B_{g,1}, C_1$.]
	{\includegraphics[height=100pt]{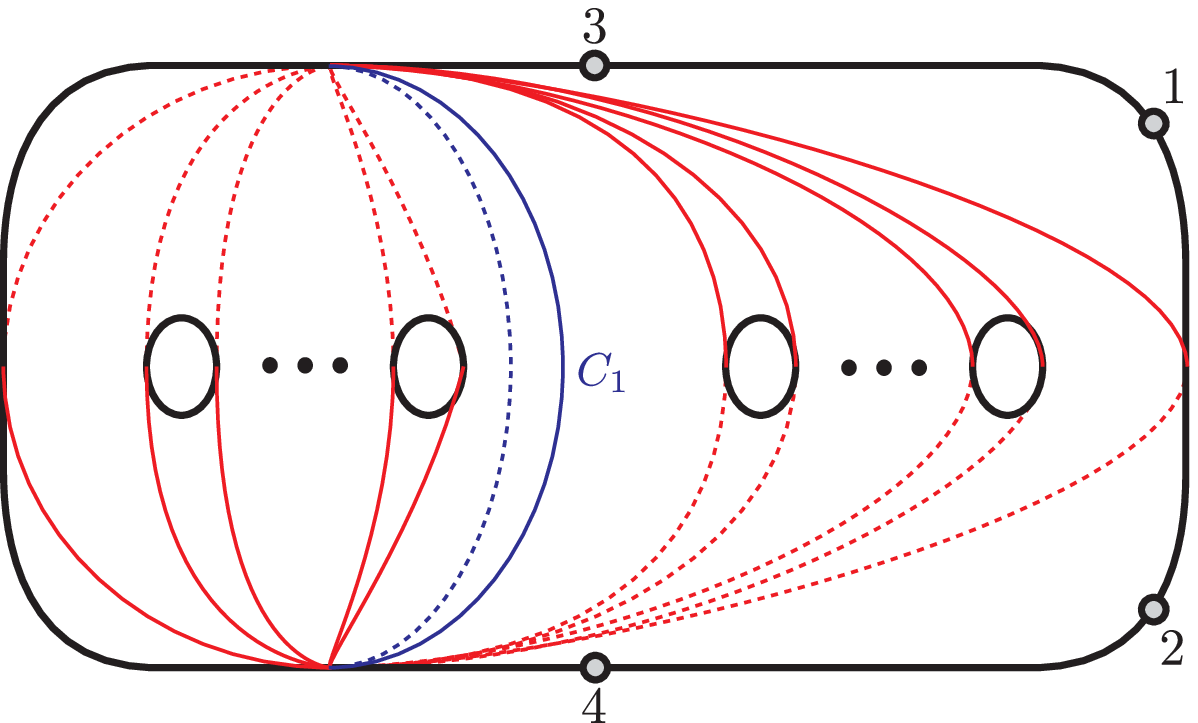}} 
	\hspace{.8em}	
	\subfigure[$B_{0,2}, B_{1,2}, \cdots, B_{g,2}, C_2$.]
	{\includegraphics[height=100pt]{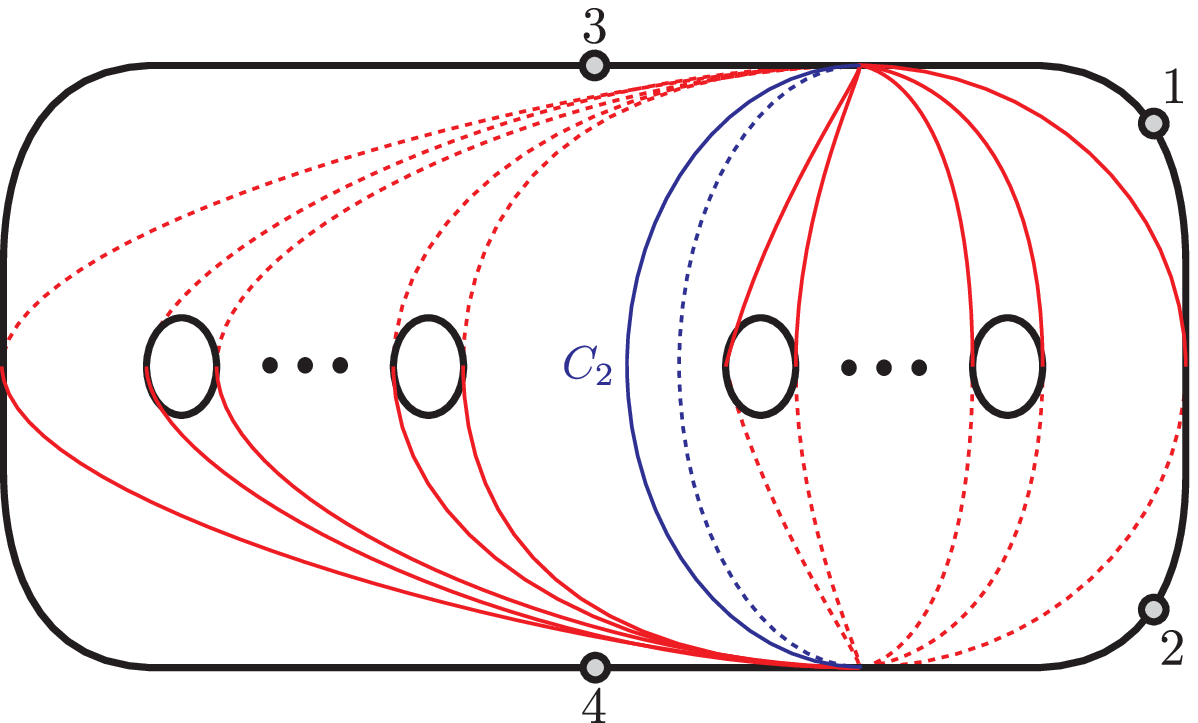}} 
	\vspace{-\baselineskip}
	\caption{The curves for $W_{\mathrm{\rmI}B}$ of even $g$.} \label{F:WIBeven}	
\end{figure}
\begin{figure}[htbp]
	\centering
	\subfigure[$B_{0,1}, B_{1,1}, \cdots, B_{g,1}, C_1$.]
	{\includegraphics[height=100pt]{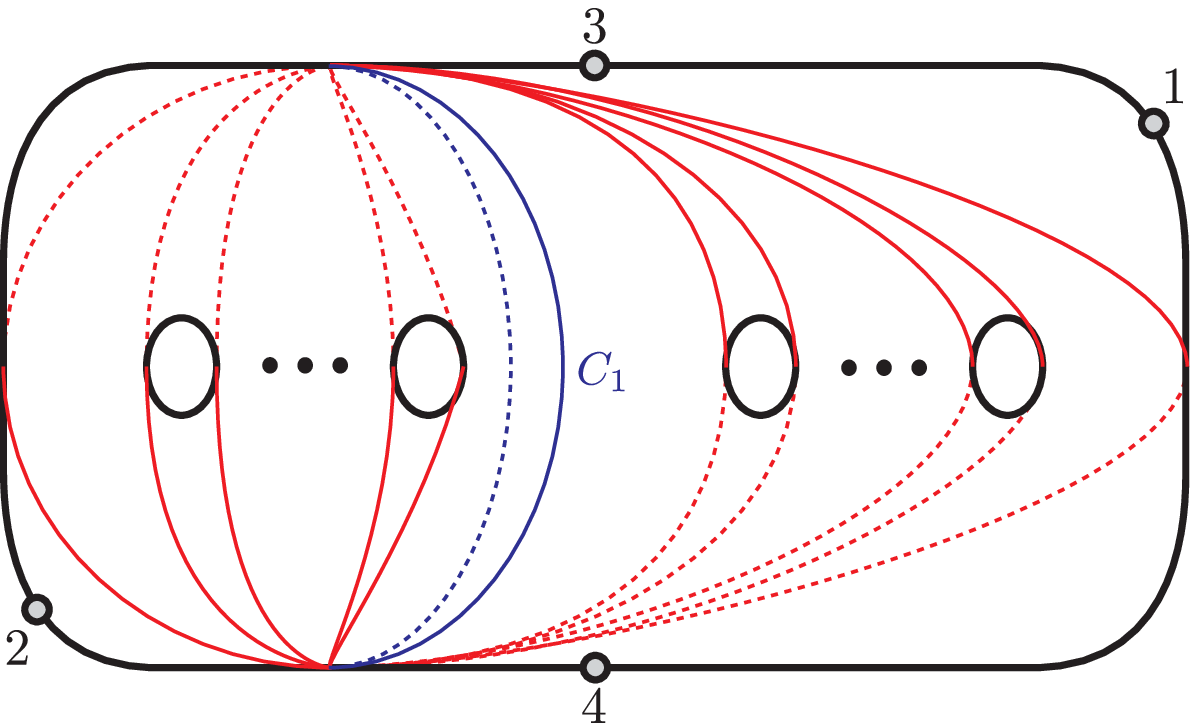}} 
	\hspace{.8em}	
	\subfigure[$B_{0,2}, B_{1,2}, \cdots, B_{g,2}, C_2$.]
	{\includegraphics[height=100pt]{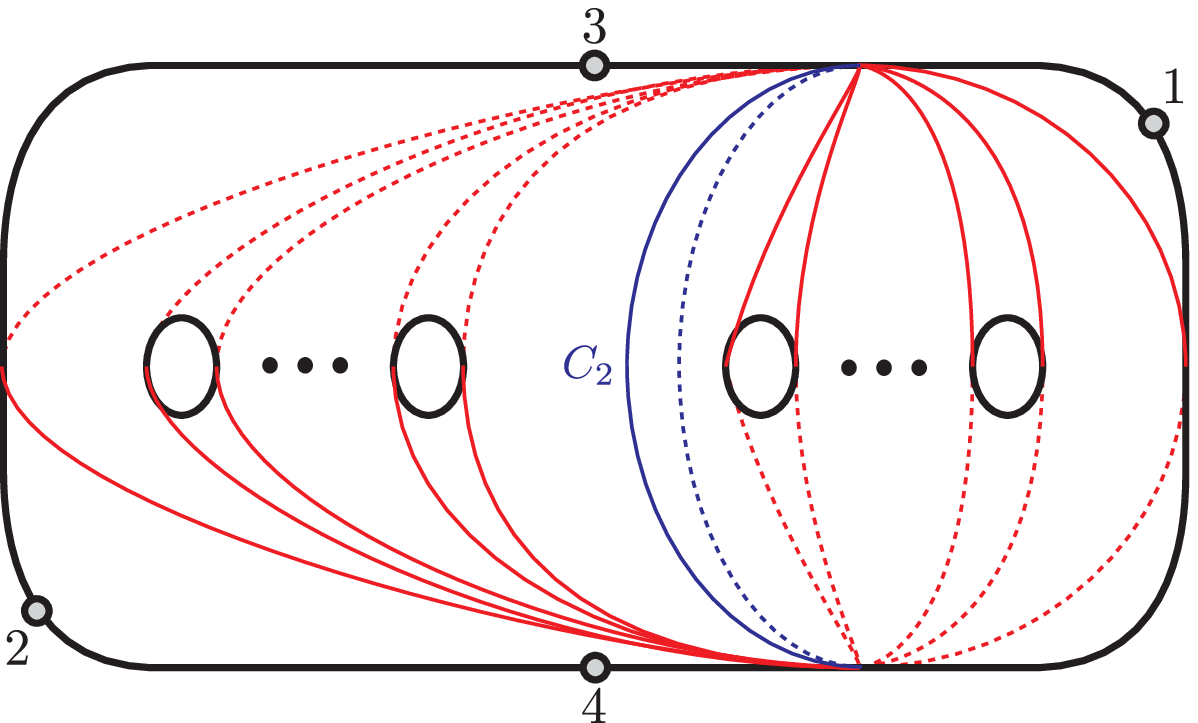}} 
	\vspace{-\baselineskip}
	\caption{The curves for $W_{\mathrm{\rmII}A}$ of even $g$.} \label{F:WIIAeven}	
\end{figure}
\begin{figure}[htbp]
	\centering
	\subfigure[$B_{0,1}, B_{1,1}, \cdots, B_{g,1}, C_1$.]
	{\includegraphics[height=100pt]{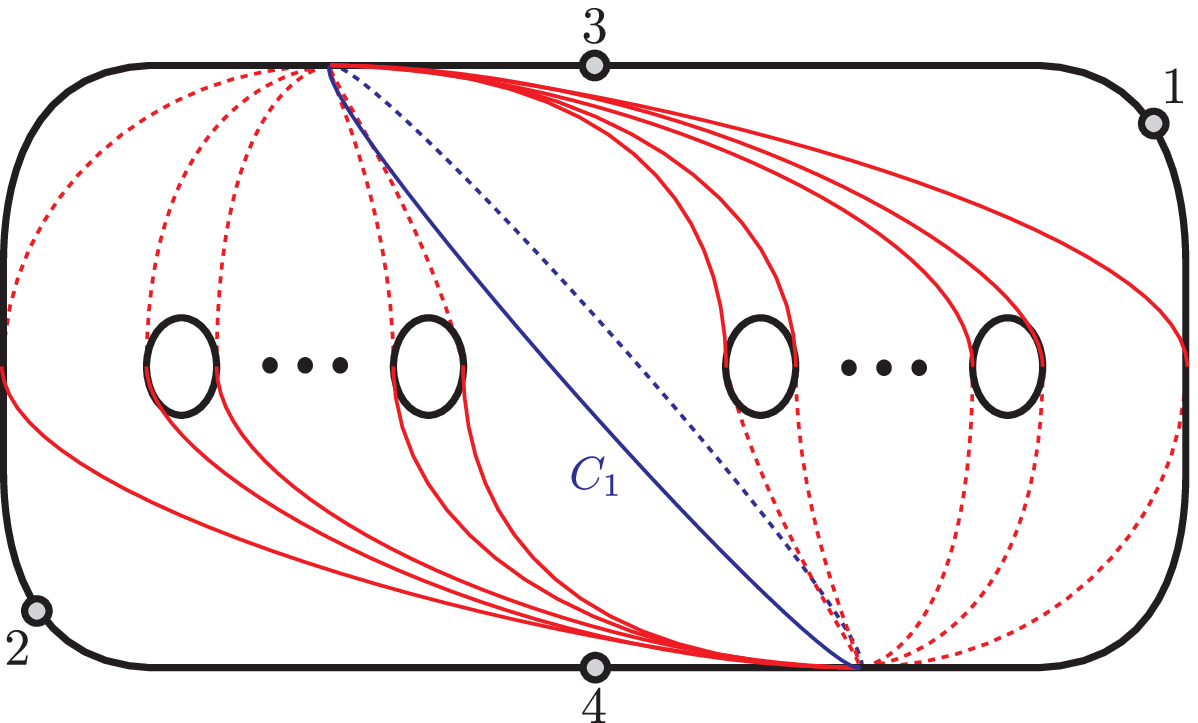}} 
	\hspace{.8em}	
	\subfigure[$B_{0,2}, B_{1,2}, \cdots, B_{g,2}, C_2$.]
	{\includegraphics[height=100pt]{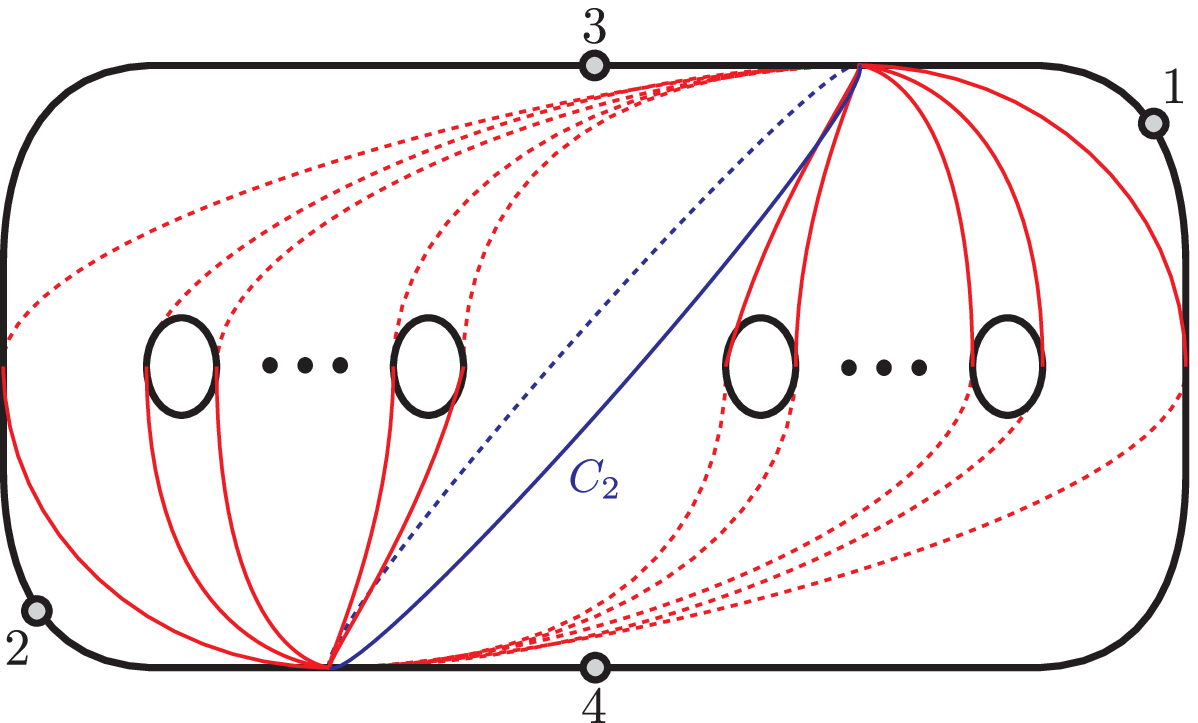}} 
	\vspace{-\baselineskip}
	\caption{The curves for $W_{\mathrm{\rmII}B}$ of even $g$.} \label{F:WIIBeven}	
\end{figure}

In what follows we further construct other new relations in the same way as~\eqref{eq:WIAeven}, so the procedures will be simplified.
We will again use the same symbols as those in~\eqref{eq:WIAeven}.

\subsubsection{Type $\mathrm{\rmI}B$}
Consider the relation $W_{\mathrm{\rmI}} = t_{\delta_1} \cdots t_{\delta_4} t_{\tilde{\delta}_5} t_{\tilde\delta_6}$ in $\Mod(\Sigma_g^4)$ and apply the lantern relations with $\{ a_1, x_1, b_1, y_1 \}$ and $\{ a_3, x_2, b_3, y_2 \}$:
\begin{align}
	t_{\delta_1} \cdots t_{\delta_4} 
	&= t_{B_0} \cdots t_{B_{g-1}} t_{a_1} t_{a_3} t_{b_1} t_{b_3} t_{\tilde{\delta}_5}^{-1} \cdot t_{B_0} \cdots t_{B_{g-1}} t_{x_1} t_{x_2} t_{y_1} t_{y_2} t_{\tilde{\delta}_6}^{-1}  \notag \\
	&= t_{B_0} \cdots t_{B_{g-1}} \underline{t_{a_1} t_{x_1} t_{b_1} t_{y_1} t_{\tilde{\delta}_5}^{-1}} \cdot {}_{t_{a_3}t_{b_3} t_{x_1}^{-1}t_{y_1}^{-1}}(t_{B_0} \cdots t_{B_{g-1}}) \; \underline{t_{a_3} t_{x_2} t_{b_3} t_{y_2} t_{\tilde{\delta}_6}^{-1}} \notag \\
	&= t_{B_{0,1}} t_{B_{1,1}} \cdots t_{B_{g,1}} t_{C_1} \cdot 
	t_{B_{0,2}} t_{B_{1,2}} \cdots t_{B_{g,2}} t_{C_2}, \label{eq:WIBeven}
\end{align}	
where the resulting curves are shown in Figure~\ref{F:WIBeven}.
Let $W_{\mathrm{\rmI}B}$ denote the last word.
\qed

\subsubsection{Type $\mathrm{\rmII}A$}
Consider the relation $W_{\mathrm{\rmII}} = t_{\delta_1} \cdots t_{\delta_4} t_{\tilde{\delta}_5} t_{\tilde\delta_6}$ in $\Mod(\Sigma_g^4)$ and apply the lantern relations with $\{ a_1, x_1, b_3, y_2 \}$ and $\{ a_3, x_2, b_2, y_1 \}$:
\begin{align}
	t_{\delta_1} \cdots t_{\delta_4} 
	&= t_{B_0^{\prime}} \cdots t_{B_{g-1}^{\prime}} t_{a_1} t_{a_3} t_{b_2} t_{b_3} t_{\tilde{\delta}_5}^{-1} \cdot t_{B_0^{\prime}} \cdots t_{B_{g-1}^{\prime}} t_{x_1} t_{x_2} t_{y_1} t_{y_2} t_{\tilde{\delta}_6}^{-1}  \notag \\
	&= {}_{t_{b_2}^{-1} t_{y_2}}(t_{B_0^{\prime}} \cdots t_{B_{g-1}^{\prime}}) \; \underline{t_{a_1} t_{x_1} t_{b_3} t_{y_2} t_{\tilde{\delta}_5}^{-1}} \cdot {}_{t_{a_3} t_{x_1}^{-1}}(t_{B_0^{\prime}} \cdots t_{B_{g-1}^{\prime}}) \; \underline{t_{a_3} t_{x_2} t_{b_2} t_{y_1} t_{\tilde{\delta}_6}^{-1}} \notag \\
	&= t_{B_{0,1}} t_{B_{1,1}} \cdots t_{B_{g,1}} t_{C_1} \cdot 
	t_{B_{0,2}} t_{B_{1,2}} \cdots t_{B_{g,2}} t_{C_2}, \label{eq:WIIAeven}
\end{align}	
where the resulting curves are shown in Figure~\ref{F:WIIAeven}.
Let $W_{\mathrm{\rmII}A}$ denote the last word.
\qed

\subsubsection{Type $\mathrm{\rmII}B$}
Consider the relation $W_{\mathrm{\rmII}} = t_{\delta_1} \cdots t_{\delta_4} t_{\tilde{\delta}_5} t_{\tilde\delta_6}$ in $\Mod(\Sigma_g^4)$ and apply the lantern relations with $\{ a_1, x_1, b_2, y_1 \}$ and $\{ a_3, x_2, b_3, y_2 \}$:
\begin{align}
	t_{\delta_1} \cdots t_{\delta_4} 
	&= t_{B_0^{\prime}} \cdots t_{B_{g-1}^{\prime}} t_{a_1} t_{a_3} t_{b_2} t_{b_3} t_{\tilde{\delta}_5}^{-1} \cdot t_{B_0^{\prime}} \cdots t_{B_{g-1}^{\prime}} t_{x_1} t_{x_2} t_{y_1} t_{y_2} t_{\tilde{\delta}_6}^{-1}  \notag \\
	&= t_{B_0^{\prime}} \cdots t_{B_{g-1}^{\prime}} \underline{t_{a_1} t_{x_1} t_{b_2} t_{y_1} t_{\tilde{\delta}_5}^{-1}} \cdot {}_{t_{a_3}t_{b_3} t_{x_1}^{-1}t_{y_1}^{-1}}(t_{B_0^{\prime}} \cdots t_{B_{g-1}^{\prime}}) \; \underline{t_{a_3} t_{x_2} t_{b_3} t_{y_2} t_{\tilde{\delta}_6}^{-1}} \notag \\
	&= t_{B_{0,1}} t_{B_{1,1}} \cdots t_{B_{g,1}} t_{C_1} \cdot 
	t_{B_{0,2}} t_{B_{1,2}} \cdots t_{B_{g,2}} t_{C_2}, \label{eq:WIIBeven}
\end{align}	
where the resulting curves are shown in Figure~\ref{F:WIIBeven}.
Let $W_{\mathrm{\rmII}B}$ denote the last word.
\qed

\subsubsection*{Summary}
Each of the four relations, $W_{\mathrm{\rmI}A}=t_{\delta_1} \cdots t_{\delta_4}$, 
$W_{\mathrm{\rmI}B}=t_{\delta_1} \cdots t_{\delta_4}$, $W_{\mathrm{\rmII}A}=t_{\delta_1} \cdots t_{\delta_4}$ and $W_{\mathrm{\rmII}B}=t_{\delta_1} \cdots t_{\delta_4}$, is a lift of the MCK relation $W=1$ for even $g$.
Besides, two of them, $W_{\mathrm{\rmI}A}=t_{\delta_1} \cdots t_{\delta_4}$ and 
$W_{\mathrm{\rmI}B}=t_{\delta_1} \cdots t_{\delta_4}$, are further lifts of the known lift of the MCK relation $W=t_{\delta_1} t_{\delta_2}$.

\begin{rmk}
	There is one more known lift of Matsumoto's relation $W=1$ of genus $2$ to $\Sigma_2^2$ which was found by Baykur-Hayano~\cite[Lemma 4.6]{BaykurHayano2016}.
	The relation $(t_{d_4} t_{d_3} t_{d_2})^2 t_{d_+} t_{d_-} = t_{\delta_1} t_{\delta_2}$ in question can be seen as a lift of $W=1$ as follows:
	\begin{align*}
		t_{\delta_1} t_{\delta_2} 
		&= (\underline{t_{d_4} t_{d_3} t_{d_2}})^2 t_{d_+} t_{d_-}  \\
		&\sim (\underline{t_{B_2} t_{t_{d_4}(d_2)}} t_{C})^2 t_{B_0} t_{d_-} \\
		&\sim t_{B_0} (t_{B_1} t_{B_2} t_{C}) \underline{(t_{B_1} t_{B_2} t_{C}) t_{d_-}} \\
		&\sim (t_{B_0} t_{B_1} t_{B_2} t_{C})^2,
	\end{align*}
	where we put $B_0=d_+$, $C=d_2$, $B_2 = t_{d_4}(d_3)$ and $B_1 = t_{B_2} t_{d_4} (d_2)$ and used an observation that $t_{B_1} t_{B_2} t_{C}(d_-) = B_0$.
	By the symbol ``$\sim$'' we mean the Hurwitz equivalence.
	The last factorization actually coincides with the factorization $W_{\mathrm{\rmII}A} = W_{\mathrm{\rmII}B} = t_{\delta_1} t_{\delta_2}$ where $\delta_3$ and $\delta_4$ have been capped off.
	It follows simultaneously that $W_{\mathrm{\rmII}A}=t_{\delta_1} \cdots t_{\delta_4}$ and $W_{\mathrm{\rmII}B}=t_{\delta_1} \cdots t_{\delta_4}$ are further lifts of Baykur-Hayano's lift.
\end{rmk}

\section{Sections of the MCK Lefschetz fibration and the corresponding Lefschetz pencils} \label{S:sections-LPs}

\subsection{$(-1)$-sections of the MCK Lefschetz fibration} \
Each of the lifts of the MCK relation constructed in Sections~\ref{subsect:odgenus} and~\ref{subsect:evengenus} implies the existence of four or eight disjoint $(-1)$-sections of the MCK Lefschetz fibration.
By employing the classification of symplectic $4$-manifolds of symplectic Kodaira dimension $-\infty$, we can indeed deduce that those sections are maximal, that is, there are no more $(-1)$-sections disjoint from the above ones.

It is a fundamental theorem in the theory of Lefschetz pencils that an oriented closed $4$-manifold admits a symplectic structure if and only if it admits a Lefschetz pencil~\cite{Donaldson1999, GompfStipsicz1999}.
In the case of the MCK Lefschetz fibration $f_{W}: X_W \rightarrow S^2$ the total space $X_W$ is (blow up of) a ruled surface (a $S^2$-bundle over a compact Riemann surface), which is obviously symplectic.
This class of symplectic $4$-manifolds has been well-understood in terms of the \textit{symplectic Kodaira dimension} $\kappa_{sym}$, an analogue of the usual Kodaira dimension for algebraic surfaces, which is a numerical invariant taking the values in $\{ -\infty, 0, 1, 2\}$ (for this notion, see~\cite{OzbagciStipsicz2004book}).
That is to say, a minimal symplectic $4$-manifold with $\kappa_{sym} = -\infty$ is diffeomorphic to a rational or ruled surface and vice versa. 
Since symplectic Kodaira dimension depends only on oriented diffeomorphism types and preserves under blow up operation, $\Sigma \times S^2 \# k\overline{\mathbb{CP}}{}^{2}$ also has $\kappa_{sym} = -\infty$, where $\Sigma$ is an oriented closed surface.
From the above classification, it is easy to see that if $\Sigma \times S^2 \# k\overline{\mathbb{CP}}{}^{2} \cong X \# k\overline{\mathbb{CP}}{}^{2}$ then $X$ has to be a minimal ruled surface unless $\Sigma = S^2$ \footnote{Note that $S^2 \times S^2 \# \overline{\mathbb{CP}}{}^{2} \cong S^2 \tilde{\times} S^2 \# \overline{\mathbb{CP}}{}^{2} \cong \mathbb{CP}{}^{2} \# 2\overline{\mathbb{CP}}{}^{2}$.}.
Therefore $\Sigma \times S^2 \# k\overline{\mathbb{CP}}{}^{2}$ cannot have more than $k$ disjoint $(-1)$-spheres.
Applying this argument to $X_W = \Sigma_{g/2} \times S^2 \# 4\overline{\mathbb{CP}}{}^{2}$ for even $g \geq 2$ or $X_W = \Sigma_{(g-1)/2} \times S^2 \# 8\overline{\mathbb{CP}}{}^{2}$ for odd $g \geq 3$, we conclude the following:

\begin{theorem}
	The genus-$g$ Matsumoto-Cadavid-Korkmaz Lefschetz fibration  $f_W : X_W \rightarrow S^2$ admits four disjoint $(-1)$-sections when $g$ is even and eight disjoint $(-1)$-sections when $g$ is odd.
	When $g \geq 2$, this number of disjoint $(-1)$-sections is the maximal for $f_W$.
\end{theorem}
\begin{remark} \label{rmk:KO}
	When $g=1$, the MCK Lefschetz fibration is nothing but the elliptic Lefschetz fibration $E(1) \rightarrow S^2$ and this fibration comes from an algebraic pencil of curves that has nine base points.
	Hence, it naturally has \textit{nine} disjoint $(-1)$-sections.
	A set of nine disjoint $(-1)$-sections was also explicitly given by Korkmaz-Ozbagci~\cite{KorkmazOzbagci2008} in the form of a monodromy factorization of the boundary multi-twist.
	We can also achieve such nine disjoint $(-1)$-sections from $W_{\mathrm{\rmI}A} = t_{\delta_1} \cdots t_{\delta_8}$ or $W_{\mathrm{\rmII}A}  = t_{\delta_1} \cdots t_{\delta_8}$ of genus $1$ by further performing a lantern breeding, say, with respect to $\{ a_3, b_4, \delta_3 \}$ or $\{ a_4, b_4, \delta_1 \}$, respectively \footnote{In fact, it is possible (though not obvious) to show that $W_{\mathrm{\rmI}A} = t_{\delta_1} \cdots t_{\delta_8}$ and $W_{\mathrm{\rmII}A}  = t_{\delta_1} \cdots t_{\delta_8}$ are Hurwitz equivalent to Korkmaz-Ozbagci's $8$-holed torus relation, to which they performed a lantern breeding to obtain the ninth $(-1)$-section.}. 
	Compare this fact with Table~\ref{T:diffeotypeodd}; each of types $\mathrm{\rmI}A$ and $\mathrm{\rmII}A$ has an extra $(-1)$-sphere ($S^2 \tilde \times S^2 \cong \mathbb{CP}{}^{2} \# \overline{\mathbb{CP}}{}^{2}$) while the others do not.
\end{remark}
\begin{remark}
	From the monodromy factorizations of the boundary multi-twist, we can explicitly locate the sections in a handlebody diagram of the MCK Lefschetz fibration by following the method explained in Section~$4$ of~\cite{KorkmazOzbagci2008}, though we do not dare to draw them here.
\end{remark}

\subsection{The supporting Lefschetz pencils} \
By blowing down the maximal disjoint $(-1)$-sections discussed above, we obtain several minimal Lefschetz pencils that support the MCK Lefschetz fibration.
We study some topological aspects of those pencils here. 
For a type $T \in \{ \mathrm{\rmI}A, \mathrm{\rmI}B, \mathrm{\rmII}A, \mathrm{\rmII}B \}$ for even $g$ we denote by $f_{W_T} : X_{W_T} \setminus B_{W_T} \rightarrow S^2$ the Lefschetz pencil corresponding to $W_T = t_{\delta_1} \cdots t_{\delta_4}$ and call it \textit{the type $T$ MCK Lefschetz pencil}.
The same notation for odd $g$ is used as well.

\subsubsection{The diffeomorphism types} \
We first determine the diffeomorphism types of the total spaces of the MCK Lefschetz pencils.
Since the manifolds are ruled surfaces, which are topologically the product (the trivial $S^2$-bundle) $\Sigma \times S^2$ or the nontrivial $S^2$-bundle $\Sigma \tilde \times S^2$ where $\Sigma$ is the base surface, it is sufficient to determine whether they are spin or not.

\paragraph{\it Nonspin} \
Some of the pencils are indeed easy to judge \textit{nonspin}.
Let us consider a Lefschetz pencil with a fixed monodromy $t_{a_n} \cdots t_{a_1} = t_{\delta_1} \cdots t_{\delta_k}$.
Suppose that a subsurface $S$ of the reference fiber is bounded by $p$ vanishing cycles $a_{i_1}, \cdots, a_{i_p}$ and contains $q$ base points corresponding to $\delta_{j_1}, \cdots, \delta_{j_q}$. 
Then by clustering the Lefschetz critical points corresponding to $a_{i_1}, \cdots, a_{i_p}$ on the same singular fiber, we can obtain a fiber component $S^{\prime}$ that contains $p$ critical points and $q$ base points.
The self-intersection of this closed surface $S^{\prime}$ turns out to be $q-p$.

Now consider the type $\mathrm{\rmI}B$ MCK Lefschetz pencil for even $g$, which corresponds to $W_{\mathrm{\rmI}B}=t_{\delta_1} \cdots t_{\delta_4}$.
See Figure~\ref{F:WIBeven}.
The vanishing cycle $C_1$ bounds a subsurface of genus $g/2$ that contains no base points.
We can thus obtain an embedded surface of self-intersection $-1$, which is odd.
It follows that the total space $X_{W_{\mathrm{\rmI}B}}$ cannot be spin.
Similarly, we can find embedded surfaces of odd self-intersection numbers for some of the other pencils.
For $W_{\mathrm{\rmII}B}=t_{\delta_1} \cdots t_{\delta_4}$ for even $g$, $\{ C_1, \delta_2, \delta_4 \}$ yields a surface of self-intersection $1$.
For $W_{\mathrm{\rmI}A}=t_{\delta_1} \cdots t_{\delta_8}$ for odd $g$, $\{ a_3, b_4, \delta_3 \}$ yields a surface of self-intersection $-1$.
For $W_{\mathrm{\rmII}A}=t_{\delta_1} \cdots t_{\delta_8}$ for odd $g$, $\{ a_4, b_4, \delta_1 \}$ yields a surface of self-intersection $-1$.
Consequently, those pencils are nonspin.

\paragraph{\it Spin} \
In contrast, to judge a Lefschetz pencil to be spin would need more work to do.
We employ the following criteria, which was established in~\cite{BaykurHayanoMonden_preprint} generalizing a theorem for Lefschetz fibrations by Stipsicz~\cite{Stipsicz2001}: 
\begin{theorem}[\cite{BaykurHayanoMonden_preprint}] \label{Thm:BHM}
	Let $f : X \setminus B \rightarrow S^2$ be a Lefschetz pencil and $t_{a_n} \cdots t_{a_2} t_{a_1} = t_{\delta_1} t_{\delta_2} \cdots t_{\delta_k}$ the corresponding monodromy factorization. 
	The manifold $X$ admits a spin structure if and only if there exists a quadratic form $q : H_1(\Sigma_g^k; \Z / 2\Z) \rightarrow Z / 2Z $ with respect to the intersection pairing $\langle \; , \; \rangle$ of $H_1(\Sigma_g^k; \Z / 2\Z)$ \footnote{Namely, $q$ is a map such that $q(0)=0$ and $\langle x,y \rangle =q(x+y) + q(x) + q(y)$.} such that 
	\begin{enumerate}
		\item[(A)] $q(a_i) = 1$ for any $i$, and 
		\item[(B)] $q(\delta_j) = 1$ for some $j$.
	\end{enumerate}
\end{theorem}
\begin{figure}[htbp]
	\centering
	\subfigure[For Type $\mathrm{\rmI}A$ of even $g$. \label{F:Spin_even_IA}]
	{\includegraphics[height=75pt]{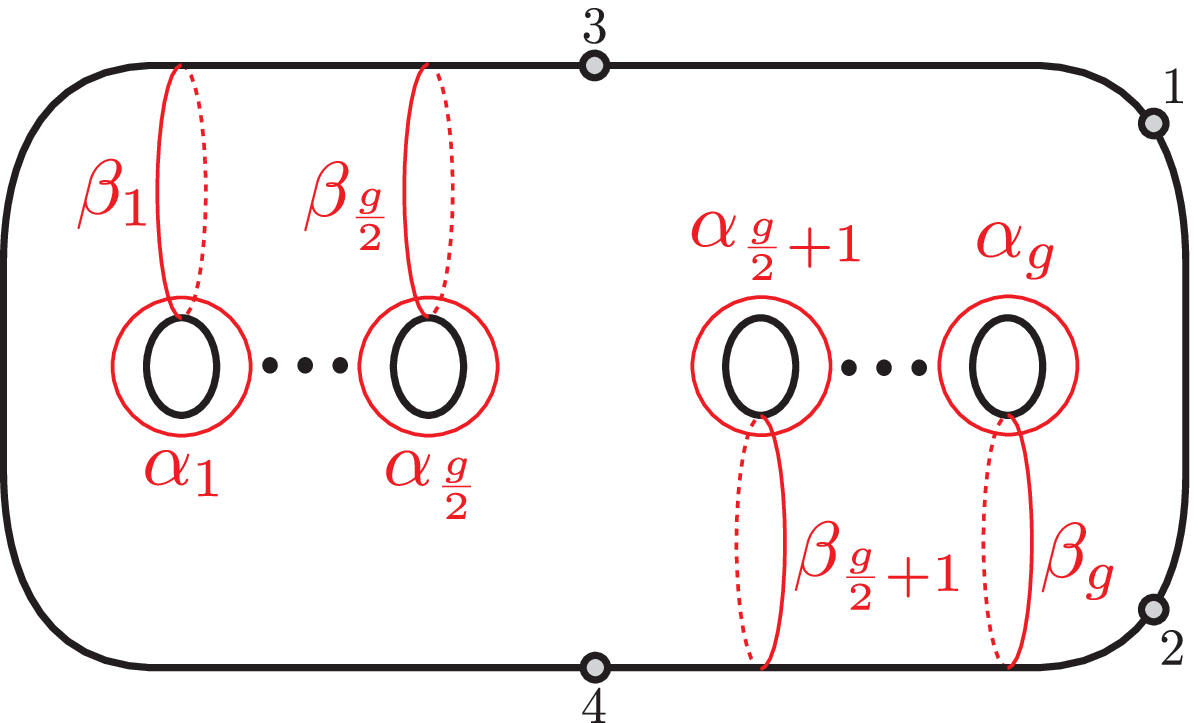}} 
	\hspace{-5pt}	
    \setcounter{subfigure}{2}
  	\subfigure[For Type $\mathrm{\rmI}B\sharp$ of odd $g$. \label{F:Spin_odd_IBs}]
	{\includegraphics[height=75pt]{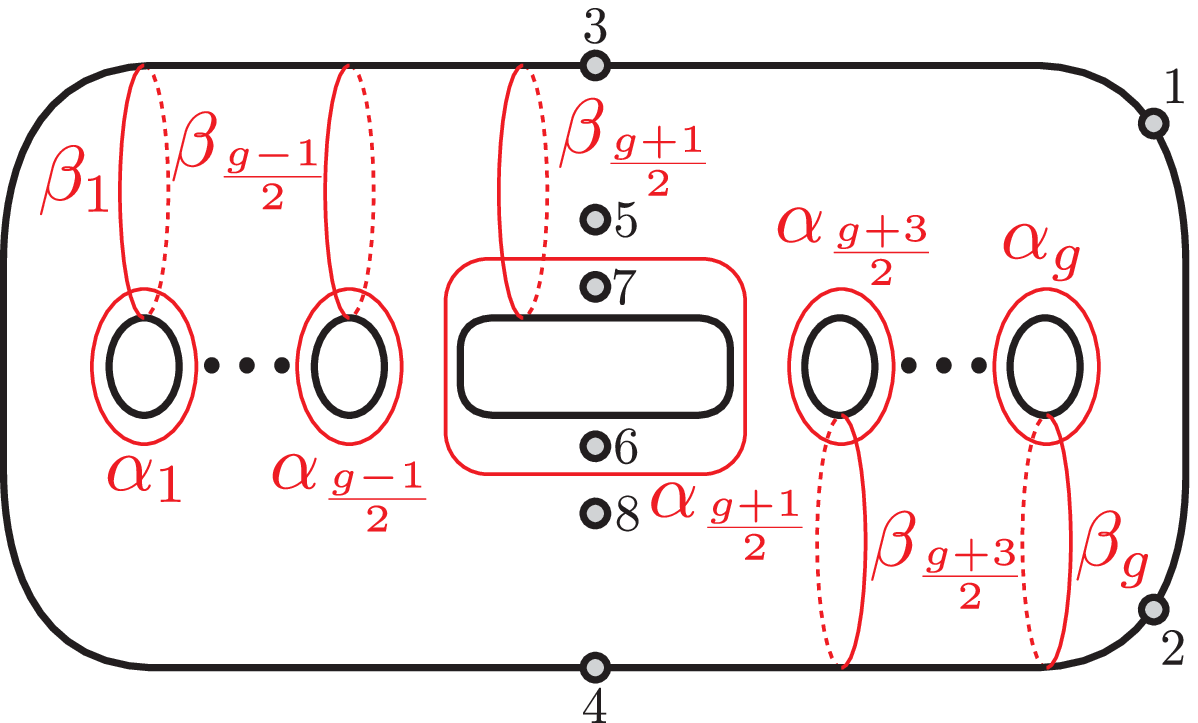}} 
	\hspace{-5pt}	
    \setcounter{subfigure}{4}
    \subfigure[For Type $\mathrm{\rmII}B\sharp$ of odd $g$. \label{F:Spin_odd_IIBs}]
    {\includegraphics[height=75pt]{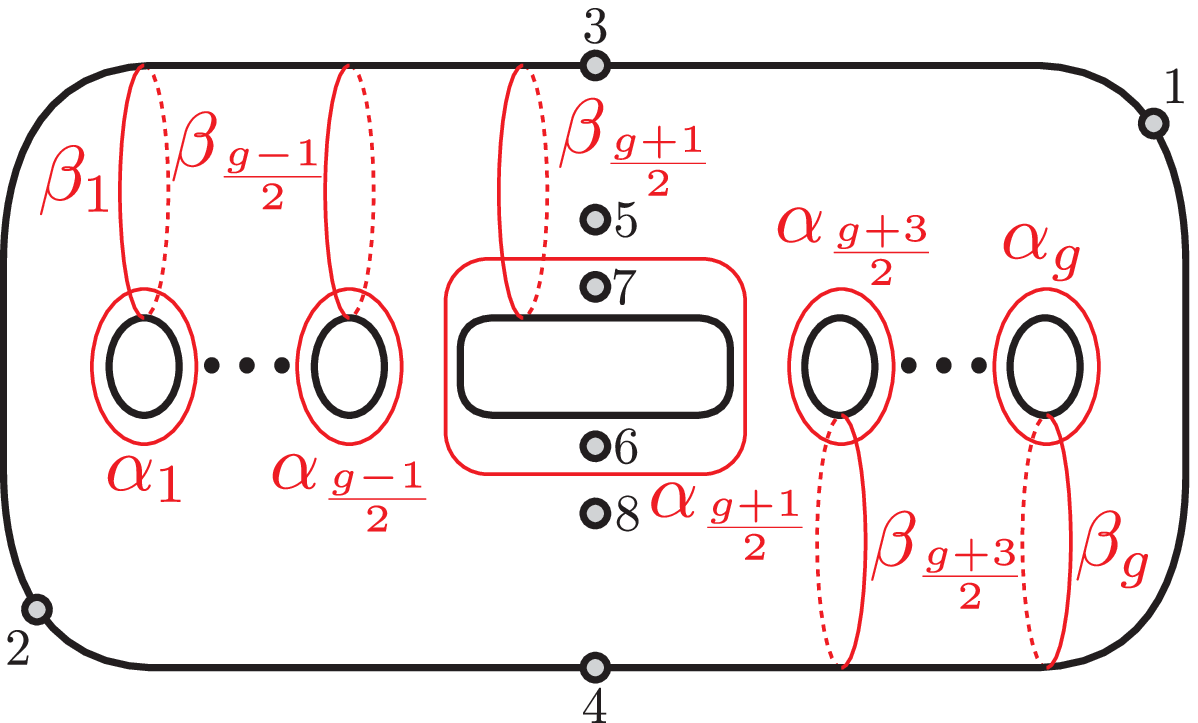}} 
    \hspace{-5pt}	
    \setcounter{subfigure}{1}
    \subfigure[For Type $\mathrm{\rmII}A$ of even $g$. \label{F:Spin_even_IIA}]
    {\includegraphics[height=75pt]{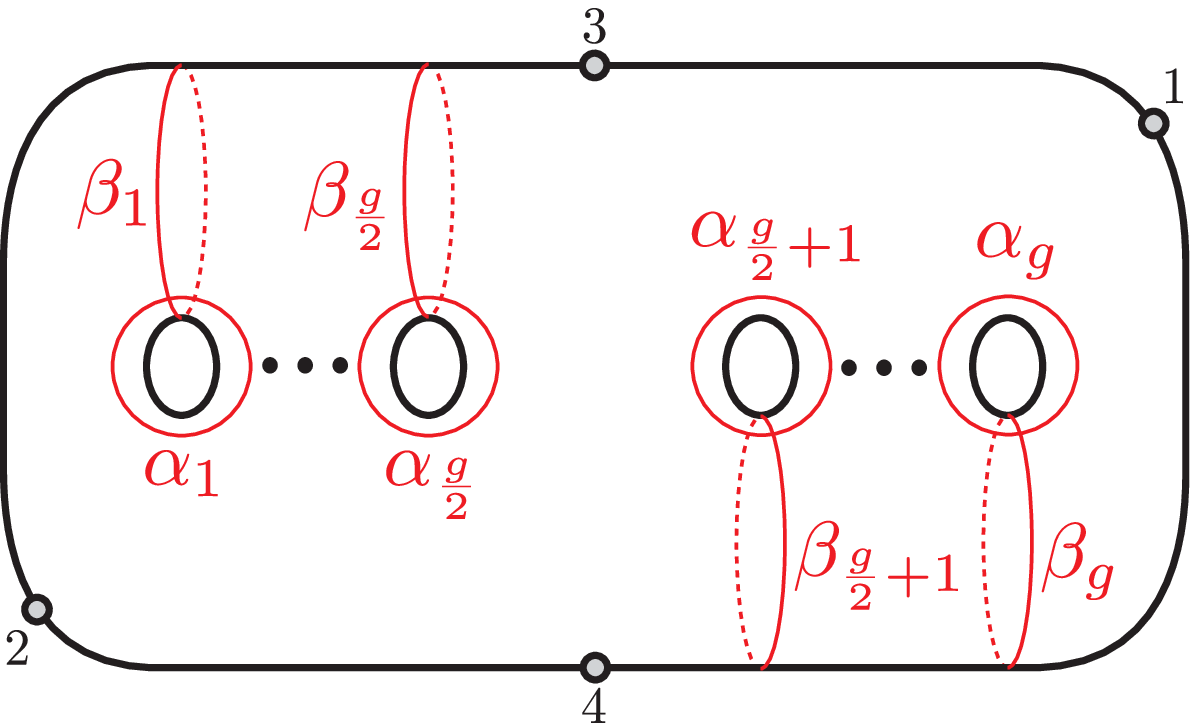}} 
    \hspace{-5pt}	
	\setcounter{subfigure}{3}
	\subfigure[For Type $\mathrm{\rmI}B\flat$ of odd $g$. \label{F:Spin_odd_IBf}]
	{\includegraphics[height=75pt]{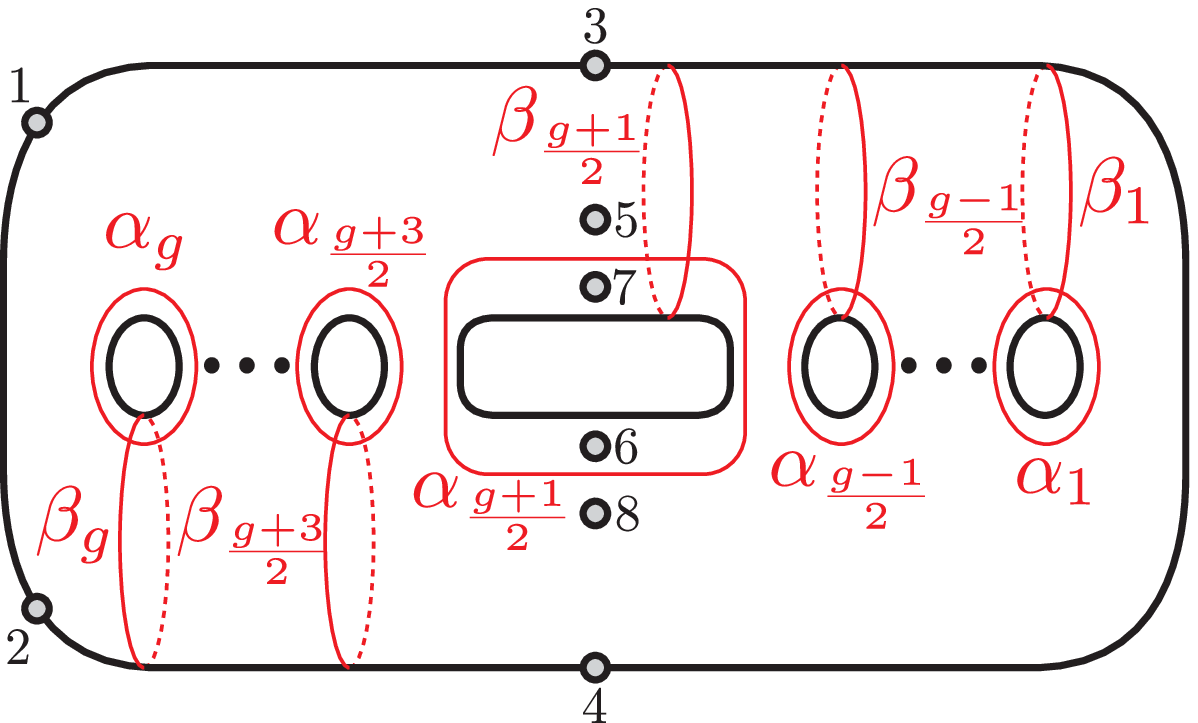}} 
	\hspace{-5pt}	
    \setcounter{subfigure}{5}
	\subfigure[For Type $\mathrm{\rmII}B\flat$ of odd $g$. \label{F:Spin_odd_IIBf}]
	{\includegraphics[height=75pt]{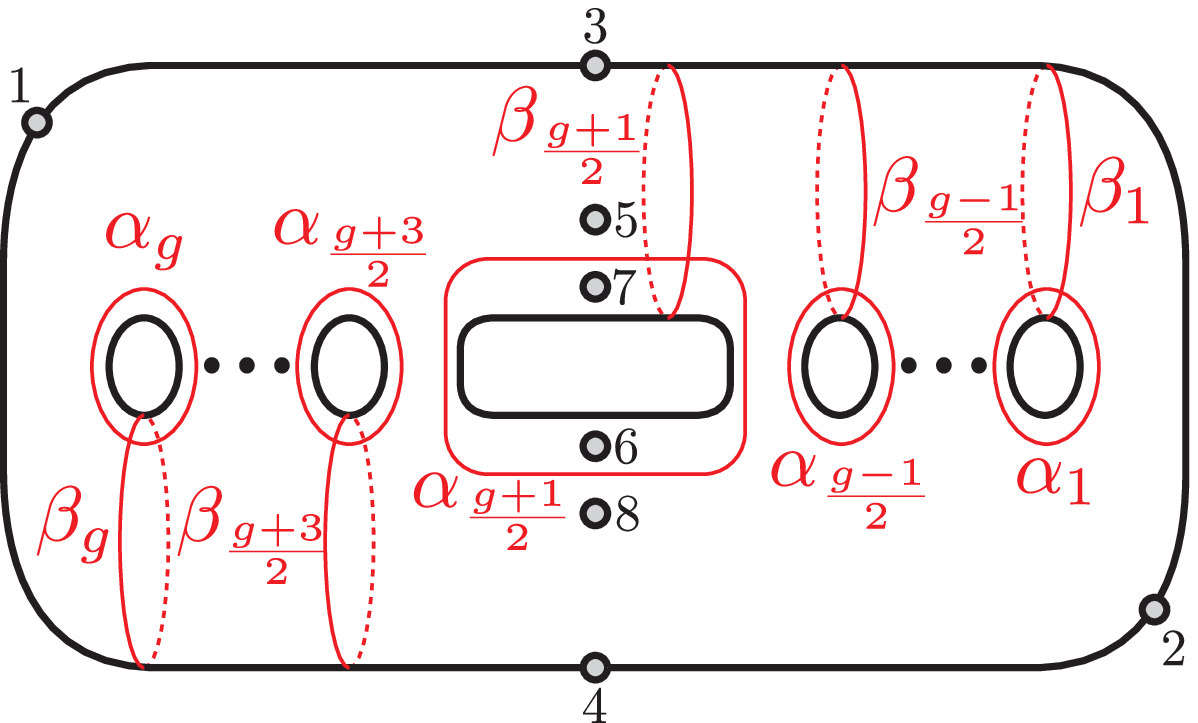}} 
	\caption{Generators of $H_1(\Sigma_g^4; \Z / 2\Z)$ or $H_1(\Sigma_g^8; \Z / 2\Z)$.} 	
\end{figure}

\paragraph{\it Type $\mathrm{\rmI}A$ of even $g$} \
To begin with, we consider $W_{\mathrm{\rmI}A}=t_{\delta_1} \cdots t_{\delta_4}$ for even $g$.
Take the generators $\alpha_1, \cdots, \alpha_g, \beta_1, \cdots, \beta_g, \delta_1, \cdots, \delta_4$ (which include an extra generator) of $H_1(\Sigma_g^4; \Z / 2\Z)$ as in Figure~\ref{F:Spin_even_IA}.
The $\Z/2\Z$ homology classes of the vanishing cycles are calculated as follows:
\begin{align*}
	\quad &B_{0,1} = \alpha_1 + \cdots + \alpha_g + \delta_2; \\
	&B_{1,1} = B_{0,1} + \beta_1 + \beta_g + \delta_2 = \alpha_1 + \cdots + \alpha_g + \beta_1 + \beta_g; \\
	&B_{2i,1} = B_{2i-1,1} + \alpha_i + \alpha_{g+1-i} \quad \text{for $2i=2, 4, \cdots, g$};\\
	&B_{2i+1,1} =  B_{2i,1} + \beta_i + \beta_{i+1} + \beta_{g-i} + \beta_{g+1-i} \quad \text{for $2i+1=3, 5, \cdots, g-1$}; \\
	\intertext{hence,} 
	&B_{2i,1} = \alpha_{i+1} + \cdots + \alpha_{g-i} + \beta_i + \beta_{g+1-i} \quad \text{for $2i=2, 4, \cdots, g-2$}; \\
	&B_{g,1} = \beta_{g/2} + \beta_{g/2+1}; \\
	&B_{2i+1,1} = \alpha_{i+1} + \cdots + \alpha_{g-i} + \beta_{i+1} + \beta_{g-i} \quad \text{for $2i+1=3, 5, \cdots, g-1$}; \\
	\intertext{and}
	&C_1 = \delta_4; \\
	&B_{i,2} = B_{i,1} + \delta_3 + \delta_4 \quad \text{for $i= 1, \cdots, g$}; \\
	&C_2 = \delta_3.
\end{align*}
Put $q(\alpha_i)=0$ for all $i$, $q(\beta_i)=1$ for $i=1, \cdots, g/2$, $q(\beta_i)=0$ for $i=g/2+1, \cdots, g$ and $q(\delta_j)=1$ for all $j$ and extend $q$ to a quadratic form on $H_1(\Sigma_g^4; \Z / 2\Z)$.
(This extension is possible since $q$ is consistent with the only defining relation $\delta_1 + \cdots + \delta_4 = 0$.)
Then it is easy to check that $q(B_{i,j})=q(C_j)=1$ for $i=1, \cdots, g$ and $j=1,2$.
Thus, $q$ satisfies $(A)$ and $(B)$ in Theorem~\ref{Thm:BHM}, which implies that the total space $X_{\mathrm{\rmI}A}$ is spin.
\qed

\paragraph{\it Type $\mathrm{\rmII}A$ of even $g$} \
Take the generators of $H_1(\Sigma_g^4; \Z / 2\Z)$ as in Figure~\ref{F:Spin_even_IIA}.
The $\Z/2\Z$ homology classes of the vanishing cycles are
\begin{align*}
	\quad &B_{0,1} = \alpha_1 + \cdots + \alpha_g + \delta_4; \\
	&B_{1,1} = B_{0,1} + \beta_1 + \beta_g = \alpha_1 + \cdots + \alpha_g + \beta_1 + \beta_g + \delta_4; \\
	&B_{2i,1} = B_{2i-1,1} + \alpha_i + \alpha_{g+1-i} \\
	&\qquad = \alpha_{i+1} + \cdots + \alpha_{g-i} + \beta_i + \beta_{g+1-i} + \delta_4 \quad \text{for $2i=2, 4, \cdots, g-1$}; \\
	&B_{g,1} = \beta_{g/2} + \beta_{g/2+1} + \delta_4; \\
	&B_{2i+1,1} =  B_{2i,1} + \beta_i + \beta_{i+1} + \beta_{g-i} + \beta_{g+1-i} \\
	&\qquad = \alpha_{i+1} + \cdots + \alpha_{g-i} + \beta_{i+1} + \beta_{g-i} + \delta_4 \quad \text{for $2i+1=3, 5, \cdots, g-2$}; \\
	&C_1 = \delta_2; \\
	&B_{i,2} = B_{i,1} + \delta_3 + \delta_4 \quad \text{for $i= 1, \cdots, g$}; \\
	&C_2 = \delta_1.
\end{align*}
Put $q(a_i)=q(\beta_i)=0$ for all $i$ and $q(\delta_j)=1$ for all $j$, then $q$ satisfies $(A)$ and $(B)$.
\qed

\paragraph{\it Type $\mathrm{\rmI}B\sharp$ of odd $g$} \
Take the generators of $H_1(\Sigma_g^8; \Z / 2\Z)$ as in Figure~\ref{F:Spin_odd_IBs}.
The $\Z/2\Z$ homology classes of the vanishing cycles are
\begin{align*}
	\quad &B_{0,1} = \alpha_1 + \cdots + \alpha_g + \delta_2 + \delta_4 + \delta_5 + \delta_8; \\
	&B_{1,1} = B_{0,1} + \beta_1 + \beta_g + \delta_2 = \alpha_1 + \cdots + \alpha_g + \beta_1 + \beta_g + \delta_4 + \delta_5 + \delta_8; \\
	&B_{2i,1} = B_{2i-1,1} + \alpha_i + \alpha_{g+1-i} \\
	&\qquad = \alpha_{i+1} + \cdots + \alpha_{g-i} + \beta_i + \beta_{g+1-i} + \delta_4 + \delta_5 + \delta_8 \quad \text{for $2i=2, 4, \cdots, g-1$}; \\
	&B_{2i+1,1} =  B_{2i,1} + \beta_i + \beta_{i+1} + \beta_{g-i} + \beta_{g+1-i} \\
	&\qquad = \alpha_{i+1} + \cdots + \alpha_{g-i} + \beta_{i+1} + \beta_{g-i} + \delta_4 + \delta_5 + \delta_8 \quad \text{for $2i+1=3, 5, \cdots, g-2$}; \\
	&B_{g,1} = \alpha_{(g+1)/2} + \beta_{(g+1)/2} + \beta_{(g+1)/2} + \delta_5 + \delta_6 = \alpha_{(g+1)/2} + \delta_5 + \delta_6; \\
	&a_1 = \beta_{(g+1)/2} + \delta_5; \quad a_2 = \beta_{(g+1)/2} + \delta_7; \\
	&b_1 = \beta_{(g+1)/2} + \delta_6; \quad b_2 = \beta_{(g+1)/2} + \delta_8; \\
	&B_{i,2} = B_{i,1} + \delta_3 + \delta_4 \quad \text{for $i= 1, \cdots, g$}; \\
	&a_3 = \beta_{(g+1)/2} + \delta_3; \quad a_4 = \beta_{(g+1)/2} + \delta_3 + \delta_5 + \delta_7; \\
	&b_3 = \beta_{(g+1)/2} + \delta_4; \quad b_4 = \beta_{(g+1)/2} + \delta_4 + \delta_6 + \delta_8.
\end{align*}
Put $q(\alpha_i)=1$ for all $i$, $q(\beta_i)=0$ for $i=1, \cdots, (g+1)/2$, $q(\beta_i)=1$ for $i=(g+3)/2, \cdots, g$ and $q(\delta_j)=1$ for all $j$, then $q$ satisfies $(A)$ and $(B)$.
\qed

\paragraph{\it Type $\mathrm{\rmI}B\flat$ of odd $g$} \
Take the generators of $H_1(\Sigma_g^8; \Z / 2\Z)$ as in Figure~\ref{F:Spin_odd_IBf}.
The $\Z/2\Z$ homology classes of the vanishing cycles are
\begin{align*}
	\quad &B_{0,2} = \alpha_1 + \cdots + \alpha_g + \delta_2 + \delta_4 + \delta_7 + \delta_6; \\
	&B_{1,2} = B_{0,2} + \beta_1 + \beta_g + \delta_2 = \alpha_1 + \cdots + \alpha_g + \beta_1 + \beta_g + \delta_4 + \delta_7 + \delta_6; \\
	&B_{2i,2} = B_{2i-1,2} + \alpha_i + \alpha_{g+1-i} \\
	&\qquad = \alpha_{i+1} + \cdots + \alpha_{g-i} + \beta_i + \beta_{g+1-i} + \delta_4 + \delta_7 + \delta_6 \quad \text{for $2i=2, 4, \cdots, g-1$}; \\
	&B_{2i+1,2} =  B_{2i,2} + \beta_i + \beta_{i+1} + \beta_{g-i} + \beta_{g+1-i} \\
	&\qquad = \alpha_{i+1} + \cdots + \alpha_{g-i} + \beta_{i+1} + \beta_{g-i} + \delta_4 + \delta_7 + \delta_6 \quad \text{for $2i+1=3, 5, \cdots, g-2$}; \\
	&B_{g,2} = \alpha_{(g+1)/2} + \beta_{(g+1)/2} + \beta_{(g+1)/2} + \delta_7 + \delta_8 = \alpha_{(g+1)/2} + \delta_7 + \delta_8; \\
	&a_1 = \beta_{(g+1)/2} + \delta_3 + \delta_7; \quad a_2 = \beta_{(g+1)/2} + \delta_3 + \delta_5; \\
	&b_1 = \beta_{(g+1)/2} + \delta_4 + \delta_8; \quad b_2 = \beta_{(g+1)/2} + \delta_4 + \delta_6; \\
	&B_{i,1} = B_{i,2} + \delta_3 + \delta_4 \quad \text{for $i= 1, \cdots, g$}; \\
	&a_3 = \beta_{(g+1)/2} + \delta_5 + \delta_7; \quad a_4 = \beta_{(g+1)/2}; \\
	&b_3 = \beta_{(g+1)/2} + \delta_6 + \delta_8; \quad b_4 = \beta_{(g+1)/2}.
\end{align*}
Put $q(\alpha_i)=1$ for all $i$, $q(\beta_i)=1$ for $i=1, \cdots, (g+1)/2$, $q(\beta_i)=0$ for $i=(g+3)/2, \cdots, g$ and $q(\delta_j)=1$ for all $j$, then $q$ satisfies $(A)$ and $(B)$.
\qed

\paragraph{\it Type $\mathrm{\rmII}B\sharp$ of odd $g$} \
Take the generators of $H_1(\Sigma_g^8; \Z / 2\Z)$ as in Figure~\ref{F:Spin_odd_IIBs}.
The $\Z/2\Z$ homology classes of the vanishing cycles are
\begin{align*}
	\quad &B_{0,1} = \alpha_1 + \cdots + \alpha_g + \delta_5 + \delta_8; \\
	&B_{1,1} = B_{0,1} + \beta_1 + \beta_g = \alpha_1 + \cdots + \alpha_g + \beta_1 + \beta_g + \delta_5 + \delta_8; \\
	&B_{2i,1} = B_{2i-1,1} + \alpha_i + \alpha_{g+1-i} \\
	&\qquad = \alpha_{i+1} + \cdots + \alpha_{g-i} + \beta_i + \beta_{g+1-i} + \delta_5 + \delta_8 \quad \text{for $2i=2, 4, \cdots, g-1$}; \\
	&B_{2i+1,1} =  B_{2i,1} + \beta_i + \beta_{i+1} + \beta_{g-i} + \beta_{g+1-i} \\
	&\qquad = \alpha_{i+1} + \cdots + \alpha_{g-i} + \beta_{i+1} + \beta_{g-i} + \delta_5 + \delta_8 \quad \text{for $2i+1=3, 5, \cdots, g-2$}; \\
	&B_{g,1} = \alpha_{(g+1)/2} + \beta_{(g+1)/2} + \beta_{(g+1)/2} + \delta_2 + \delta_4 + \delta_5 + \delta_6 = \alpha_{(g+1)/2} + \delta_2 + \delta_4 + \delta_5 + \delta_6; \\
	&a_1 = \beta_{(g+1)/2} + \delta_5; \quad a_2 = \beta_{(g+1)/2} + \delta_7; \\
	&b_1 = \beta_{(g+1)/2} + \delta_2 + \delta_4 + \delta_6; \quad b_2 = \beta_{(g+1)/2} + \delta_2 + \delta_4 + \delta_8; \\
	&B_{i,2} = B_{i,1} + \delta_3 + \delta_4 \quad \text{for $i= 1, \cdots, g$}; \\
	&a_3 = \beta_{(g+1)/2} + \delta_3; \quad a_4 = \beta_{(g+1)/2} + \delta_3 + \delta_5 + \delta_7; \\
	&b_3 = \beta_{(g+1)/2} + \delta_2; \quad b_4 = \beta_{(g+1)/2} + \delta_2 + \delta_6 + \delta_8.
\end{align*}
Put $q(\alpha_i)=1$, $q(\beta_i)=0$ for all $i$ and $q(\delta_j)=1$ for all $j$, then $q$ satisfies $(A)$ and $(B)$.
\qed

\paragraph{\it Type $\mathrm{\rmII}B\flat$ of odd $g$} \
Take the generators of $H_1(\Sigma_g^8; \Z / 2\Z)$ as in Figure~\ref{F:Spin_odd_IIBf}.
The $\Z/2\Z$ homology classes of the vanishing cycles are
\begin{align*}
	\quad &B_{0,2} = \alpha_1 + \cdots + \alpha_g + \delta_6 + \delta_7; \\
	&B_{1,2} = B_{0,2} + \beta_1 + \beta_g = \alpha_1 + \cdots + \alpha_g + \beta_1 + \beta_g + \delta_6 + \delta_7; \\
	&B_{2i,2} = B_{2i-1,2} + \alpha_i + \alpha_{g+1-i} \\
	&\qquad = \alpha_{i+1} + \cdots + \alpha_{g-i} + \beta_i + \beta_{g+1-i} + \delta_6 + \delta_7 \quad \text{for $2i=2, 4, \cdots, g-1$}; \\
	&B_{2i+1,2} =  B_{2i,2} + \beta_i + \beta_{i+1} + \beta_{g-i} + \beta_{g+1-i} \\
	&\qquad = \alpha_{i+1} + \cdots + \alpha_{g-i} + \beta_{i+1} + \beta_{g-i} + \delta_6 + \delta_7 \quad \text{for $2i+1=3, 5, \cdots, g-2$}; \\
	&B_{g,2} = \alpha_{(g+1)/2} + \beta_{(g+1)/2} + \beta_{(g+1)/2} + \delta_2 + \delta_4 + \delta_7 + \delta_8 = \alpha_{(g+1)/2} + \delta_2 + \delta_4 + \delta_7 + \delta_8; \\
	&a_1 = \beta_{(g+1)/2} + \delta_7; \quad a_2 = \beta_{(g+1)/2} + \delta_5; \\
	&b_1 = \beta_{(g+1)/2} + \delta_2 + \delta_4 + \delta_8; \quad b_2 = \beta_{(g+1)/2} + \delta_2 + \delta_4 + \delta_6; \\
	&B_{i,1} = B_{i,2} + \delta_3 + \delta_4 \quad \text{for $i= 1, \cdots, g$}; \\
	&a_3 = \beta_{(g+1)/2} + \delta_3; \quad a_4 = \beta_{(g+1)/2} + \delta_3 + \delta_5 + \delta_7; \\
	&b_3 = \beta_{(g+1)/2} + \delta_2; \quad b_4 = \beta_{(g+1)/2} + \delta_2 + \delta_6 + \delta_8.
\end{align*}
Put $q(\alpha_i)=1$, $q(\beta_i)=0$ for all $i$ and $q(\delta_j)=1$ for all $j$, then $q$ satisfies $(A)$ and $(B)$.
\qed

We summarize the diffeomorphism types of the MCK Lefschetz pencils in Tables~\ref{T:diffeotypeeven} and~\ref{T:diffeotypeodd}.
\begin{remark} \label{rmk:Kodaira-infty}
	We observe from the Tables that our pencils cover all the diffeomorphism types of minimal ruled surfaces.
	In addition, as observed in Remark~\ref{rmk:KO} the type $\mathrm{\rmI}A$ (or $\mathrm{\rmII}A$) genus-$1$ Lefschetz pencil has one more $(-1)$-section, which can be blown down to produce the minimal rational surface $\mathbb{CP}{}^{2}$.
	In this sense, we can claim that \textit{the MCK Lefschetz pencils exhaust all the diffeomorphism types of minimal symplectic $4$-manifolds with symplectic Kodaira dimension $-\infty$}.
\end{remark}
\renewcommand{\arraystretch}{1.2}
\begin{table}[htpb]
	\caption{For even $g$.} \label{T:diffeotypeeven}
	\begin{tabular}{c||cccc} \hline
		type & monodromy & spin/nonspin & total space & separating cycles\\ \hline
		$\mathrm{\rmI}A$ & $W_{\mathrm{\rmI}A}=t_{\delta_1} \cdots t_{\delta_4}$ & spin & $\Sigma_{g/2} \times S^2$ & $(g/2;1), (g/2;1)$ \\
		$\mathrm{\rmI}B$ & $W_{\mathrm{\rmI}B}=t_{\delta_1} \cdots t_{\delta_4}$ & nonspin & $\Sigma_{g/2} \tilde\times S^2$ & $(g/2;0), (g/2;2)$ \\
		$\mathrm{\rmII}A$ & $W_{\mathrm{\rmII}A}=t_{\delta_1} \cdots t_{\delta_4}$ & spin & $\Sigma_{g/2} \times S^2$ & $(g/2;1), (g/2;1)$ \\
		$\mathrm{\rmII}B$ & $W_{\mathrm{\rmII}B}=t_{\delta_1} \cdots t_{\delta_4}$ & nonspin & $\Sigma_{g/2} \tilde\times S^2$ & $(g/2;2), (g/2;2)$ \\ \hline
	\end{tabular}
\end{table}
\begin{table}[htpb]
	\caption{For odd $g$.} \label{T:diffeotypeodd}
	\begin{tabular}{c||cccc} \hline
		type & monodromy & spin/nonspin & total space \\ \hline
		$\mathrm{\rmI}A$ & $W_{\mathrm{\rmI}A}=t_{\delta_1} \cdots t_{\delta_8}$ & nonspin & $\Sigma_{(g-1)/2} \tilde\times S^2$  \\
		$\mathrm{\rmI}B\sharp$ & $W_{\mathrm{\rmI}B\sharp}=t_{\delta_1} \cdots t_{\delta_8}$ & spin & $\Sigma_{(g-1)/2} \times S^2$ \\
		$\mathrm{\rmI}B\flat$ & $W_{\mathrm{\rmI}B\flat}=t_{\delta_1} \cdots t_{\delta_8}$ & spin & $\Sigma_{(g-1)/2} \times S^2$ \\
		$\mathrm{\rmII}A$ & $W_{\mathrm{\rmII}A}=t_{\delta_1} \cdots t_{\delta_8}$ & nonspin & $\Sigma_{(g-1)/2} \tilde\times S^2$ \\
		$\mathrm{\rmII}B\sharp$ & $W_{\mathrm{\rmII}B\sharp}=t_{\delta_1} \cdots t_{\delta_8}$ & spin & $\Sigma_{(g-1)/2} \times S^2$ \\ 
		$\mathrm{\rmII}B\flat$ & $W_{\mathrm{\rmII}B\flat}=t_{\delta_1} \cdots t_{\delta_8}$ & spin & $\Sigma_{(g-1)/2} \times S^2$ \\ \hline
	\end{tabular}
\end{table}

\subsubsection{The isomorphism classes}
From the above results, it follows immediately that the type $T_A$ Lefschetz pencil is not isomorphic to the type $T_B$ Lefschetz pencil for $T_A \in \{ \mathrm{\rmI}A, \mathrm{\rmII}A \}$ and $T_B \in \{ \mathrm{\rmI}B, \mathrm{\rmII}B \}$ when $g$ is even, or $T_B \in \{ \mathrm{\rmI}B\sharp, \mathrm{\rmI}B\flat, \mathrm{\rmII}B\sharp, \mathrm{\rmII}B\flat \}$ when $g$ is odd.
When $g$ is even, we can furthermore distinguish type $\mathrm{\rmI}B$ and $\mathrm{\rmII}B$ since $W_{\mathrm{\rmI}B}$ has a separating cycle that bounds $\Sigma_{g/2}^1$ in the (holed) reference fiber while $W_{\mathrm{\rmII}B}$ does not (Hurwitz equivalence for Lefschetz pencils preserves the topological type of a vanishing cycle).
In Table~\ref{T:diffeotypeeven} we also list the topological types of separating cycles;
$(h;l)$ represents a separating cycle that separates $\Sigma_g^4$ into $\Sigma_h^l$ and $\Sigma_{g-h}^{4-l}$.
In conclusion, we showed the following result.
\begin{theorem} \label{thm:multiplepencils}
	The genus-$g$ Matsumoto-Cadavid-Korkmaz Lefschetz fibration has at least three nonisomorphic supporting minimal Lefschetz pencils when $g$ is even and at least two such pencils when $g$ is odd.
\end{theorem}

For the special case where $g=1$, we can actually show that 
\begin{enumerate}
	\item[(A)] $W_{\mathrm{\rmI}A} = t_{\delta_1} \cdots t_{\delta_8}$ and $W_{\mathrm{\rmII}A}  = t_{\delta_1} \cdots t_{\delta_8}$ are mutually Hurwitz equivalent,
	\item[(B)] $W_{\mathrm{\rmI}B\sharp} = t_{\delta_1} \cdots t_{\delta_8}$, $W_{\mathrm{\rmI}B\flat} = t_{\delta_1} \cdots t_{\delta_8}$, $W_{\mathrm{\rmII}B\sharp} = t_{\delta_1} \cdots t_{\delta_8}$ and $W_{\mathrm{\rmII}B\flat}  = t_{\delta_1} \cdots t_{\delta_8}$ are mutually Hurwitz equivalent.
\end{enumerate}
We give procedures of Hurwitz moves.
\paragraph{\it (A)} \
We start from the factorization~\eqref{eq:WIAodd} $W_{\mathrm{\rmI}A} = t_{\delta_1} \cdots t_{\delta_8}$.
\begin{align*}
	t_{\delta_1} \cdots t_{\delta_8} &= W_{\mathrm{\rmI}A} 
	= t_{B_{0,1}} t_{B_{1,1}} t_{a_1} t_{a_2} t_{b_1} t_{b_2} t_{B_{0,2}} t_{B_{1,2}} \underline{t_{a_3} t_{a_4} t_{b_3} t_{b_4}} \\
	&\sim t_{B_{0,1}} t_{B_{1,1}} t_{a_1} t_{a_2} t_{b_1} t_{b_2} \underline{t_{B_{0,2}} t_{B_{1,2}} t_{b_3} t_{b_4}} t_{a_3} t_{a_4} \\
	&\sim t_{B_{0,1}} t_{B_{1,1}} t_{a_1} t_{a_2} t_{b_1} t_{b_2} \underline{t_{b_3^{\prime}} t_{b_4^{\prime}} t_{B_{0,2}} t_{B_{1,2}}} t_{a_3} t_{a_4} \\
	&\sim t_{B_{0,1}} t_{B_{1,1}} t_{a_1} t_{a_2} t_{b_1} t_{b_2} \underline{ t_{B_{1,2}^{\prime}} t_{t_{b_3^{\prime}} t_{b_4^{\prime}}(B_{1,2})}} \; \underline{t_{b_3^{\prime}} t_{b_4^{\prime}} t_{a_3} t_{a_4}} \\
	&\sim t_{B_{0,1}} t_{B_{1,1}} t_{a_1} t_{a_2} t_{b_1} t_{b_2} t_{B_{0,2}^{\prime}} t_{B_{1,2}^{\prime}} t_{a_3} t_{a_4} t_{b_3^{\prime}} t_{b_4^{\prime}},
\end{align*}
where $b_i^{\prime} = t_{B_{0,2}} t_{B_{1,2}}(b_i)$ for $i=3,4$, $B_{1,2}^{\prime} = t_{b_3^{\prime}} t_{b_4^{\prime}} (B_{0,2})$ and $B_{0,2}^{\prime} = t_{B_{1,2}^{\prime}} t_{b_3^{\prime}} t_{b_4^{\prime}}(B_{1,2})$.
One can easily check that the last expression coincides with the factorization~\eqref{eq:WIIAodd} $W_{\mathrm{\rmII}A} = t_{\delta_1} \cdots t_{\delta_8}$ after a small readjustment of the reference fiber $\Sigma_1^8$.
\qed

\paragraph{\it (B)} \
We first start from the factorization~\eqref{eq:WIBsharpodd} $W_{\mathrm{\rmI}B\sharp} = t_{\delta_1} \cdots t_{\delta_8}$ to give the factorization~\eqref{eq:WIIBsharpodd} $W_{\mathrm{\rmII}B\sharp} = t_{\delta_1} \cdots t_{\delta_8}$. 
\begin{align*}
	t_{\delta_1} \cdots t_{\delta_8} &= W_{\mathrm{\rmI}B\sharp} 
	= t_{B_{0,1}} t_{B_{1,1}} \underline{t_{a_1} t_{a_2} t_{b_1} t_{b_2}} t_{B_{0,2}} t_{B_{1,2}} t_{a_3} t_{a_4} t_{b_3} t_{b_4}  \\
	&\sim  \underline{t_{B_{0,1}} t_{B_{1,1}} t_{b_1} t_{b_2}} t_{a_1} t_{a_2} t_{B_{0,2}} t_{B_{1,2}} t_{a_3} t_{a_4} t_{b_3} t_{b_4}\\
	&\sim \underline{t_{b_1^{\prime}} t_{b_2^{\prime}} t_{B_{0,1}} t_{B_{1,1}}} t_{a_1} t_{a_2} t_{B_{0,2}} t_{B_{1,2}} t_{a_3} t_{a_4} t_{b_3} t_{b_4} \\
	&\sim \underline{ t_{B_{1,1}^{\prime}} t_{t_{b_1^{\prime}} t_{b_2^{\prime}}(B_{1,1})}} \; \underline{t_{b_1^{\prime}} t_{b_2^{\prime}} t_{a_1} t_{a_2}} t_{B_{0,2}} t_{B_{1,2}} t_{a_3} t_{a_4} t_{b_3} t_{b_4} \\
	&\sim t_{B_{0,1}^{\prime}} t_{B_{1,1}^{\prime}} t_{a_1} t_{a_2} t_{b_1^{\prime}} t_{b_2^{\prime}} t_{B_{0,2}} t_{B_{1,2}} t_{a_3} t_{a_4} t_{b_3} t_{b_4} = W_{\mathrm{\rmII}B\sharp}, 
\end{align*}
where $b_i^{\prime} = t_{B_{0,1}} t_{B_{1,1}}(b_i)$ for $i=1,2$, $B_{1,1}^{\prime} = t_{b_1^{\prime}} t_{b_2^{\prime}} (B_{0,1})$ and $B_{0,1}^{\prime} = t_{B_{1,1}^{\prime}} t_{b_1^{\prime}} t_{b_2^{\prime}}(B_{1,1})$.
Secondly, starting from the factorization~\eqref{eq:WIBflatodd} $W_{\mathrm{\rmI}B\flat} = t_{\delta_1} \cdots t_{\delta_8}$, the exactly same procedure of Hurwitz moves as in (A) gives the factorization~\eqref{eq:WIIBflatodd} $W_{\mathrm{\rmII}B\flat} = t_{\delta_1} \cdots t_{\delta_8}$.
Finally, we start from the factorization~\eqref{eq:WIIBsharpodd} $W_{\mathrm{\rmII}B\sharp} = t_{\delta_1} \cdots t_{\delta_8}$ to give the factorization~\eqref{eq:WIIBflatodd} $W_{\mathrm{\rmII}B\flat} = t_{\delta_1} \cdots t_{\delta_8}$.
\begin{align*}
	t_{\delta_1} \cdots t_{\delta_8} &= W_{\mathrm{\rmII}B\sharp} 
	= t_{B_{0,1}} \underline{t_{B_{1,1}} t_{a_1} t_{a_2} t_{b_1} t_{b_2} t_{B_{0,2}}} t_{B_{1,2}} t_{a_3} t_{a_4} t_{b_3} t_{b_4}  \\
	&= \underline{t_{B_{0,1}}} t_{B_{0,1}^{\prime}} t_{B_{1,1}} t_{a_1} t_{a_2} t_{b_1} t_{b_2}  \underline{t_{B_{1,2}} t_{a_3} t_{a_4} t_{b_3} t_{b_4}}  \\
	&= t_{B_{0,1}^{\prime}} t_{B_{1,1}} t_{a_1} t_{a_2} t_{b_1} t_{b_2} t_{B_{0,2}^{\prime}} t_{B_{1,2}} t_{a_3} t_{a_4} t_{b_3} t_{b_4} = W_{\mathrm{\rmII}B\flat},
\end{align*}
where $B_{0,1}^{\prime} = t_{B_{1,1}} t_{a_1} t_{a_2} t_{b_1} t_{b_2} (B_{0,2})$ and $B_{0,2}^{\prime} = t_{B_{1,2}} t_{a_3} t_{a_4} t_{b_3} t_{b_4} (B_{0,1}) $.
\qed
\begin{remark} \label{Rem:g=1Steinfilling}
	In fact, the two Hurwitz inequivalent factorizations for $g=1$ have much simpler expressions discovered by the author in~\cite{Hamada_torusrel_preprint}.
	It is possible to show that (A) $W_{\mathrm{\rmI}A} = t_{\delta_1} \cdots t_{\delta_8}$ and $W_{\mathrm{\rmII}A}  = t_{\delta_1} \cdots t_{\delta_8}$ are Hurwitz equivalent to the factorization $A_8$ in~\cite{Hamada_torusrel_preprint}, and (B) $W_{\mathrm{\rmI}B\sharp} = t_{\delta_1} \cdots t_{\delta_8}$, $W_{\mathrm{\rmI}B\flat} = t_{\delta_1} \cdots t_{\delta_8}$, $W_{\mathrm{\rmII}B\sharp} = t_{\delta_1} \cdots t_{\delta_8}$ and $W_{\mathrm{\rmII}B\flat}  = t_{\delta_1} \cdots t_{\delta_8}$ are Hurwitz equivalent to the factorization $B_8$ in~\cite{Hamada_torusrel_preprint}.
	There is a very reasonable evidence that may explain why we have only two distinct factorizations for $g=1$, which is well described in~\cite{Ozbagci2015}.
	As the monodromy of an open book, the boundary multi-twist $t_{\delta_{1}} \cdots t_{\delta_{8}}$ in $\Mod(\Sigma_1^8)$ produces the contact $3$-manifold $(Y_8,\xi_8)$ that is given as the boundary of the symplectic $D^2$-bundle over $T^2$ with Euler number $-8$.
	While the symplectic $D^2$-bundle naturally gives a Stein filling of $(Y_8,\xi_8)$, there are \textit{exactly two} more Stein fillings of $(Y_8,\xi_8)$ according to Ohta and Ono~\cite{OhtaOno2003}.
	Those Stein fillings are realized as positive allowable Lefschetz fibrations over $D^2$; 
	the obvious Dehn twist factorization $t_{\delta_{1}} \cdots t_{\delta_{8}}$ gives the symplectic $D^2$-bundle, and the factorizations $A_8$ and $B_8$ provide the other two Stein fillings.
\end{remark}
\begin{remark} \label{Rem:classifyMCKLP}
	We could not distinguish (nor identify) all the types of the MCK Lefschetz pencils for general $g \geq 2$, though we believe that they are distinctive.
	To distinguish them we need a more subtle invariant of Lefschetz pencils which is not invariant for Lefschetz fibrations obtained by blowing up at base points.
\end{remark}

\section{Final Remarks}

\subsection*{Questions} \
We would like to pose some questions regarding the MCK Lefschetz pencils as well as the problem that we have just mentioned in Remark~\ref{Rem:classifyMCKLP}.
Since the constructions of the pencils produced in this paper are purely combinatorial the geometric meanings of them are not clear.
However, the high symmetricity of the vanishing cycles and base points of our pencils might suggest that some neat geometric structures are hidden in the background.
Besides, from the original geometric construction of Matsumoto's Lefschetz fibration (MCK of genus $2$) we can observe that at least a ``half" of the fibration, namely, the Lefschetz fibration over the disk corresponding to the subword $t_{B_0} t_{B_1} t_{B_2} t_{C}$, is holomorphic.
Considering those suggestive evidence, it would be reasonable to ask the following:
let $T \in \{ \mathrm{\rmI}A, \mathrm{\rmI}B, \mathrm{\rmII}A, \mathrm{\rmII}B \}$ for even $g$ or $T \in \{ \mathrm{\rmI}A, \mathrm{\rmI}B\sharp, \mathrm{\rmI}B\flat, \mathrm{\rmII}A, \mathrm{\rmII}B\sharp, \mathrm{\rmII}B\flat \}$ for odd $g$. 
\begin{question} \label{q:holomorphic}
	Is the type $T$ MCK Lefschetz pencil holomorphic?
\end{question}
A weakened version of Question~\ref{q:holomorphic} may be still interesting:
let $g$ be the genus, $n$ the number of critical points and $b$ the number of base points of a Lefschetz pencil.
\begin{question}
	Is there a holomorphic Lefschetz pencil on a ruled surface that has the same data $(g,n,b)$ as those of the type $T$ MCK Lefschetz pencil?
\end{question}
By Gompf's observation, each type of the MCK Lefschetz pencil provides a symplectic structure on the total space.
Comparing type $\mathrm{\rmI}A$ and type $\mathrm{\rmII}B$ for even $g$, we observe that the subtle difference of the locations of base points (the pictures are the same only except $\delta_2$!) matters a great deal; they changes the topology of the total spaces.
Therefore it would be interesting to investigate what kind of geometric structures, especially symplectic structures, reflects the difference of the MCK Lefschetz pencils even when the total spaces are the same.
We note, however, that the symplectic structure for an $S^2$-bundle over a Riemann surface is unique up to diffeomorphism and deformation according to~\cite{LiLiu1995}.
Therefore the different MCK Lefschetz pencils on the same ruled surface define the same symplectic structure after all.

\subsection*{Applications} \
As the MCK Lefschetz fibration has been a great source to create new Lefschetz fibrations, 
we may expect that the MCK Lefschetz pencils are a good source to create new Lefschetz \textit{pencils} ---especially, \textit{minimal} ones--- as well.
In fact, some of the pencils have already been used to produce interesting pencils, mainly by the so-called \textit{breeding} operation (see \cite{Baykur_preprint}).

In~\cite{Baykur_preprint}, Baykur made the most of our pencils.
He used $W_{\mathrm{\rmII}A}$ of genus $2$ to give a family of genus-$3$ symplectic Calabi-Yau (SCY, for short) Lefschetz pencils that contains all rational homology classes of $T^2$-bundles over $T^2$.
Independently of this work, the author and Hayano~\cite{HamadaHayano_preprint} similarly used $W_{\mathrm{\rmII}A}$ of genus $2$ to realize all homeomorphism classes of $T^2$-bundles over $T^2$ admitting sections, as genus-$3$ SCY Lefschetz pencils.
They also showed that the pencil homeomorphic to the four-torus $T^4$ is in fact diffeomorphic to $T^4$.
They further generalized the construction to gain a genus-$g$ SCY Lefschetz pencil for arbitrary $g\geq3$ whose total space is homeomorphic (diffeomorphic when $g$ is odd) to $T^4$.
As another generalization, Baykur~\cite{Baykur_preprint} utilized $W_{\mathrm{II}A}$ of even genus to construct a genus-$g$ Lefschetz fibrations(pencils) for any odd $g\geq3$ with $b_1=g+1$, which is the largest among the known genus-$g$ Lefschetz fibrations(pencils).
Even non-maximal lifts are also useful.
Capping off boundary components, Baykur~\cite{Baykur_preprint} did breeding using $W_{\mathrm{II}A} = t_{\delta_1} t_{\delta_2}$ (with $\delta_3$ and $\delta_4$ capped off) and $W_{\mathrm{II}A} = t_{\delta_1}$ (with $\delta_2, \delta_3, \delta_4$ capped off) of genus $2$, along with a genus-$2$ smallest Lefschetz pencil, to construct genus-$3$ pencils on exotic $\mathbb{CP}{}^{2} \# p\overline{\mathbb{CP}}{}^{2}$ for $p=7,8,9$.
By varying this construction, he also found infinitely many genus-$3$ noncomplex Lefschetz pencils.
One more application can be found in~\cite{BaykurHamada_preprint}, where $W_{\mathrm{\rmII}A} = t_{\delta_1} t_{\delta_2} t_{\delta_4}$ (with $\delta_3$ capped off) of genus $2$ was utilized to find a set of maximal (which is three) disjoint $(-1)$-sections of the genus-$2$ Lefschetz fibration discovered in~\cite{BaykurKorkmaz_preprint} as one having the smallest possible number (which is seven) of critical points.

As illustrated in those examples, the MCK Lefschetz pencils appear well-suited for the breeding operation.
Besides the above, we can find various other configurations in the relations we have constructed, with which one may conveniently carry out breedings.
We expect further applications will be followed using them.


\appendix


\noindent
{\it Acknowledgements.} 
The author would like to thank Susumu Hirose, Naoyuki Monden and Ryoma Kobayashi for the fruitful discussions and hospitality during his visit to Tokyo University of Science in July 2013, during which he came up with the very first idea of this paper.
The author also thanks Kenta Hayano for informing him of Theorem~\ref{Thm:BHM}.




\begin{thebibliography}{99999}

\bibitem{AkhmedovOzbagci2014}
A.~Akhmedov and B.~Ozbagci, 
\textit{Singularity links with exotic Stein fillings}, 
J.\ Singul.\ \textbf{8} (2014), 39-–49.

\bibitem{AkhmedovSaglam2015}
A.~Akhmedov and K.~N.~Saglam,
\textit{New exotic $4$-manifolds via Luttinger surgery on Lefschetz fibrations}, 
Internat.\ J.\ Math.\ \textbf{26} (2015), 1550010, 21 pp. 

\bibitem{AkhmedovSakalli2016}
A.~Akhmedov and S.~Sakall\i,
\textit{On the geography of simply connected nonspin symplectic $4$-manifolds with nonnegative signature}, 
Topology Appl.\ \textbf{206} (2016), 24-–45. 

\bibitem{Baykur2012}
R.~\.{I}.~Baykur, 
\textit{Non-holomorphic surface bundles and Lefschetz fibrations}, 
Math.\ Res.\ Lett.\ \textbf{19} (2012), 567-–574. 

\bibitem{Baykur2016}
R.~\.{I}.~Baykur, 
\textit{Minimality and fiber sum decompositions of Lefschetz fibrations}, 
Proc.\ Amer.\ Math.\ Soc.\ \textbf{144} (2016), 2275-–2284.

\bibitem{Baykur_preprint}
R.~\.{I}.~Baykur, 
\textit{Small symplectic Calabi-Yau surfaces and exotic $4$-manifolds via genus-$3$ pencils},
preprint; \verb|https://arxiv.org/abs/1511.05951|.

\bibitem{BaykurKorkmaz_preprint}
R.~\.{I}.~Baykur and M.~Korkmaz, 
\textit{Small Lefschetz fibrations and exotic $4$-manifolds},
Math.\ Ann.\ (2016), doi:10.1007/s00208-016-1466-2.

\bibitem{BaykurHamada_preprint}
R.~\.{I}.~Baykur and N.~Hamada, 
\textit{Exotic rational surfaces via genus-$5$ pencils},
in preparation.

\bibitem{BaykurHayano2016}
R.~\.{I}.~Baykur and K.~Hayano, 
\textit{Multisections of Lefschetz fibrations and topology of symplectic $4$-manifolds}
Geom.\ Topol.\ \textbf{20} (2016), 2335--2395. 

\bibitem{BaykurHayano_preprint}
R.~\.{I}.~Baykur and K.~Hayano, 
\textit{Hurwitz equivalence for Lefschetz fibrations and their multisections}, 
Real and Complex Singularities, 1--24, Contemp.\ Math., \textbf{675}, Amer.\ Math.\ Soc., Providence, RI, 2016.

\bibitem{BaykurHayanoMonden_preprint}
R.~\.{I}.~Baykur, K.~Hayano and N.~Monden, 
\textit{Unchaining surgery and symplectic $4$-manifolds}, 
in preparation.

\bibitem{Cadavid1998}
C.~Cadavid, 
\textit{On a remarkable set of words in the mapping class group},
Thesis~(Ph.D.), The University of Texas at Austin, 1998.

\bibitem{Donaldson1999}
S.~K.~Donaldson, 
\textit{Lefschetz pencils on symplectic manifolds},
J.\ Differential Geom.\ \textbf{53} (1999), 205-–236. 

\bibitem{EKKOS2002}
H.~Endo, M.~Korkmaz, D.~Kotschick, B.~Ozbagci and A.~Stipsicz, 
\textit{Commutators, Lefschetz fibrations and the signatures of surface bundles}, 
Topology \textbf{41} (2002), 961--977.

\bibitem{FarbMargalit2012}
B.~Farb and D.~Margalit, 
\textit{A primer on mapping class groups}, 
Princeton Math.\ Ser., Vol.\ 49, Princeton Univ.\ Press, Princeton, NJ, 2012. 

\bibitem{GompfStipsicz1999}
R.~E.~Gompf and A.~I.~Stipsicz, 
\textit{$4$-manifolds and Kirby calculus}, 
Graduate Studies in Math., Vol.\ 20, Amer.\ Math.\ Soc., Providence, RI, 1999.

\bibitem{Gurtas_preprint2005}
Y.~Z.~Gurtas, 
\textit{Positive Dehn twist expressions for some elements of finite order in the mapping class group},
preprint; \verb|https://arxiv.org/abs/math/0501385|.

\bibitem{Hamada_torusrel_preprint}
N.~Hamada,
\textit{Simple expressions for the holed torus relations},
preprint; \\ \verb|https://arxiv.org/abs/1701.02171|.

\bibitem{HamadaHayano_preprint}
N.~Hamada and K.~Hayano, 
\textit{Topology of holomorphic Lefschetz pencils on the four-torus},
preprint; \verb|https://arxiv.org/abs/1603.08284|.

\bibitem{HamadaKobayashiMonden_preprint}
N.~Hamada, R.~Kobayashi and N.~Monden, 
\textit{Non-holomorphic Lefschetz fibrations with $(-1)$-sections}, 
preprint; \verb|https://arxiv.org/abs/1609.02420|.


\bibitem{Kobayashi2016}
R.~Kobayashi,
\textit{On genera of Lefschetz fibrations and finitely presented groups}, 
Osaka J.\ Math.\ \textbf{53} (2016), 351–-376. 

\bibitem{KobayashiMonden2016}
R.~Kobayashi and N.~Monden, 
\textit{Lefschetz pencils and finitely presented groups}, 
Pacific J.\ Math.\ \textbf{282} (2016), 359–-388. 

\bibitem{Korkmaz2001}
M.~Korkmaz, 
\textit{Noncomplex smooth $4$-manifolds with Lefschetz fibrations},
Internat.\ Math.\ Res.\ Notices\ \textbf{2001} (2001), 115--128.

\bibitem{Korkmaz2009}
M.~Korkmaz, 
\textit{Lefschetz fibrations and an invariant of finitely presented groups},
Int.\ Math.\ Res.\ Not.\ IMRN\ \textbf{2009} (2009), 1547--1572. 


\bibitem{KorkmazOzbagci2008}
M.~Korkmaz and B.~Ozbagci, 
\textit{On sections of elliptic fibrations}, 
Michigan Math.\ J.\ \textbf{56} (2008), 77--87.


\bibitem{LiLiu1995}
T.-J.~Li and A.~Liu,
\textit{Symplectic structure on ruled surfaces and a generalized adjunction formula},
Math.\ Res.\ Lett.\ \textbf{2} (1995), 453–-471. 

\bibitem{Matsumoto1996}
Y.~Matsumoto, 
\textit{Lefschetz fibrations of genus two -- a topological approach},
Topology and Teichm\"{u}ller spaces, Proceedings of the 37th Taniguchi Symposium, pp.~123--148, 
World Scientific, River Edge, NJ, 1996.	


\bibitem{OhtaOno2003}
H.~Ohta and K.~Ono, 
\textit{Symplectic fillings of the link of simple elliptic singularities},
J.\ Reine Angew.\ Math.\ \textbf{565} (2003), 183–-205. 

\bibitem{OkudaTakamura_preprint}
T.~Okuda and S.~Takamura, 
\textit{Sequences of degenerations of propeller surfaces and their splittings}, 
in preparation.

\bibitem{Onaran2010}
S.~\c{C}.~Onaran, 
\textit{On sections of genus two Lefschetz fibrations},
Pacific J.\ Math.\ \textbf{248} (2010), 203–-216. 


\bibitem{Ozbagci2015}
B.~Ozbagci, 
\textit{On the topology of fillings of contact $3$-manifolds},
Geometry \& Topology Monographs Vol.\ \textbf{19} (2015), 73--123.

\bibitem{OzbagciStipsicz2000}
B.~Ozbagci and A.~I.~Stipsicz, 
\textit{Noncomplex smooth $4$-manifolds with genus-$2$ Lefschetz fibrations},
Proc.\ Amer.\ Math.\ Soc.\ \textbf{128} (2000), 3125–-3128. 

\bibitem{OzbagciStipsicz2004}
B.~Ozbagci and A.~I.~Stipsicz, 
\textit{Contact $3$-manifolds with infinitely many Stein fillings},
Proc.\ Amer.\ Math.\ Soc.\ \textbf{132} (2004), 1549-–1558.

\bibitem{OzbagciStipsicz2004book}
B.~Ozbagci and A.~I.~Stipsicz, 
\textit{Surgery on contact $3$-manifolds and Stein surfaces},
Bolyai Soc.\ Math.\ Stud., 13.\ Springer-Verlag, Berlin, 2004.

\bibitem{ParkYun2009}
J.~Park and K.-H.~Yun, 
\textit{Nonisomorphic Lefschetz fibrations on knot surgery $4$-manifolds},
Math.\ Ann.\ \textbf{345} (2009), 581-–597. 

\bibitem{ParkYun2011}
J.~Park and K.-H.~Yun, 
\textit{Lefschetz fibration structures on knot surgery $4$-manifolds},
Michigan Math.\ J.\ \textbf{60} (2011), 525-–544.



\bibitem{Stipsicz2001}
A.~I.~Stipsicz,
\textit{Spin structures on Lefschetz fibrations},
Bull.\ London Math.\ Soc.\ \textbf{33} (2001), 466–-472. 

\bibitem{Stipsicz2002}
A.~I.~Stipsicz,
\textit{Surface bundles with nonvanishing signature},
Acta Math.\ Hungar.\ \textbf{95} (2002), 299-–307. 

\bibitem{Tanaka2012}
S.~Tanaka, 
\textit{On sections of hyperelliptic Lefschetz fibrations},
Algebr.\ Geom.\ Topol.\ \textbf{12} (2012), 2259--2286.

\bibitem{Yun2006}
K.-H.~Yun, 
\textit{On the signature of a Lefschetz fibration coming from an involution}, 
Topology Appl.\ \textbf{153} (2006), 1994-–2012. 







\end{thebibliography}
\end{document}